\pdfoutput = 1
\documentclass[11pt]{amsart}
\usepackage{graphicx}
\usepackage{amssymb,amsmath,latexsym} 
\newcommand    {\comment}[1]   {{}}

\newcommand{\labell}[1] {\label{#1}}

\newcommand{\ts} {\textstyle}
\renewcommand{\Tilde}{\widetilde}
\renewcommand{\Hat}{\widehat}
\newcommand{\QED}{\hfill$\Box$}
\newcommand{\less}{{\smallsetminus}}
\newcommand{\NI}{{\noindent}}
\newcommand{\SSS}{{\smallskip}}
\newcommand{\MS}{{\medskip}}

\newcommand{\la}{{\lambda}}
\newcommand{\De}{{\Delta}}

\newcommand{\Ga}{{\Gamma}}
\newcommand{\ga}{{\gamma}}

\newcommand{\PU}{{PU}}

\newcommand\fk {{\mathfrak k}}
\newcommand\ft {{\mathfrak t}}
\newcommand\fg {{\mathfrak g}}
\newcommand\p {{\partial}}

\newcommand{\Uu}{{\mathcal U}}

\newcommand{\Hh}{{\mathcal H}}
\newcommand{\Aa}{{\mathcal A}}

\newcommand{\Cc}{{\mathcal C}}
\newcommand{\Oo}{{\mathcal O}}
\newcommand{\Ii}{{\mathcal I}}

\newcommand{\La}{{\Lambda}}
\newcommand{\CP}{{\mathbb CP}}
\newcommand{\Proj}{{\mathbb P}}
\newcommand{\C}{{\mathbb C}}
\newcommand{\Q}{{\mathbb Q}}
\newcommand{\R}{{\mathbb R}}
\newcommand{\Z}{{\mathbb Z}}

\newcommand{\om}{{\omega}}
\newcommand{\ka}{{\kappa}}
\newcommand{\al}{{\alpha}}

\newcommand{\io}{{\iota}}
\newcommand{\be}{{\beta}}
\newcommand{\Symp}{{\rm Symp}}

\newcommand{\Ham}{{\rm Ham}}

\newcommand{\SU}{{SU}}
\newcommand{\U}{{U}}

\newcommand{\Aut}{{\rm Aut}}
\newcommand{\Isom}{{\rm Isom}}

\newcommand{\eps}{{\epsilon}}
\newcommand{\id}{{\rm id}}

\newcommand{\Po}{{{\mathcal P}_{\om}}}
\newcommand{\ssminus}{{{\smallsetminus}}}

\newtheorem{theorem}{Theorem}[section]
\newtheorem{thm}[theorem]{Theorem}
\newtheorem{corollary}[theorem]{Corollary}
\newtheorem{cor}[theorem]{Corollary}

\newtheorem{lemma}[theorem]{Lemma}
\newtheorem{proposition}[theorem]{Proposition}
\newtheorem{prop}[theorem]{Proposition}
\newtheorem{definition}[theorem]{Definition}
\newtheorem{defn}[theorem]{Definition}
\newtheorem{rmk}[theorem]{Remark}
\newtheorem{example}[theorem]{Example}
\numberwithin{figure}{section}
\numberwithin{equation}{section}
\begin{document}   
    
\title{Polytopes with mass linear functions, part I}
\author{Dusa McDuff}\thanks{First author partially supported by NSF grant DMS 0604769, and second by NSF grant DMS 0707122.}
\address{Department of Mathematics,
Barnard College, Columbia University, 3990 Broadway, New York, NY 110027}
\email{dusa@math.columbia.edu}

\author{Susan Tolman}
\address{Department of Mathematics,
 University of Illinois at Urbana--Champaign, Urbana, IL  61801,
USA}
\email{tolman@math.uiuc.edu}
%\urladdr{http://www.}
\date{June 25, 2008,  revised August 10, 2009}

\begin{abstract} 

Let $\De$ be a  $n$-dimensional polytope that is simple,  that is, exactly $n$ facets meet at each vertex.  
An affine function  is ``mass linear''
on $\De$ if its value on the center of mass of $\De$
depends linearly on the positions of the supporting  hyperplanes.
On the one hand, we show that 
certain types of symmetries of $\Delta$ give rise
to nonconstant mass linear functions on $\Delta$.
On the other hand, 
we show that most polytopes do not admit 
any nonconstant mass linear functions. 
Further,
if every affine function is mass linear on
$\Delta$, then $\Delta$ is 
a product of
simplices.
Our main result is a classification 
of all 
smooth polytopes of dimension $\leq 3$
which admit nonconstant mass linear functions.
In particular there is only 
one family of smooth  $3$-dimensional polytopes -- and no polygons --
that admit  `` essential mass linear functions'',  that is,
mass linear functions that do not arise from the symmetries  described above.
In part II, we  will complete 
this classification in 
the $4$-dimensional case.

These results have geometric implications.
Fix  a  symplectic toric manifold $(M,\omega,T,\Phi)$
with moment polytope $\Delta = \Phi(M)$.  
Let $\Symp_0(M,\omega)$ denote 
the identity component of the group of symplectomorphisms of  $(M,\om)$.
Any 
linear function $H$ on $\Delta$
generates a 
Hamiltonian $\R$ action on $M$ 
whose closure is a subtorus $T_H$ of $T$.
We show that  if  the map   
$\pi_1(T_H)\to \pi_1(\Symp_0(M,\omega))$ 
has finite image, then $H$ is mass linear.   
Combining this fact and the claims described above, we prove that
in most cases the induced map 
$\pi_1(T) \to \pi_1(\Symp_0(M,\omega))$ 
is an injection. 
Moreover, 
the map does not have finite image unless $M$ is 
a product of projective spaces. 
Note also that  there is a natural maximal 
compact   connected subgroup 
$\Isom_0(M) \subset \Symp_0(M,\om)$;
there is a natural compatible complex structure $J$ on $M$, and 
$\Isom_0(M)$ is the identity component of the
group  of
symplectomorphisms that also preserve this structure. 
We prove that 
if the polytope $\De$  supports no
essential mass linear functions, then the induced map 
$\pi_1(\Isom_0(M)) \to \pi_1(\Symp_0(M,\omega))$ 
is injective.
Therefore, this map is injective for all $4$-dimensional symplectic toric 
manifolds and is
injective in the $6$-dimensional case 
unless $M$ 
is a $\CP^2$ bundle over $\CP^1$.
\end{abstract}

\maketitle

\begin{center}
with Appendix written with V. Timorin\\ \MS
\end{center}   

\tableofcontents

%%%%%%%%%%%%%%%%%%%%%%%%%%%%%%%%%%%%%%%%%%%%%%%%%%%%%%%%%%%%%%
\section{Introduction}
%%%%%%%%%%%%%%%%%%%%%%%%%%%%%%%%%%%%%%%%%%%%%%%%%%%%%%%%%%%%%%

To begin, we shall give some basic definitions; we then state the main results. These are elaborated in \S\ref{ss:add}; their geometric implications are 
described in \S\ref{ss:geom}.

Let $\ft$ be a real vector space,
let $\ft^*$ denote the dual space,  and 
let $\langle \, \cdot \, , \, \cdot \, \rangle \colon \ft \times \ft^* \to \R$ 
denote the natural pairing.
Let $A$ be an affine space modeled on $\ft^*$.
A {\bf (convex) polytope} $\Delta \subset A$ 
is the bounded intersection of a finite set of affine half-spaces.
Hence, $\Delta$ can  be written
\begin{equation}\labell{eq:De}
\De =   
\bigcap_{i = 1}^N \{ x \in A \mid h_i(x) \leq \kappa_i  \},
\end{equation}
where the $h_i$ are affine functions on $A$ and the {\bf support numbers}
$\kappa_i$ lie in $\R$
for all $1 \leq i \leq N$;
the {\bf outward conormals}  
are the unique\footnote
{
The ambiguity in choosing $h_i$ is discussed in Remarks \ref{rmk:choices} and \ref{rmk:normalize}.}
vectors $\eta_i \in \ft$ so that $h_i(y) - h_i(x) =  \langle \eta_i, y-x \rangle$
for all $x$ and $y$ in $A$.   When the $\{h_i\}$ are fixed and the dependence on 
$\{\ka_i\}$ is important we write $\De: = \De(\ka)$.

In  this paper, we will always assume that
$A$ is the  affine span of $\Delta$
and that the affine  span of each facet  $F_i := 
\Delta \cap \{ x \in A \mid h_i( x)  = \kappa_i  \}$ 
is a hyperplane.
In particular, each 
$F_i$ is nonempty.
Moreover, the polytopes we consider are {\bf simple},
that is, $\dim \ft$ facets
meet at every vertex. 

We  sometimes need stronger assumptions:
Given an integer lattice $\ell \subset \ft$, a
polytope $\Delta \subset A$ is {\bf rational}
if we can choose each outward conormal to lie in  
$\ell$.
A rational polytope is {\bf smooth} 
(or Delzant)
if, for each vertex of $\De$,
the primitive outward conormals 
to the facets which meet at that vertex form a basis for 
$\ell$. 
Here, a vector $\eta \in \ell$ is {\bf primitive}
if it is not a positive integer multiple of any other lattice element.

The {\bf chamber} $\Cc_\De$ of $\Delta: = \De(\ka)$
is the connected component that 
contains $\kappa$
of the set of all $\kappa' \in \R^N$ such that $\De(\ka')$ is simple.
Note that  $\Cc_\De$ is open.
Moreover, for every $\kappa' \in \Cc_\De$ the polytope $\Delta(\kappa')$
is {\bf analogous}\footnote
{
Two arbitrary polytopes are {\bf analogous} if there is  combinatorial equivalence
 (cf. Definition \ref{def:combeq})  between them in which  corresponding facets are parallel.  In the case at hand, these conditions reduce to the ones mentioned above.}
to $\Delta$, that is, 
it has the property that
the intersection 
$\bigcap_{i \in I} F'_i$ is empty exactly if the intersection 
$\bigcap_{i \in I} F_i$ is empty for all $I \subset \{1,\ldots,N\}$.
Since $\Delta$ is simple, $\Cc_\De$ is a nonempty 
connected
open set.

\begin{figure}[htbp] 
   \centering
\includegraphics[width=1.8 in]{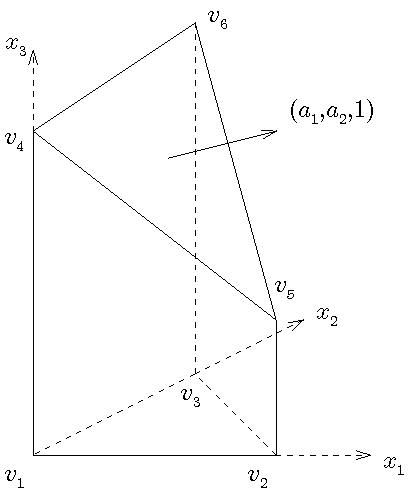} 
   \caption{The polytope $Y_a$ constructed in Example~\ref{ex:1}.}
   \label{fig:0}
\end{figure}

\begin{example}\labell{ex:1} \rm
Let $\{e_i\}_{i=1}^n$ denote the standard basis for $\R^n$.
Given $a = (a_1,a_2) \in \R^2$,
define
\begin{gather*}
\eta_1 = -e_1, \  \eta_2 = -e_2, \  \eta_3 = e_1 + e_2, \ 
\eta_4 = - e_3, \ \ \eta_5 = e_3 + a_1 e_1 + a_2 e_2, \ \  \mbox{and} \\
\Cc_a = \left\{ \kappa \in \R^5 \left| \
 \sum_{i=1}^3 \kappa_i > 0 \mbox{ and }
\kappa_4 + \kappa_5 >  - a_1 \kappa_1  - a_2 \kappa_2   
+  \max ( 0,  a_1 ,a_2) \sum_{i=1}^3 \kappa_i \right. \right\}.
\end{gather*}
Given $\kappa \in \Cc_a$, let 
\vspace{-.1in}
$$Y = Y_a(\kappa) =
\bigcap_{i = 1}^5 \{ x \in (\R^3)^* \mid \langle \eta_i, x\rangle 
\leq \kappa_i  \};$$
it is illustrated in Figure~\ref{fig:0}. 
One can simplify the formulas here by translating $Y$ so that 
its first vertex is
at the origin.  Then $\ka_1=\ka_2=\ka_4=0$,
and we see that 
$Y_a(\ka)$ really depends 
only on two parameters, $\ka_3$ and $\ka_5$.
The vertices now are
$$
 (0,0,0),\;\; (\ka_3,0,0), \;\;(0,\ka_3,0),\;\; (0,0,\ka_5),\;\;
 (\ka_3,0,\ka_5-a_1\ka_3),\;\;(0,\ka_3,\ka_5-a_2\ka_3).
 $$
The polytope  $Y_a(\kappa)$ is simple;  it is smooth
exactly if $a \in \Z^2$; and
its chamber  is $\Cc_a$.
\end{example}

Let $c_\De \colon  \Cc_\De  \to A$ 
denote the {\bf center of mass} function, which assigns 
to $\ka : = (\kappa_1,\dots,\ka_N)\in \Cc_\De$
the center of mass of $\De(\ka)$ 
(with respect to any affine identification of $A$
with  Euclidean space).

\begin{defn}\labell{def:kalin}  
Let $\Delta \subset A$ be a simple polytope.
Let $c_\De \colon  \Cc_\De  
\to A$ 
denote the center of mass function.
An affine function $H \colon A \to \R$ is called
 {\bf mass linear} with respect to $\Delta$ if the composite
%an affine function $H \colon A \to \R$ so that
$$
H \circ c_{\Delta}  \colon \Cc_\De \to \R 
$$ 
is  linear, that is,
$$ H( c_{\Delta}(\kappa))=  \sum_{i=1}^N   \gamma_i \kappa_i + C \quad \forall \, \kappa \in \Cc_\Delta, $$
where $C \in \R$ and $\gamma_i \in \R$ for all $1 \leq i \leq N$.
For short, we shall say that  $H$ is  a {\bf mass linear function} on $\De$;
 we sometimes also call  $(\De,H)$  a {\bf  mass linear pair}.

We call $\gamma_i$ the {\bf coefficient} of the support number $\kappa_i$
(in $H \circ c_\Delta$).
If, additionally,  $\gamma_i \in \Z$ for all $1 \leq i \leq N$, then we say that
$H$ is a {\bf mass linear function with integer coefficients} on $\De$.
\end{defn}

In practice, we usually choose an 
identification of $A$ with the linear space  $\ft^*$
and consider 
those
$H$ which are elements 
of the dual space $\ft$;  
cf. Remarks~\ref{rmk:affine}~and~\ref{rmk:choices}.

\comment{The referee objected to "mass linear function on $\Delta$".
In some ways, I see the point.  
It is a function on all of $\ft^*$, not just $\Delta$.
On the other hand, it isn't really mass linear "on" $\Delta$,
either. It is math linear *with respect to* $\Delta$.
Anyway, my main point is that we use this phrase *many* times  
(including the definition above).  Do we want to change all of them?
Just change the one the referee pointed out?}

Our primary goal in this paper is to analyze   
mass linear functions on simple polytopes.
On the one hand, we show that certain
types of symmetries of $\Delta$ give rise
to nonconstant  
 functions that are mass linear on $\Delta$.
We call nonconstant  mass linear functions 
that arise in this way 
{\bf inessential}  (see
Definition~\ref{def:iness});
the others are called
{\bf essential.}
On the other hand, we show that   
most polytopes do not admit any nonconstant 
mass linear functions.
Further,
if every affine function is mass linear
on $\Delta$, then $\Delta$ is affine
equivalent to a product of simplices.
Our main result is a classification 
of all  $2$-dimensional simple polytopes and  $3$-dimensional
smooth polytopes  which admit nonconstant mass linear functions.
In particular, the polytopes
constructed in Example~\ref{ex:1} are the only
smooth
polytopes 
of dimension $\le 3$ 
which admit essential  mass linear functions.
In part II, we  will complete this classification in the $4$-dimensional case.

Our first result holds both for an arbitrary polygon and for a smooth polygon, though, as we explain in Remark \ref{rmk:affine} below, the words 
used have different interpretations in the two cases. This remark also gives the definition of affine equivalence.

\begin{thm}\labell{thm:2} 
A simple $2$-dimensional polytope supports  a nonconstant 
mass linear function
exactly  if it is  
a triangle, a trapezoid or a parallelogram. All such functions  are inessential. 
The same statement holds in the smooth case.
\end{thm}

\begin{theorem}\labell{thm:3dim}
Let $\Delta$ be a smooth $3$-dimensional polytope.
If there exists an essential mass linear function on
$\Delta$, then $\Delta$ is affine equivalent to 
the polytope $Y_a(\kappa)$ constructed in Example~\ref{ex:1} for some $a \in \Z^2$
and $\kappa \in \Cc_Y$.
Further, for any  $a \in \Z^2$ 
there exists an essential mass linear function on the polytope $Y_a$
exactly if  
$a_1\neq 0, a_2\neq 0$ and $a_1\neq a_2$.
\end{theorem}

We prove the second statement in Theorem \ref{thm:3dim}
by  a direct calculation; see  
Proposition~\ref{prop:Mab}.
 However the first claim 
 relies on many of our other results.

These results have geometric
implications.  
Fix a  symplectic toric manifold $(M,\omega,T,\Phi)$ with moment polytope  $\Delta = \Phi(M)$.  (See \S\ref{ss:geom}.)
Let $\Symp_0(M,\omega)$ denote the identity component of the 
group of symplectomorphisms of
$(M,\omega)$.  
Any linear function $H$ on $\Delta$
generates a 
Hamiltonian $\R$ action on $M$ 
whose closure is a subtorus $T_H$ of $T$.
We show that  if  the map  $\pi_1(T_H)\to 
\pi_1(\Symp_0(M,\omega))$ 
has finite image, then $H$ is mass linear.   
Since (as we explained above) most polytopes do not admit mass linear functions,
this implies that
in most cases the induced map  
$\pi_1(T) \to \pi_1(\Symp_0(M,\omega))$ 
is an injection.
Moreover,  
since  products of simplices are the only polytopes for which every affine function
is mass linear,
the map does not have finite image unless $M$ is 
a product of projective spaces.

Note also that  there is a natural maximal compact
connected
subgroup $\Isom_0(M) \subset \Symp_0(M,\om)$;
there is a natural compatible complex structure $J$ 
on $M$ and 
hence an  associated metric $g_J: = \om(\cdot,J\cdot)$, and
we define $\Isom_0(M)$ to be the identity component of the
group  of
isometries of the Riemannian manifold $(M,g_J)$. 
It is well known that
$\Isom_0(M) \subset \Symp_0(M,\om)$.  
Theorem~\ref{thm:geo} below states that
if $\Delta$ admits no essential
mass linear functions, then the natural map
from $\pi_1(\Isom_0(M))$ to $\pi_1(\Symp_0(M,\omega))$
is an injection. 
Since the polytopes constructed in Example~\ref{ex:1}
are exactly the moment polytopes of $\CP^2$ bundles over $\CP^1$ (see Example~\ref{ex:ZK}),
Theorems~\ref{thm:2} and \ref{thm:3dim} have the following 
consequence.

\begin{theorem}\labell{thm:23dimg}
Let $(M,\omega,T,\Phi)$ be a symplectic toric manifold,
and let $\Isom_0(M)$ 
be
the identity component of the associated K\"ahler
isometry group.
If $M$ is $4$-dimensional, or if $M$ is $6$-dimensional but
is not a  $\CP^2$ bundle over $\CP^1$,
then the natural map from $\pi_1(\Isom_0(M))$ to $\pi_1(\Symp_0(M,\omega))$
is an injection.
\end{theorem}

We end this subsection with some 
important
technical remarks.

\begin{rmk} \rm \labell{rmk:affine}
The group of affine transformations of $A$ acts 
on simple polytopes $\Delta \subset A$ and 
affine functions $H \colon A \to \R$ by
$a \cdot (\De, H) = (a(\De), H \circ a^{-1})$.
In this case, we will say that $(\Delta, H)$ and $(a(\Delta), H \circ a^{-1})$
are {\bf affine equivalent}.
When we restrict to smooth polytopes 
we will only allow
affine transformations $a \colon A \to A$ such that
the {\bf associated linear transformation} $\Hat{a} \colon \ft^* \to \ft^*$ 
induces an automorphism of the 
integer
 lattice;
here,
$\Hat{a}$ is the unique linear map  so
that $\Hat{a}(x - y) = a(x) - a(y)$ for all $x$ and $y$ in $\ft^*$.

Except for the concepts of ``volume'' and ``moment''
considered in \S\ref{ss:volume} below,
all the notions  in this paper are affine invariant.
For example, 
since $a(c_\Delta) = c_{a(\Delta)}$,  
an affine function
$H$ is mass linear 
on $\Delta$ exactly if $ H \circ a^{-1}$ 
is mass linear on $a(\De)$.
Therefore,  
by identifying $A$ with $\ft^*$, we may restrict our
attention to polytopes $\Delta \subset \ft^*$.
Given two   affine equivalent polytopes $\De$ and $\De'$, we will 
usually simply say that they are the same 
and write $\De=\De'$.
For example, 
whether we work with simple or smooth polytopes, 
there is a unique equivalence class of $k$-dimensional polytopes
with $k+1$ facets, 
namely the $k$-simplex $\De_k$. (See Example \ref{ex:simplex}.)  
Hence we say that every such polytope {\em is} $\De_k$.
These conventions mean that, for example,
the interpretation of the word \lq\lq is" in Theorem~\ref{thm:2}
depends on whether we are considering smooth polytopes or simple ones.
\end{rmk}

There is another much weaker notion that is sometimes useful.

\begin{defn}\labell{def:combeq}
 Two polytopes $\Delta$ and $\Delta'$ are
said to be {\bf combinatorially equivalent} 
  if there exists a bijection of facets
$F_i \leftrightarrow F'_i$
so that $\bigcap_{i \in I} F_i \neq \emptyset$ exactly if
$\bigcap_{i \in I} F'_i \neq \emptyset$ for all $I \subset \{1,\ldots,N\}$.
\end{defn}

We will never use this equivalence relation without explicitly saying so.
Note 
that corresponding  facets of combinatorially equivalent polytopes need not be parallel.  
For example, any two polygons with the same number of edges are combinatorially equivalent.

\begin{rmk}\rm \labell{rmk:choices}
Each outward conormal of a simple polytope $\Delta \subset A$
is only well defined up to multiplication
by a  positive constant.
(See Remark~\ref{rmk:normalize} for an important normalization convention.)
The chamber $\Cc_\Delta$ of $\Delta$ depends on these
choices but the notion of 
mass linearity  
does not -- it is
well defined.
On the other hand,
if  $H$ is mass linear the coefficients of the support numbers depend
on the choice of the outward conormals.

Even once these outward conormals are determined, the
affine functions $h_i$ 
of Equation (\ref{eq:De})
are only well defined up to addition by  a
constant. However, these choices have no significant consequence.
In particular, if $\Delta \subset \ft^*$ we may assume that
the $h_i$ are homogeneous.
Similarly,  we may
restrict our attention to homogeneous linear functions $H \in \ft$.
 \end{rmk}

%%%%%%%%%%%%%%%%%%%%%%%%%%%%%
\subsection{Additional   results on polytopes}\labell{ss:add}
%%%%%%%%%%%%%%%%%%%%%%%%%%%%%%%%%

Let us begin with a few definitions.

\begin{defn}\labell{def:symmetric}
Fix an affine function $H \colon A \to \R$ and a simple polytope 
$\Delta\subset A$. 
We say that the facet $F$ is 
{\bf symmetric} (or {\bf $\mathbf H$-symmetric}\footnote{When there is no danger of confusion, we will omit the $H$.})
if $\langle H, c_\Delta \rangle \colon \Cc_\Delta \to \R$
does not depend on the support number of $F$. Otherwise we say that $F$ is 
{\bf asymmetric}  (or {\bf $\mathbf H$-asymmetric}).
More generally, we say that a face  $f$ is {\bf symmetric} 
if it is the intersection of symmetric  facets.
\end{defn}

\begin{defn}\labell{def:perv}  
We say that a 
facet $F$ of a simple  
polytope $\Delta\subset A$ is  
{\bf pervasive} if it meets (i.e. has nonempty intersection with)
every other facet of $\Delta$. 
We say that $F$ is {\bf flat} if there is a hyperplane in $\ft$ that 
contains the conormal of every 
facet (other than $F$ itself)
that meets $F$.
\end{defn}

For example, in the polytopes constructed
in Example~\ref{ex:1}, the triangular facets are flat but the quadrilaterals
are not;
the pervasive facets are the quadrilaterals.

In \S\ref{s:main}, we analyze the key features of symmetric and asymmetric
facets.
On the one hand, we prove
that if $H$ is mass linear on  $\Delta$, then $H$ is also mass
linear on each symmetric face $f$ of $\Delta$.
This allows us to analyze mass linear functions on polytopes which
have symmetric faces by ``induction'' on the dimension of the polytope.  
On the other hand, we prove that if  $H$ is a mass linear function
on $\Delta$, 
then asymmetric facets are very special; they
must be either pervasive  or flat.
As we show in \S\ref{ss:asym}, this implies
that most polytopes do not support
any nonconstant mass linear functions.

\begin{thm}\labell{thm:easy} If $\Delta \subset A$  
is a simple polytope that 
contains no pervasive facets and no flat facets then every mass linear 
function
on $\Delta$ is constant. 
\end{thm}

On the other hand, there are some polytopes with 
nonconstant 
mass linear functions.  
The simplest type of such functions -- 
inessential functions --   arise
from symmetries of the polytope.

\begin{defn}\labell{def:robust}
Let $\Delta \subset A$ 
be a simple polytope.  A  {\bf symmetry} of $\Delta$  
is an  affine map $a \colon A \to A$ so
that $a(\Delta) = \Delta$.
When $\Delta$ is a smooth polytope 
we will only allow
affine transformations $a \colon A \to A$ such that
the associated linear transformation  
induces an automorphism of the 
integer lattice; cf. Remark~\ref{rmk:affine}.
A {\bf robust symmetry} of $\Delta$
is a linear map $\Hat{a} \colon \ft^* \to \ft^*$
so that for all $\kappa \in \Cc_\Delta$ there exists
an affine map  $a_\kappa \colon A \to A$
such that
\begin{enumerate}
\item $a_\kappa$ is a symmetry of $\Delta(\kappa)$, and
\item $\Hat{a}$ is the linear map associated to $a_\kappa$.
\end{enumerate}
Finally, we say that two facets $F_i$ and $F_j$ are 
{\bf equivalent}, denoted $F_i \sim F_j$, if there exists
a robust symmetry $\Hat{a}$ so that each $a_{\kappa}$ takes 
$F_i(\ka)$ to $F_j(\ka)$
\end{defn}

It is clear that this defines an equivalence relationship
because the robust symmetries form a subgroup of the
group of linear transformations. 
As an example, consider the square 
$$
\Delta = \{x \in (\R^2)^* \mid -1 \leq x_1 \leq 1 \mbox{ and } -1 \leq x_2 \leq  1\}.
$$
The
map $(x_1,x_2) \mapsto  (x_2, -x_1)$ is a symmetry
of $\Delta$
which is not robust.
In contrast,
the  map $(x_1,x_2) \mapsto (-x_1,-x_2)$ is a robust symmetry of $\Delta$.
Therefore, each edge is equivalent to the opposite 
edge.

In \S\ref{ss:ines} we give several  easily verified  criteria 
for facets to be equivalent.
In Proposition~\ref{prop:preaut}, 
we show that in the smooth case
 equivalent
facets can also be characterized 
in terms of the symmetries of the associated toric manifold $M$.
Moreover, in this case, $F_i\sim F_j$ exactly if the corresponding divisors in $M$ are homologous; cf. Remark~\ref{rmk:Cox}.

\begin{rmk}\rm\label{rmk:normalize}
Let $\Delta \subset A$ be a simple polytope.
We can (and always will) choose the outward conormals 
so that for every every robust symmetry
$\Hat{a} \colon \ft^* \to \ft^*$ of $\Delta$ and every facet $F_i$, 
there exists a facet $F_j$ so that $\Hat{a}^*(\eta_i) = \eta_j$.
Since the robust symmetries form a group, the fact that this is possible
follows easily
from Lemma~\ref{le:symmetries}.
If $\Delta \subset \ft^*$ is a smooth polytope, then  
by definition $\Hat{a}^*$ takes each primitive outward conormal to another primitive outward
conormal. Hence, when we restrict to smooth polytopes we can (and will) always choose
the unique primitive outward conormal to each facet.
\end{rmk}

\begin{definition}\labell{def:iness}
Let $\Delta  =
\bigcap_{i = 1}^N \{ x \in \ft^* \mid \langle \eta_i, x\rangle \leq \kappa_i  \}$
be a simple polytope; 
let $\Ii$ denote the set of equivalence
classes of facets.
We say that 
 $H \in \ft$ 
is {\bf inessential} iff
\begin{equation}\labell{eq:iness}
H = \sum \be_i \eta_i,
 \quad \mbox{where} \quad
\be_i\in \R
\ \forall \, i \quad \mbox{and}
\quad \sum_{i \in I} \beta_i = 0
\quad \forall \, I \in \Ii. 
\end{equation}
Otherwise, we say that $H$ is {\bf essential}.
\end{definition}

If each equivalence class in $\Ii$ has just one element then 
every inessential 
function on $\De$ vanishes.
More generally, the set of inessential functions 
is an $(N - |\Ii|)$-dimensional subspace of $\ft$;  
see  Lemma~\ref{le:rational2}.  

\begin{example}\rm\labell{ex:simplex}
As an example, (here and subsequently)  
let $\Delta_k$ denote the standard $k$-simplex
$$ 
\Delta_k = \{ x \in (\R^{k})^* \mid x_i \geq 0 \ \forall \, i 
 \mbox{ and } \sum x_i \leq 1 \};
$$
the outward conormals are $\{\eta_i\}_{i=1}^{k+1}$,
where $\eta_i = -e_i$ for all $i \in \{1,\ldots,k\}$
and $\eta_{k+1} = \sum_{i=1}^k e_i$.
The linear map 
$(x_1,\ldots,x_k) \mapsto (-\sum x_i, x_1,\dots,x_{k-1}) $
is a robust symmetry of $\Delta$
whose transpose permutes the conormals $\eta_i$ 
cyclically.
Hence, all  $k+1$ facets of $\De_k$ are equivalent, and every $H \in \R^k$
is inessential. 
For example, $e_j = -\eta_j + \frac{1}{k+1} \sum_{i=1}^{k+1} \eta_i$
for all 
$j\le k$.
Moreover,  a direct computation shows that 
$$
\langle e_j , c_{\Delta_k}(\kappa) \rangle = -\kappa_j + {\ts \frac{1}{k+1}}
\sum_{i=1}^{k+1} \kappa_i.
$$
\end{example}

\begin{example}\labell{ex:square}\rm
Consider a product of simplices $\Delta = \Delta_{{k}}
\times \Delta_{k'} \subset (\R^{{k}})^* \times
(\R^{k'}) = (\R^{{k}+k'})^*.$
This polytope has ${k} + k' + 2$ facets.
The first  ${k}+1$ facets are of the form $F_i \times \Delta_{k'}$, and
are all equivalent to each other.
The remaining $k' +1$ facets are also equivalent, and are of the form
$\Delta_{{k}} \times F_i$.
Again, every $H \in \R^{{k} + k'}$ is inessential and
mass linear;
\begin{equation*}
\langle e_j, c_\Delta(\kappa) \rangle =
\begin{cases}
-\kappa_j +  \frac{1}{{k}+1} \sum_{i=1}^{{k}+1} \kappa_i & 
\mbox{if } j \leq {k} \\
-\kappa_{j+1} +  \frac{1}{k'+1} \sum_{i=1}^{ k' +1} \kappa_{k+1+i} & 
\mbox{otherwise}
\end{cases}
\end{equation*}
A nearly identical remark applies if $\Delta$ is the product
of three or more simplices.

Here, the {\bf product} of simple polytopes
$\Hat\Delta \subset \Hat A$ and $\Tilde \Delta \subset  \Tilde A$
is $\Delta = \Hat \Delta \times \Tilde \Delta \subset
\Hat A \times \Tilde A$.  This is a simple polytope;
it is smooth exactly if $\Hat \Delta$ and $\Tilde \Delta$
are smooth.
\end{example}

\begin{example}\labell{ex:1b}\rm
The triangular facets of the polytope $Y = Y_a(\kappa)$
defined in Example~\ref{ex:1}
are equivalent because they are exchanged by the affine
reflection $a_\kappa \colon (\R^3)^* \to (\R^3)^*$ defined by
$$
a_\kappa(x_1,x_2,x_3) = (x_1,x_2, - a_1 x_1 - a_2 x_2 -  x_3 + \kappa_5 - \kappa_4).
$$
Hence, 
$H =  \eta_4 - \eta_5 =  -a_1 e_1 - a_2 e_2 -2 e_3$ is inessential.
\end{example}

Although the definition of an inessential function looks somewhat opaque
at first glance, this is one of our most important notions. As we shall see in
Corollary \ref{cor:aut},
in the smooth case
 inessential functions have
 a very natural  interpretation in terms of the geometry
  of the associated toric manifold.   Formula (\ref{eq:iness}) is also closely tied to
 algebraic properties of robust symmetries; cf. Remark~\ref{rmk:p}.

We now show that
every inessential function is mass 
linear; cf. Example~\ref{ex:simplex}.

\begin{proposition}\labell{prop:inessential}
Let $\Delta \subset \ft^*$
be a simple polytope;
let $\Ii$ denote the set of equivalence
classes of facets.
If 
$H \in \ft$ is inessential,
write  $H =  \sum \be_i \eta_i$, where the $\be_i$ satisfy the conditions of 
Equation (\ref{eq:iness}). Then 
$$
\langle H, c_{\De}(\ka) \rangle = \sum \be_i \kappa_i.
$$
\end{proposition}

\begin{proof}
Assume that $F_i \sim F_j$.
By definition, for all $\kappa \in \Cc_\Delta$
there exists an affine transformation $a_\kappa \colon \ft^* \to \ft^*$
so that $a_\kappa$ takes $\Delta(\kappa)$ to itself 
and takes $F_i$ to $F_j$.
Since $a_\kappa(F_i) = F_j$,
$$ 
\langle \eta_i, x \rangle - \kappa_i = \langle \eta_j, a_\kappa(x)
\rangle - \kappa_j \quad \forall \, x \in \ft^*;
$$
see Remark~\ref{rmk:normalize}.
On the other hand, since $a_\kappa(\Delta(\kappa)) = \Delta(\kappa)$, $a_\kappa$
fixes the center of mass $c_{\Delta}(\kappa)$.
Hence, 
$$\langle \eta_i - \eta_j, c_{\Delta}(\kappa) \rangle  = \kappa_i - \kappa_j.$$
Finally, every inessential $H$ is a linear combination
of terms of the form $\eta_i - \eta_j$, where $F_i \sim F_j$. 
\end{proof}

Note that, in general, even if $H = \sum_i \beta_i \eta_i$
is mass linear the function $\langle H, c_\Delta(\kappa) \rangle$ need not
equal $\sum \beta_i \kappa_i$; indeed, the  coefficients of
the $\eta_i$  are not
uniquely determined by $H$ while the coefficients of $\ka_i$ are. 
In contrast,
as we shall see in 
Lemma~\ref{le:Hsum}, if $H$ is mass linear and 
$\langle H, c_\Delta(\kappa) \rangle=\sum \al_i\ka_i$ then 
$H=\sum \al_i\eta_i$.

In \S\ref{s:ines}, we analyze polytopes with inessential
functions.
In particular, we show that there are
exactly two types of polytopes which  admit inessential functions:
bundles over the $k$-simplex (see Definition~\ref{def:bund})
and $k$-fold expansions (see Definition~\ref{def:expand}).
We then show that it is possible to
use  Proposition~\ref{prop:inessential} to 
reduce the number of asymmetric
facets that we need to consider;
we can  subtract an inessential function $H'$ from a mass linear 
function $H$ to find a new mass linear function $\Tilde{H} = H - H'$ with
fewer asymmetric facets.
For example,
if $H$ is a mass linear function 
on a simple polytope $\Delta$ then,
after possibly  subtracting an inessential function, 
we may assume that every asymmetric facet is pervasive
(Proposition~\ref{prop:flat}).
This has the following consequence. 

\begin{prop}\label{prop:mid} If $\Delta$ is a simple 
polytope that 
contains no pervasive facets then every mass linear 
function on $\Delta$ is inessential.
\end{prop}

In Example~\ref{ex:square}, we saw that every affine function
on a product of simplices is mass linear.
Our next theorem shows that these are the only
polytopes with this property.
The equivalence of conditions (iii) and (iv) below is reminiscent of 
Nill's Proposition~2.18 in \cite{NILL}.  
However  because Nill works only with rational polytopes the results 
are somewhat different.

\begin{thm} \labell{thm:allmass} Let $\De$ be a  simple (or smooth)
polytope.  The following conditions are equivalent:
\begin{enumerate}
\item [(i)]  every $H \in \ft$ is mass linear on $\Delta$;
\item [(ii)]   every $H \in \ft$ is inessential  on $\Delta$;
\item [(iii)]   $\De$ is a product of simplices;
\item [(iv)]  
$\bigcap_{i \in I} F_i = \emptyset$ for every
equivalence class of facets 
$I$; in particular, $|I|>1$.
\end{enumerate}
\end{thm}

In \S\ref{s:23dim},
we analyze mass linear functions on low dimensional polytopes.
In particular, we prove
Theorems~\ref{thm:2} and~\ref{thm:3dim}.

In part II, we will give a complete list $4$-dimensional
polytopes which admit essential mass linear functions.
As described above, the main technique is to use ``induction''
by looking at lower dimensional symmetric faces.
In the Appendix, 
in joint work with Timorin,
we prove the key result that will
enable us to begin this process.

\begin{thm} \labell{thm:4dim}
Let $H \in \ft$ be a mass linear function on a smooth
$4$-dimensional polytope $\Delta \subset \ft^*$. 
Then there exists an inessential function $H' \in \ft$
so that the mass linear function $\Tilde H = H - H'$
has the following property: at least one facet of
$\Delta$ is $\Tilde H$-symmetric.
\end{thm}

To prove this, we show that  every $4$-dimensional
simple polytope which admits a  mass linear function with
no symmetric facets is combinatorially
equivalent to the product of simplices.

%%%%%%%%%%%%%%%%%%%%%%%%%%%%%%%%%%%%%%%%%%%%%%%%%%%%%%%%%%%%%%
\subsection{Geometric motivation}\labell{ss:geom}
%%%%%%%%%%%%%%%%%%%%%%%%%%%%%%%%%%%%%%%%%%%%%%%%%%%%%%%%%%%%%%

Our
motivation for studying 
mass linear functions on polytopes
is geometric; we wish to
understand the fundamental group of the group 
of symplectomorphisms of a symplectic toric manifold. 
For other approaches to this question see Januszkiewicz--K\c edra~\cite{JK}, K\c edra--McDuff~\cite{KM}, McDuff--Tolman~\cite{MT1}, and Vi\~ na~\cite{Vi}. 
Additionally,  Entov--Polterovich~\cite{EP} give some interesting 
applications of compressible  
subtori of the group of symplectomorphisms; see Definition~\ref{def:comp}.

Let $(M,\om,T, \Phi)$ 
be a symplectic toric manifold,
where $M$ is a compact connected 
manifold
of dimension $2n$,
$\om$ is a symplectic form on $M$,
$T = \ft/\ell$ is a compact $n$-dimensional torus,
and $\Phi \colon   M \to \ft^*$ is a moment map for an effective 
$T$ action on $(M,\om)$.
Then the {\bf moment polytope}  $\Delta = \Phi(M) \subset \ft^*$ is a smooth polytope. 
Conversely, every smooth polytope 
$\Delta \subset \ft^*$ determines a canonical
symplectic  toric manifold $(M_\Delta, \om_\De,  \Phi_\Delta)$
with moment polytope $\Delta = \Phi_\Delta(M_\Delta)$.
(See \S\ref{s:geometry} for more details.)
Finally, any two symplectic toric manifolds with the
same moment polytope are equivariantly symplectomorphic \cite{Del}.

For example, the symplectic  
manifold associated to 
the standard  $n$-simplex
$\Delta_n$
is  complex projective space $\CP^n$ with
the standard (i.e. diagonal) $(S^1)^n$ action and suitably normalized
Fubini-Study form.
More generally, a product of simplices corresponds
to a product of projective spaces.
Similarly, the polytopes
constructed in Example~\ref{ex:1} correspond to $\CP^2$ bundles over $\CP^1$;
cf. Example~\ref{ex:ZK}.

Let $\Symp_0(M,\omega)$ denote the identity component of
the group of symplectomorphisms of $(M,\omega)$.
Since the torus $T$ is a subgroup of $\Symp_0(M,\omega)$, each $H \in \ell$
induces a circle $\Lambda_H \subset
T \subset 
 \Symp_0(M,\omega)$.
More generally, each $H \in \ft$ generates a Hamiltonian $\R$-action
on $M$ whose closure is a subtorus 
$T_H \subset T \subset \Symp_0(M,\omega)$.
The next proposition  is proved in \S\ref{s:geometry}
by interpreting the quantity 
$\langle H, c_\Delta(\kappa) \rangle$ as the value of 
Weinstein's action homomorphism
$\Aa_H: \pi_1(\Symp_0(M_\Delta,\omega_\Delta))\to \R/\Po$.

\begin{prop}\labell{prop:symp} 
Let $(M,\omega,T,\Phi)$ be a symplectic toric manifold
with moment polytope $\Delta  \subset \ft^*$;
fix $H \in \ell$.
If  the image of $\Lambda_H$
in $\pi_1(\Symp_0(M,\omega))$ is trivial,  then
$H$ is a mass linear function with integer coefficients, that is,
there exist integers $\alpha_1,\dots,\alpha_N$ so that
$$
\langle H, c_\Delta(\kappa) \rangle  =  
\sum \al_i \kappa_i
$$
for all $\ka$ in
the chamber $\Cc_\Delta$.
\end{prop}

\begin{defn}
\labell{def:comp}
Given  a topological group $G$, a 
subtorus   $K \subset G$ is
{\bf compressible} in $G$ if the 
natural map $\pi_1(K) \to \pi_1(G)$ 
has finite image.
\end{defn}

 \begin{cor}\labell{cor:symp} 
Let $(M,\omega,T, \Phi)$ be a symplectic toric manifold
with moment polytope 
$\Delta  \subset \ft^*$. 
If $K \subset T$ is compressible in $\Symp_0(M,\omega)$, 
then every $H \in \fk \subset \ft$ is mass linear.
 \end{cor}

\begin{thm}\labell{thm:geo}
Let $(M,\omega,T,\Phi)$ be a symplectic toric manifold with moment
polytope $\Delta$.
\begin{itemize}
\item [(i)]
The map $\pi_1(T) \to \pi_1(\Symp_0(M,\omega))$ is an injection if
there are no mass linear functions $H \in \ell$ on $\Delta$.

\item [(ii)]
The torus $T$ is compressible in $\Symp_0(M,\om)$ exactly if
$(M,\om)$ is a product of projective spaces with a product symplectic form.
In particular, the image of $\pi_1(T)$ in  
$\pi_1\bigl(\Symp_0(M,\om)\bigr)$ is never zero.
\end{itemize}
\end{thm}

\begin{proof}
Part (i) is an immediate consequence of Proposition \ref{prop:symp}.
The first claim of part  (ii) follows by combining 
Theorem~\ref{thm:allmass} with Corollary~\ref{cor:symp}.
The second claim follows from Proposition~\ref{prop:symp}
and Example~\ref{ex:square}.
\end{proof}

\begin{rmk}\labell{rmk:nomlin}\rm 
More generally, since
\lq\lq most" polytopes have no
nonconstant  mass linear functions, the above theorem implies that
in most cases the map
 $\pi_1(T) \to \pi_1(\Symp_0(M,\omega))$ is injective.  For example,
 in the $4$-dimensional case Theorem \ref{thm:2} implies that it is injective 
if $\De$ has five or more facets,  while in 
the  $6$-dimensional case 
Theorem~\ref{thm:3d1}
implies that it is injective 
unless $M$ is a  bundle over a projective space or is a smooth pencil 
of $4$-dimensional toric manifolds  that intersect along a $\C P^1$. 
(Here, we are using Remark~\ref{rmk:geobun} and 
Remark~\ref{rmk:geoexpand}.)
Further,  
for each facet $F_i$, let $x_i\in H^2(M;\Z)$ denote
the Poincar\'e dual to the divisor $\Phi^{-1}(F_i)$. 
Then Theorem \ref{thm:easy} implies the following: the map is injective 
if, for each $i$,
$x_i^2\ne 0$ but there is $j\ne i$ such that $x_i \cdot x_j = 0.$ 
Thus we can always make the map injective by equivariantly blowing $M$ up at a suitable 
number of fixed points. 
\end{rmk}

Next, we claim that, given any smooth polytope $\Delta \subset \ft^*$,
there is also a canonical  
complex structure $J_\De$ on the associated symplectic
toric manifold $M_\De$;  this complex structure
is $T$-equivariant and is compatible with $\om_\Delta$.
By the discussion above, this implies that there is a 
compatible complex structure $J$ on every symplectic
toric manifold $(M,\omega,T,\Phi)$; moreover, up
to $T$-equivariant symplectomorphism, this complex
structure is canonical. 
For example, the canonical complex structure on
$\CP^n$ is the standard one.

If there are robust symmetries of
the moment polytope $\Delta = \Phi(M)$,
then 
there are extra symmetries of 
$(M,\om)$
which preserve the canonical K\"ahler structure 
$(\om,J, g_J)$, where $g_J$ denotes 
the associated metric $\om(\cdot,J\cdot)$.
We shall see in Proposition~\ref{prop:preaut}  that these robust symmetries
 in some sense generate $\Isom_0(M)$, the identity component of the 
group of isometries of the Riemannian manifold $(M,g_J)$.
Moreover,  $\Isom_0(M)$
is a maximal
connected compact Lie subgroup of $\Symp_0(M,\omega)$.
For example, if $M= \CP^k$, $\Isom_0(M)$ 
is the projective unitary group $\PU(k+1)$. 
Finally, note that $T$ is also a subgroup of $\Isom_0(M)$.

\begin{lemma}\labell{le:aut}   
Let $(M,\omega,T,\Phi)$ be a symplectic toric manifold,
and let $\Isom_0(M)$ be the identity component of the associated K\"ahler isometry group.
The map $\pi_1(T) \to \pi_1(\Isom_0(M))$
is surjective.

Moreover, 
let $\Delta 
= \bigcap_{i = 1}^N 
\{ x \in \ft^* \mid \langle \eta_i,x\rangle \leq \kappa_i  \}$
be the moment polytope.
Let $\Ii$ denote the set of equivalence classes
of facets  of $\Delta$,
and fix 
$H \in \ell$. The image of $\Lambda_H$ in 
$\pi_1(\Isom_0(M))$ is trivial
exactly if $H$  can be written
$$
H =  \sum \be_i \eta_i, \quad \mbox{where}  \ 
\be_i \in \Z \ \forall\, i \,  \mbox{and} \ \sum_{i \in I}  \be_i =  
0
\ \forall \, I \in \Ii.
$$
\end{lemma}

The above lemma motivated the definition of inessential function. Its proof,  given in 
\S\ref{s:geometry}, is independent of the 
intervening sections, except for the discussion of the equivalence relation in \S\ref{ss:ines}.

\begin{cor}\labell{cor:aut} 
Let $(M,\omega,T,\Phi)$ be a symplectic toric manifold,
and let $\Isom_0(M)$ be the identity component of the associated K\"ahler isometry group.
Given $H \in \ft$,
the torus $T_H$ is compressible in $\Isom_0(M)$
exactly if $H$ is inessential on $\Delta$. 
\end{cor}

\begin{proof} 
The Lie algebra $\ft_H$ of $T_H$ is the smallest rational
subspace of $\ft$ containing 
$H$.
Since the set of inessential $H' \in \ft$  is clearly
a rational subspace of $\ft$, this implies that
$H$ is mass linear exactly if every 
$H' \in \ell\cap \ft_H$ is mass linear. 
Therefore
it suffices to prove the result for $H\in \ell$.  But    
by Lemma~\ref{le:rational2},  
 $H \in \ell$ is inessential
exactly if 
$H$ can be written
$
H =  \sum \be_i \eta_i,$ where
$\be_i \in \Q $ for all $i$ and $\sum_{i \in I}  \be_i =  
0 $ for all $I \in \Ii.$
By Lemma~\ref{le:aut}, this implies that $H$ is inessential
exactly if there exists a natural
number $m$ such that $\La_{mH}$ contracts in
$\Isom_0(M)$, that is, exactly if $\La_H$ has finite
order in $\Isom_0(M)$.
\end{proof}

Our final result concerns the map
$$
\rho:   \pi_1(\Isom_0(M))\to \pi_1(\Symp_0(M,\omega))
$$
induced by inclusion.

\begin{prop}\labell{prop:symp2}
Let $(M,\omega,T,\Phi)$ be a symplectic toric manifold,
and let $\Isom_0(M)$  be the identity component of 
the associated K\"ahler isometry group.
If there are no essential mass linear functions $H \in \ell$
on the moment polytope $\Delta \subset \ft^*$, 
then the
map  $\rho: \pi_1(\Isom_0(M))\to \pi_1(\Symp_0(M,\omega))$
is an injection.
\end{prop}

\begin{proof}
Assume 
that $\rho$ is not injective.
By Lemma~\ref{le:aut}, the map 
$\pi_1(T) \to \pi_1(\Isom_0(M))$  
is surjective.  
Therefore,  
there exists
$H \in \ell$ such that
the image of $\Lambda_H$ in $\pi_1(\Symp_0(M,\omega))$ vanishes
but the image of $\Lambda_H$ in $\pi_1(\Isom_0(M))$ 
does not.
 By Proposition~\ref{prop:symp}, 
$H$ is mass linear; moreover,
$\langle H, c_\Delta \rangle =
\sum \alpha_i \kappa_i$, where $\alpha_i \in \Z$ for all $i$.

Since we have assumed that $\De$ has no essential mass linear functions, 
$H$ must be inessential.
Hence, 
 we may write $H = \sum\be_i\eta_i$ 
 where $\sum_{i \in I} \be_i  = 0$
for each equivalence class of facets $I$.
By Proposition~\ref{prop:inessential} we must have
$\langle H, c_\Delta(\kappa) \rangle = \sum\be_i\ka_i$.   
But then $\be_i = \al_i$ because the variables $\ka_i$ are
independent.
In particular $\be_i\in \Z$.
 Lemma~\ref{le:aut} then implies that
the image of $\Lambda_H$ in $\pi_1(\Isom_0(M))$  vanishes.
This gives a contradiction.
\end{proof}

\begin{example}\labell{ex:cpnsymp}\rm
If $M = \CP^k$, then 
$H = e_1 $ 
generates a circle $\La_H$
with order $k+1$ in 
$\pi_1(\Isom_0(\C P^k)) = \pi_1(PU(k+1))=\Z/(k+1)$.  
It follows from Proposition \ref{prop:symp} 
(and Example~\ref{ex:simplex})
that $\La_H$ also has order
$k+1$ in $\pi_1(\Symp_0(\C P^k))$,  a result first proved in
Seidel~\cite{Sei}; see also \cite{MT1}.
\end{example}

\begin{rmk}\rm  
The converse of Proposition~\ref{prop:symp} 
is open, that is,
we do not know if the following statement holds:
\begin{quote}{\it 
If $H \in \ell$ is a mass linear function 
with integer coefficients, 
then
the image of $\Lambda_H$
in $\pi_1(\Symp_0(M_\Delta,\omega_\Delta))$ is trivial.}
\end{quote}
This   seems very likely to be 
the 
case, at least 
in dimension $3$,
and is the subject of
ongoing research.

If the statement above did hold, it would
imply 
that
the converse to part (i) of Theorem~\ref{thm:geo} is also true,
and  hence
that
\begin{quote}
{\it The map  $\rho: \pi_1(\Isom_0(M))\to \pi_1(\Symp_0(M,\omega))$
is an injection exactly if 
there are no essential mass linear functions $H \in \ell$ on $\Delta$.}
\end{quote}
To see this, note that if there were an essential  mass linear $H \in \ell$
then by  Remark~\ref{rmk:rational1} it would have rational coefficients.
Hence some  multiple would have 
integer
coefficients. 
By our assumption and by Lemma~\ref{le:aut}, 
this would  give an element of $\pi_1(T)$ which lies in the 
kernel of  the homomorphism $\pi_1(T)\to \pi_1(\Symp_0(M,\om))$ 
but does not lie in the kernel of the homomorphism $\pi_1(T) \to \pi_1( \Isom_0(M))$.
\end{rmk}

\begin{example}\rm 
Let $M$ be the one point blow up 
of $\C P^3$.
  Then
$\Isom_0(M) = \U(3)$, and Theorem~\ref{thm:23dimg} 
implies that $\pi_1(\U(3)) = \Z$ injects into $\pi_1(\Symp_0(M,\om))$,
 a result already noted by Vi\~na~\cite{Vi}.  In this case, we also understand the behavior of the higher homotopy groups, at least rationally;  $\pi_3(\U(3))=\Z$ injects into $\pi_3(\Symp_0(M,\om))$ by Reznikov~\cite{Rez} (cf. also ~\cite{KM}), while
 $\pi_5(\U(3))\otimes \Q$ injects into $\pi_5(\Symp_0(M,\om))\otimes \Q$
 because it has nontrivial image in $\pi_5(\CP^3)\otimes \Q$ under the
 composite of the point evaluation map $\Symp_0(M,\om)\to M$ with the blow down $M\to \C P^3$.
 \end{example} 

\begin{rmk}\rm   
If $M$ is  a  $\CP^2$ bundle over $\CP^1$,
 the map  $\pi_1(T)\to \pi_1(\Isom_0(M))$ 
is not injective: 
the action of $\SU(2)$ on the base $\C P^1$ lifts
to an action of $\SU(2)$ on $M$  that contains one of the circle factors in $T$.
On the other hand, for generic toric manifolds
$\Isom_0(M) = T$. 
 One might wonder whether the condition
$\Isom_0(M) = T$ is enough to imply that 
$\pi_1(T)$ injects into $\pi_1(\Symp_0(M,\omega))$. 
 However, this seems unlikely to be true since, 
as we shall see in Part II,
there are mass linear pairs $(H,\Delta)$ 
where $H$ is not constant and 
$\De$ has no robust symmetries, i.e. is such that
$\Isom_0(M) = T$.  
In this case, although $\La_H$ may represent a nonzero element in
$\pi_1(\Symp_0(M,\omega))$,  our methods are not 
sufficiently powerful to detect this.
 \end{rmk}

\begin{rmk}\labell{rmk:Nill}\rm
There seems to be very little work 
on simple polytopes that mentions the center of mass.
 A notable exception is
Nill's paper~\cite{NILL} on complete toric varieties.  
Interestingly enough, 
this paper is also concerned with questions 
about the symmetries of $M$, but
considers the full group of biholomorphisms $\Aut(M)$ 
 rather than $\Isom_0(M)$,  
which is one of its 
maximal connected compact subgroups.  
His work applies in the special case of smooth Fano varieties.
The  moment polytope $\De$ 
of such a variety
contains a special point $p_\De$.  
Geometers usually normalize the polytope so 
that $p_\De$ lies at the origin but this point can be described purely in terms of the moment polytope;  
see Entov--Polterovich \cite[Prop.~1.6]{EP}. 
Nill gives a direct combinatorial proof 
of the fact  that if a Fano moment polytope $\De$ has the property that 
$p_\Delta$ is equal to the center of mass $c_\De$ of $\De$,
then $\Aut(M)$ is reductive.
(For terminology, see Remark~\ref{rmk:Cox} below.)
In fact, 
the difference 
between $p_\De$ and $c_\De$
is precisely the Futaki invariant, which is now known to vanish
if and only if $(M,J)$ 
supports a K\"ahler--Einstein metric; cf. Wang--Zhu \cite{WaZ}.
Moreover the automorphism group of any K\"ahler--Einstein manifold
is  reductive; this gives an alternate proof 
of Nill's result.  

However, the {\it existence} of mass linear functions, 
which is the question that concerns us here, is not
related in any simple way to the vanishing of the Futaki invariant.
For example, the Futaki invariant of the blow up of $\C P^2$ at three points vanishes.  But the corresponding 
polytope supports no nontrivial mass linear functions
by Theorem~\ref{thm:2} 
(or by Theorem \ref{thm:easy}).
Thus $\pi_1(T^2)$ injects into $\pi_1(\Symp_0(M,\om))$  
in this case. 
\end{rmk}

\NI {\bf Organization of the paper.}
The main results in this paper are Theorems \ref{thm:2} and \ref{thm:3dim} that classify the 
low dimensional cases, Theorem \ref{thm:allmass} concerning the structure of polytopes for which all $H$ are mass linear,
and Theorems~\ref{thm:23dimg} and \ref{thm:geo} which give geometric interpretations of these results.
A more elaborate version of Theorem \ref{thm:2} is stated and proved as Proposition \ref{prop:2dim}.
Theorem \ref{thm:3dim} is proved in \S\ref{ss:3dim}, and leads into
Proposition \ref{prop:3d} which  classifies  all $3$-dimensional mass linear pairs. 
Theorem \ref{thm:allmass}  is proved in \S\ref{ss:allmass}, while all the geometric results are proved in \S\ref{s:geometry}.
\MS

\NI {\bf Acknowledgments.}  We warmly thank Vladen Timorin
for his important contribution to the Appendix.  His ideas 
also helped us clarify some of the proofs in \S\ref{s:main}.
We are also very grateful to the referee for a very careful reading 
and for
the many suggestions that helped us improve the exposition.

%%%%%%%%%%%%%%%%%%%%%%%%%%%%%%%%%%%%%%%%%%%%%%%%%%%%%%%%%%%%%%
\section{Basic ideas}\labell{s:main}
%%%%%%%%%%%%%%%%%%%%%%%%%%%%%%%%%%%%%%%%%%%%%%%%%%%%%%%%%%%%%%

In this section, we introduce the basic ideas which underlie 
the results in this paper and its sequel.
More specifically, we 
first introduce the volume and moment functions and then use them to 
analyze the key features of
symmetric faces and  asymmetric facets.
Additionally, we prove Theorem~\ref{thm:easy} in \S\ref{ss:asym}.

%%%%%%%%%%%%%%%%%%%%%%%%%%%%%%%%%%%%%%%%%%%%%%%%%%%%%%%%%%%%%%
\subsection{Notational conventions}\labell{ss:def}
%%%%%%%%%%%%%%%%%%%%%%%%%%%%%%%%%%%%%%%%%%%%%%%%%%%%%%%%%%%%%%

Unless there is explicit mention to the contrary, we 
fix an identification of the affine space $A$ with $\ft^*$ and hence
assume throughout
that 
\begin{equation}\labell{eq:De1}
\Delta =
\bigcap_{i = 1}^N \{ x \in \ft^* \mid \langle \eta_i,x\rangle \leq \kappa_i  \}
\end{equation}
is a simple polytope with facets 
$F_1,\ldots,F_N$
and outward conormals $\eta_1,\dots,\eta_N$.    
Further, for each $I\subset \{1,\dots,N\}$, we denote the intersection of the corresponding facets by
\begin{equation}\labell{eq:FI}
F_I: = \bigcap_{i\in I} F_i.
\end{equation}

We will need
to understand the faces
of a polytope as polytopes in their own right.
To this end, fix  a simple polytope $\Delta$.
Given a  (nonempty) face 
$f = F_I$, 
let $P(f) \subset \ft^*$ 
denote the smallest
affine subspace which contains $f$.
Since $\Delta$ is simple,  $P(f)$  has codimension $|I|$; indeed, 
\begin{equation}\labell{eq:Pf}
P(f) =
\bigcap_{i \in I} \{ x \in \ft^* \mid \langle \eta_i,x\rangle = \kappa_i  \}.
\end{equation}
Additionally,
 if  $F$ is a facet of $\Delta$
which does not contain $f$,
then $F \cap f$ is either empty or is a facet of $f$;
conversely, every facet of $f$ has this form.  
Therefore,
$$ 
f =
\bigcap_{j \in I_f}
\{ x \in P(f) \mid \langle \eta_j,x\rangle \leq \kappa_j  \},
$$
where  $j \in  \{1,\ldots,N\}$ lies in
$I_f$ exactly if
$j \not\in I$ and $F_j \cap f \neq \emptyset$.

Let  $\Hat \ft$ denote the quotient of $\ft$ by the
span of the $\{\eta_{i}\}_{i\in I}$, and
let $\pi \colon \ft \to \Hat \ft$ be the natural projection.
Fix $x \in P(f)$ and define an isomorphism  $j_x \colon
\Hat \ft^* \to P(f)$ by $j_x(y) = \pi^*(y) + x$.
Then
\begin{equation}\labell{eq:identifying}
f_x := 
j_x^{-1}(f) = \bigcap_{j \in I_f}
 \{ y \in \Hat\ft^* \mid 
\langle \pi(\eta_j), y \rangle \leq \kappa_j - \langle \eta_j, x \rangle \};
\end{equation}
in particular,
$f_x \subset \Hat \ft^*$ is a convex polytope with
outward conormals $\{ \pi(\eta_j)\}_{j \in I_f}.$
With this explicit description of the faces of $\De$, it is 
then easy to establish the following important facts:
\begin{itemize}
\item Every face  of a simple polytope  is 
itself a simple polytope.
\item Every face of a smooth polytope  is 
itself a smooth polytope. 
\end{itemize}

%%%%%%%%%%%%%%%%%%%%%%%%%%%%%%%%%%%%%%%%
\subsection {The volume and moment of $\Delta$.}
\labell{ss:volume}
%%%%%%%%%%%%%%%%%%%%%%%%%%%%%%%%%%%%%%%%%

The proofs in this section were greatly influenced by 
the point of view in
Timorin~\cite{Tim}.

Fix $H \in \ft$ 
and a simple polytope  
$\Delta \subset \ft^*$.
Fix an identification of $\ft^*$ with  Euclidean space, and
consider the following
polynomials of the support numbers:
$$
V=\int_\Delta dx
\quad\mbox{and} \quad 
\mu =\int_\Delta H(x)dx.
$$
That is, $V$ is the volume of $\Delta$ and 
$\mu$  is the  first
moment of $\Delta$ in the direction given by $H$;
on $\Cc_\De$  
they are polynomials of degree
at most
$\dim \Delta$ and $\dim \Delta + 1$
respectively.

More generally, for any
face $f$ of $\Delta$, fix an identification of $P(f)$ with
Euclidean space.
Let  $V_f$ denote the 
volume of $f \subset P(f)$;  analogously, let 
$\mu_f$ denote the integral of the function $H$ over the face $f$.
On $\Cc_\De$, $V_f$ and
$\mu_f$ are polynomials of degree 
at most
$\dim f$ and $\dim(f) +1$ 
respectively.
Let $\p_i$ denote 
the operator of differentiation with respect to $\kappa_i$.

\begin{rmk}\labell{rmk:Kf} \rm 
If $\Delta$ is a smooth polytope, then 
it is natural to choose identifications of $\ft^*$ and $P(f)$
which respect the lattice.  If we make this choice,
the constant $K_f$ in the proposition below is $1$ for
every non-empty face.

However, if $\Delta$ is not rational
then $P(f)$ has no lattice and there is no preferred
identification with Euclidean space.
While the center of mass of $f$ does not depend on this choice,
the constant $K_f$ in the proposition below does.
\end{rmk}

\begin{proposition}\labell{prop:VI}
\label{der}
Let $f$ be the intersection of the facets 
$F_{i_1},\dots,F_{i_k}$. 
Then
$$
\p_{i_1}\cdots\p_{i_k}V=K_f V_f
\quad \mbox{and} \quad
\p_{i_1}\cdots\p_{i_k}\mu=K_f\mu_f,
$$
where $K_f$ is a positive real number 
that depends only on $f$.
In particular, if
the intersection of facets $F_{i_1},\dots,F_{i_k}$ is empty,
then
$$
\p_{i_1}\cdots\p_{i_k}V=\p_{i_1}\cdots\p_{i_k}
\mu=0.
$$
\end{proposition}

\begin{proof}
Let  $f = F_1 \cap \cdots \cap F_k$.
Let $e_1,\ldots, e_n$ denote the standard basis of $\R^n$,
where $n = \dim \ft$. 

If $f$ is not empty, then $\eta_1,\ldots,\eta_k$ are linearly independent.
Moreover, changing the  identifications of $\ft^*$ and $P(f)$ with  Euclidean space
simply 
alters the measure by a positive constant;  this 
only affects the value of $K_f$.  Hence, we may  
identify $\ft^*$  with $\R^n$ so that $\eta_{1},\ldots,\eta_k$  are identified with
$e_1,\ldots,e_k$, 
and thus identify $P(f)$ with $\{0\} \times \R^{n-k}$. 
We can write $V$ as the integral
$$
V = \int_{-\infty}^{\kappa_{k}} \cdots \int_{-\infty}^{\kappa_1}
V(x_1,\ldots, x_k)\; dx_1 \cdots dx_k,
$$
where $V(c_1,\ldots,c_k)$ is the volume of 
$\Delta \cap \big( \{c\}\times \R^{n-k} \big)$ where $c=(c_1,\dots,c_k)$.
Therefore,  $\p_1 \cdots \p_k V = V(\kappa_1,\ldots,\kappa_k) = V_f$
by the fundamental theorem of calculus.
Similarly, we can write 
$\mu$  as the integral
$$
\mu= \int_{-\infty}^{\kappa_{k}} \cdots \int_{-\infty}^{\kappa_1}
\mu(x_1,\ldots, x_k)\; dx_1 \cdots dx_k,
$$
where 
$\mu(c_1,\ldots,c_k)$ is the  
integral of $H$ over the polytope
$\Delta \cap \big( \{c\}\times \R^{n-k} \big)$.
Hence, $\p_1 \cdots \p_k  
\mu=\mu(\kappa_1,\ldots \kappa_k) =
\mu_f$.

So suppose instead that $f$ is empty.
Let  $g = F_2 \cap \cdots \cap F_k$.
We may assume by induction
that $\p_2 \cdots \p_k V = K_g V_g$
for some positive constant $K_g$,
and so
$\p_1 \cdots \p_k V = K_g  \p_1 V_g$.
However, since $F_1$ and $g$ do not intersect, $g$
does not depend on $\kappa_1$.
Therefore, $\p_1 V_g = 0$.
A similar argument shows that
$\p_1 \cdots \p_k
\mu= 0.$
\end{proof}

Given $H \in \ft$ and a simple polytope $\Delta \subset \ft^*$,
define a 
function $\Hat H: \Cc_\Delta  \to \R$ by
\begin{equation}\labell{eq:HatH}
\Hat H(\kappa_1,\dots,\kappa_N)= 
\langle H,c_{\Delta}(\kappa_1,\dots,\kappa_N) \rangle,
\end{equation}
where $c_\Delta$ denotes the center of mass of $\Delta$.
We will repeatedly use the fact that
$$ \mu = \Hat H  V.$$

The next lemma  is a  simple consequence of the fact that,
since $V$ and 
$\mu$ are both polynomials, the function $\Hat H =  
\frac{\mu}{V}$  
is a rational function.

\begin{lemma}\labell{le:local}
Fix $H \in \ft$ and a simple polytope $\Delta \subset \ft^*$.
If the function
$\Hat H $ defined above
is linear on some open subset of $\Cc_\De$, it is linear on all of
$\Cc_\Delta$.
\end{lemma}

\begin{rmk}\labell{rmk:rational1}\rm
If $\Delta \subset \ft^*$ is a smooth polytope, then
every mass linear function $H \in \ell$
 has rational coefficients.
(Recall that $\ell$ is the integer lattice in $\ft$.)
In other words,
if $\langle H, c_\Delta(\kappa) \rangle = \sum \ga_i \, \ka_i$ 
for all $\kappa \in \Cc_\Delta$
then $\ga_i \in \Q$.

To see this, note that by Proposition~\ref{prop:VI} and Remark~\ref{rmk:Kf},
the volume $V \colon \Cc_\Delta \to \R$ is a polynomial with 
rational coefficients. 
Similarly, since  
$\langle H, v \rangle \colon \Cc_\Delta \to \R$ is a linear function  with integer coefficients 
for every vertex $v \in \Delta$, the moment 
$\mu$ is
also a polynomial with rational coefficients. So 
$
\frac{\mu}{V}$ is a rational function with
rational coefficients.
\end{rmk}

We include the following lemma for completeness, but make no use of it  in this paper.

\begin{lemma}\labell{le:rat3}
Let $\Delta \subset \ft^*$ be a smooth polytope. Then the
set of  mass linear functions on $\Delta$ is a rational subspace of $\ft$.
\end{lemma}

To prove it, we will need the following useful lemma, which
shows that the coefficients of a mass linear function $H \in \ft$  
determine the function $H$ itself; cf.  Proposition~\ref{prop:inessential}.

\begin{lemma}\labell{le:Hsum}
Fix $H \in \ft$ and 
a simple polytope $\Delta \subset \ft^*$. 
If
$ \langle H, c_\Delta(\kappa) \rangle = \sum \beta_i \kappa_i,$
then $$H = \sum \beta_i \eta_i.$$
\end{lemma}

\begin{proof}{}
Given any $\xi \in \ft^*$, let 
$\Delta' = \Delta + \xi$, the translate of $\De$ by $\xi$.
Then 
\begin{equation}\labell{shift}
\Delta' = \Delta(\kappa') := \bigcap_{i = 1}^N 
\{ x \in \ft^* \mid \langle \eta_i,x\rangle \leq \kappa'_i  \},
\quad \mbox{where} \quad \kappa'_i = \kappa_i + \langle  \eta_i,\xi \rangle 
\quad \forall \, i.
\end{equation}
Hence,
by assumption 
\begin{eqnarray*}
\langle H, c_\Delta(\kappa') \rangle 
&=& \sum \beta_i \kappa_i'\\
&=&  
\sum \beta_i \kappa_i + \sum \beta_i \langle \eta_i,\xi \rangle
\\
&=& \langle H, c_\Delta(\kappa) \rangle 
 + \sum \langle  \beta_i \eta_i ,\xi\rangle.
 \end{eqnarray*}
Since clearly
$\langle H, c_\Delta(\kappa') \rangle = \langle H, c_\Delta(\kappa) \rangle +
 \langle H, \xi \rangle,$ this implies that
$$
\langle H, \xi \rangle = 
\sum \langle  \beta_i \eta_i,\xi \rangle \quad \forall \, \xi \in \ft^*,
$$
as required.
\end{proof}

\NI
{\bf Proof of Lemma~\ref{le:rat3}}
For all $i \in \{1,\ldots,N\}$,
let $\mu_i$
denote the moment in the $\eta_i$ direction.
By the previous remark, the polynomial 
$\mu_i - V \kappa_i$
is a polynomial with rational coefficients for all $i \in \{1,\ldots,N\}$.
Therefore, the set 
$$ 
X = \left\{ \beta \in \R^N \ \left| \  
\sum_{i=1}^N  \be_i\mu_i = \sum_{i=1}^N \be_i V \ka_i
\right. \right\}
$$ 
is a rational subspace of $\R^N$.
Hence, 
to prove the result, 
it is
enough to show that  $H \in \ft$ is mass linear exactly
if $H = \sum \beta_i \eta_i$, where $\beta \in X$.
To check this, note that
$\mu_i/V = \langle \eta_i, c_\De(\ka)\rangle$. Hence
$$ 
\left\langle \sum \beta_i \eta_i, c_\Delta(\kappa) \right\rangle
= \sum \beta_i \frac{\mu_i}{V}.
$$
Therefore, if $\beta \in X$ then
$ \langle \sum \beta_i \eta_i, c_\Delta(\kappa) \rangle
= \sum \beta_i \kappa_i.$
Conversely, consider $H \in \ft$ so that
$\langle H, c_\Delta(\kappa) \rangle = \sum \beta_i \kappa_i$.
By  Lemma~\ref{le:Hsum},
this implies that
$H = \sum \beta_i \eta_i$.
But then
$$\sum\langle \beta_i \eta_i, c_\Delta(\kappa) \rangle
= \sum \beta_i \kappa_i,
$$
 which implies that $\beta \in X$.
\QED

%%%%%%%%%%%%%%%%%%%%%%%%%%%%%%%%%%%%%%%%%%%%%%%%%%%%%%%%%%%%%%
\subsection{Symmetric faces}\labell{ss:sym}
%%%%%%%%%%%%%%%%%%%%%%%%%%%%%%%%%%%%%%%%%%%%%%%%%%%%%%%%%%%%%%

We are now ready to consider
symmetric faces; see Definition~\ref{def:symmetric}.
In particular, we will  analyze the restriction
of a mass linear function to a symmetric face.

\begin{lemma}\labell{le:symcent} 
Fix $H \in \ft$
and 
let $f$ be a symmetric face of
a simple polytope  $\Delta \subset \ft^*$.
Then
$$
\langle H,c_f(\kappa) \rangle  = \langle H, c_{\Delta}(\kappa) \rangle \quad
\forall \, \kappa \in \Cc_\Delta,
$$
where 
$c_f \colon \Cc_\Delta \to \ft^*$
denotes the center of mass of $f$. 
\end{lemma}

\begin{proof}
Let $f = F_J$ be a symmetric face; 
number the facets so
that $J = \{1,\ldots,j\}$.
Since $F_i$ is symmetric
for each $1 \leq i \leq j$,
$\p_i \Hat H = 0$.
Hence, by 
Proposition~\ref{der}, 
if we apply the differential operator $\p_1 \cdots  \p_j$ to the formula 
$\mu=\Hat H V$ we obtain
$$ K_f \mu_f =   K_f \Hat H V_f.$$
Since 
$K_f \neq 0$ and 
$\mu_f= \langle H, c_f \rangle V_f$, 
this proves
that  $\langle H, c_\Delta \rangle = \langle H, c_f \rangle.$
\end{proof}

\begin{lemma}\labell{le:symasym}
Fix $H \in \ft$ and
let $\Delta \subset \ft^*$ be a simple polytope.
If $f$ is a symmetric face and $F$ an asymmetric facet,
then $f \cap F \neq \emptyset$; indeed, $F \cap f$ is a facet
of $f$.
\end{lemma}

\begin{proof}
Suppose on the contrary that $F \cap f = \emptyset$.
Since $f(\kappa)\subset \Delta(\kappa)$ 
does not depend on support number of $F$; neither does
$\langle H, c_f(\kappa) \rangle \colon
\Cc_\Delta \to \R$.
By the lemma  above, this implies that 
$\langle H, c_\Delta(\kappa) \rangle \colon \Cc_\Delta \to \R$
does not depend on  support number of $F$, that is, $F$ is symmetric. 
This is a contradiction.
Finally, we note that by 
Definition~\ref{def:symmetric}, 
$f \not\subset F$.
Hence $F \cap f$ must be a facet of $f$.
\end{proof}

\begin{prop}
\labell{prop:symface}
Let $H \in \ft$ be a mass linear function on
a simple polytope $\Delta \subset \ft^*.$ 
Let $f$ be a symmetric face of $\Delta$.
Then the restriction of $H$ to $f$ is  mass linear and
the map $F \mapsto F \cap f$ induces a one-to-one
correspondence between the asymmetric facets of $\Delta$
and the asymmetric facets of $f$.
Moreover, the coefficient of the support
number of $F$ in $\langle H, c_\Delta \rangle$
is the coefficient of the support number of $f \cap F$ in $\langle
H, c_f \rangle$.
\end{prop}

\begin{proof}
Renumber the facets so that $F_j \cap f$ is a facet of
$f$ exactly if  
$j\in \{1,\ldots,M\}$.
In the statement of this proposition, we consider $f$ as
a polytope in its own right,  that is,
$f = \bigcap_{j = 1}^M \{x \in P(f) \mid \langle \eta_j,  x \rangle \leq \kappa_k \}$,
where $P(f)$ denotes the affine span of $f$  as in  
Equation (\ref{eq:Pf}).
Hence, we consider the center of mass of $f$ as a function on 
$\Cc_f \subset \R^M$.
In contrast, in Lemma~\ref{le:symcent}, we consider $f$ as a face of the
polytope 
$\Delta$ and 
consider $c_f$
as  a function on 
$\Cc_\Delta \subset \R^N$.  
However, if $F$ is an asymmetric facet of $\Delta$
then by Lemma~\ref{le:symasym} $F \cap f$ is a facet of $f$.
Thus
the facets $F_i$ for $i>M$ are symmetric, so that 
$\langle H, c_\Delta \rangle \colon \Cc_\Delta \to \ft^*$ does
not depend on the support number $\kappa_i$.
Therefore,
Lemma~\ref{le:symcent} implies that
\begin{equation}\labell{ref}
 \langle H, c_f(\kappa_1,\ldots,\kappa_M) \rangle = \langle H, c_\Delta 
(\kappa_1,\ldots,\kappa_N)  
\rangle \quad \forall \, \kappa \in \Cc_\Delta. 
\end{equation}
Next observe that the set of $M$-tuples $(\kappa_1,\ldots,\kappa_M)$
that extend to $N$-tuples $(\kappa_1,\ldots,\kappa_N) \in \Cc_\De$
forms an open subset $U$  of $\Cc_f$ that may not be the whole of $\Cc_f$.  Nevertheless,  
by  Lemma~\ref{le:local} the fact that 
the function $\langle H, c_f \rangle$ is linear on $U$
 implies that it is linear everywhere in $\Cc_f$.  Moreover 
 its coefficients  equal those of 
$ \langle H, c_\Delta
\rangle$.  Since a facet is asymmetric if and only if the corresponding coefficient is nonzero, the result is now immediate.
\end{proof}

\begin{rmk}\rm 
The following notion is sometimes useful.
Fix $H \in \ft$ and a simple polytope  $\Delta\subset \ft^*$.
We say that a face $f$ of $\Delta$ is
{\bf centered} if 
$$
\langle H,c_f(\kappa) \rangle  = \langle H, c_{\Delta}(\kappa) \rangle \quad
\forall \, \kappa \in \Cc_\Delta.
$$
Lemma~\ref{le:symcent} shows that symmetric faces are centered.
Conversely, if $F$ is a centered facet, 
then by definition
$F$ is also symmetric.
In contrast, the converse does not hold for 
faces of higher codimension.
For example, let $\Delta_2 = 
\{ x \in \R^2 \mid x_1 \geq 0, x_2 \geq 0, \mbox{and} \ x_1 + x_2  \leq 1 \}$
 denote the standard $2$-simplex, and let $H(x) = x_1 - x_2$.
Then the vertex $(0,0)$ is centered, 
but it is not symmetric
because the facets containing it are not symmetric.
The proof of Lemma~\ref{le:symasym} can be adapted to show that
every centered face must intersect every asymmetric facet.
\end{rmk}

%%%%%%%%%%%%%%%%%%%%%%%%%%%%%%%%%%%%%%%%%%%%%%%%%%%%%%%%%%%%%%
\subsection{Asymmetric facets}\labell{ss:asym}
%%%%%%%%%%%%%%%%%%%%%%%%%%%%%%%%%%%%%%%%%%%%%%%%%%%%%%%%%%%%%%
We now consider asymmetric facets.
The main goal of this subsection
is to prove 
Theorem~\ref{thm:easy} and 
the proposition below;
see Definition~\ref{def:perv}.

\begin{prop}\labell{prop:asym}
Let $H \in \ft$ be a mass linear function
on a simple polytope $\De \subset \ft^*$.
Then every asymmetric facet is 
pervasive or  flat (or both).
\end{prop}

We will need the following lemmas.

\begin{lemma}\labell{le:flatlin}
Let $\Delta \subset \ft^*$ be a simple polytope.
The facet $F_i$ is flat if and only if the volume $V$ is a
linear function of $\kappa_i$.
\end{lemma}

\begin{proof}
If $F_i$ is flat, a straightforward computation shows that
the volume is a linear function of $\kappa_i$.
So assume that $F_i$ is not flat.
At each vertex $v$ of $F_i$ there is a unique edge $e_v$ of $\Delta$
that does not lie in $F_i$. 
Since $F_i$ is not flat, these edges cannot all  be parallel.
In particular, there is an edge $e$ joining two vertices $v$ and $v'$
so that $e_v$ is not parallel to $e_{v'}$.
Write $e = F_J \cap F_i$ for some $J = \{j_1,\ldots,j_{n-2}\} \subset
\{1,\ldots,N\}$.
Since the edges $e_v$ and $e_{v'}$ of the polygon $F_J$ are not
parallel, $V_{F_J}$ is not a linear function of $\kappa_1$;
it has degree $2$.
But by Proposition \ref{prop:VI},
$\p_{j_1} \cdots \p_{j_{n-2}} V = K_{F_J} V_{F_J}$.
Hence, 
$V$ is not a linear function of $\kappa_1$.
\end{proof}

\begin{lemma} \labell{le:twoasym}
Fix a nonzero $H \in \ft$ and  a simple
polytope $\Delta \subset \ft^*$.
Then $\Delta$ has at least two asymmetric facets.
\end{lemma}

\begin{proof}
Let $F_i$ be any facet.
Since $\Delta$ has a nonempty interior, $c_\Delta$ does not lie on $F_i$.
Therefore, by translating and 
then rescaling,
we can find a dilation of $\ft^*$ that fixes the
plane $P(F_i)$  but moves $c_\De$ to an arbitrary nearby point.
Hence $\langle H, c_\Delta \rangle$ depends on at least
one $\kappa_j$ for $j \neq i$.  
\end{proof}

\NI
{\bf Proof of Proposition~\ref{prop:asym}}. 
Let $F_1$ be an asymmetric facet which is not flat
and let $F_2$ be a facet which does not meet $F_1$.
Since $H$ is mass linear, the function
$\Hat H$  of Equation (\ref{eq:HatH}) is linear,  and so its  second derivatives all vanish.
Hence, 
since $F_1 \cap F_2 = \emptyset$,
by 
Proposition~\ref{der} 
if we apply the differential operator $\p_1 \p_2$ to the formula 
$\mu=\Hat H V$ we obtain
$$
0 = (\p_1 \Hat H) V_{F_2} +  (\p_2 \Hat H) V_{F_1}.
$$
Applying the operator $\p_1$ again to the equation above,
we obtain
$$ 
(\p_2 \Hat H )\; \p_1 V_{F_1} = 0.
$$
On the one hand, Lemma~\ref{le:symasym} implies that $F_2$ is asymmetric.
Thus $\Hat H$ depends linearly on the support number of
$F_2$, so that
$\p_2 \Hat H \neq 0$.
On the other hand, since $F_1$ is not flat Lemma~\ref{le:flatlin} implies that
$\p_1 V_{F_1} \neq 0$; 
this gives a contradiction.
\QED\MS

\NI
{\bf Proof of Theorem~\ref{thm:easy}}. 
Suppose that the simple polytope  $\De \subset \ft^*$  contains
no pervasive facets and no flat facets. It suffices to show that every
homogeneous   function $H\in \ft$ that is mass linear on $\De$ must vanish.
But if such 
an
$H$ is nonzero, then   $\Delta$ contains an asymmetric facet $F$
by
Lemma~\ref{le:twoasym}.  Moreover, $F$ must  be pervasive or flat
by  Proposition~\ref{prop:asym}, which is impossible.
\QED\MS

%%%%%%%%%%%%%%%%%%%%%%%%%%%%%%%%%%%%%%%%%%%%%%%%%%%%%%%%%%%%%%%%
\section{Polytopes with inessential functions}\labell{s:ines}
%%%%%%%%%%%%%%%%%%%%%%%%%%%%%%%%%%%%%%%%%%%%%%%%%%%%%%%%%%%%%%%%

The main results  of this section are 
Proposition~\ref{prop:Iexpan}, which 
characterizes  polytopes with nonconstant 
inessential functions,  
and Theorem~\ref{thm:allmass},
which characterizes polytopes
with the property that all linear functions are mass linear.

%%%%%%%%%%%%%%%%%%%%%%%%%%%%%%%%%%%%%%%%%%%%%%%%%%%%%%%%%%%%%%%%
\subsection{Symmetries and equivalent facets}\labell{ss:ines}

%%%%%%%%%%%%%%%%%%%%%%%%%%%%%%%%%%%%%%%%%%%%%%%%%%%%%%%%%%%%%%%%%%

In this subsection, we prove several useful criteria for determining if facets are equivalent.
To begin, we
 justify our normalization conventions by proving 
a few easy facts about symmetries: see Remark~\ref{rmk:normalize}.
Because this is not yet accomplished,
we will {\bf not} assume that our conormals are normalized in Lemma~\ref{le:symmetries}, 
though we do make this assumption in subsequent arguments.

\begin{lemma}\label{le:symmetries}
Let $a \colon \ft^* \to \ft^*$ be a symmetry of a simple polytope $\Delta \subset \ft$.
Let $\Hat{a} \colon \ft \to \ft$ be the dual to the associated linear map.
Let $\eta_1,\dots,\eta_N$ be {\bf any} (not necessarily normalized) outward conormals. 
\begin{itemize}
\item $a(F_i) = F_j$ exactly if
$\Hat{a}^*(\eta_i) = \lambda \eta_j$ for some $\lambda  > 0$.
\item $a(F_i) = F_i$ exactly if
$\Hat{a}^*(\eta_i) =  \eta_i$.
\item $a(F_i) = F_i$ for all $i$ exactly if $a = \id$.
\end{itemize}
\end{lemma}

\begin{proof}
The first claim is obvious.
To prove the second claim,
we first note that if  $\Hat a^*(\eta_i) = \la\eta_i$ for some $\la > 0$,  
then
$$ 
\langle \eta_i, \Hat a(x) \rangle = \langle \Hat a^*(\eta_i), x \rangle
= \langle \lambda \eta_i, x \rangle = \lambda \langle \eta_i,x \rangle
\quad \forall \,x\in \De.
$$
On the other hand, since $a_\kappa$ is a symmetry of $\Delta$, there exists $c \in \R$ such that $\Hat a(\Delta) = \De + c$.
Therefore
$$
\la \, \sup \bigl\{\langle \eta_i, x \rangle\,|\, x\in \De\bigr\}  = 
\sup \bigl\{\langle \eta_i, \Hat a(x) \rangle\,|\, x\in \De\bigr\}  = 
\sup \bigl\{\langle \eta_i, x \rangle\,|\, x\in \De\bigr\}  + c.
$$
Thus if $S: = \sup \bigl\{\langle \eta_i, x \rangle\,|\, x\in \De\bigr\} $ we have $\la S = S + c$.  Similarly,  $\la s = s+c$, where  $s: = 
\inf \bigl\{\langle \eta_i, x \rangle\,|\, x\in \De\bigr\} $.
Since  $S\ne s$, 
this is impossible unless $\lambda = 1$.
Finally, because $\Delta$ is bounded the $\eta_i$'s span $\ft^*$,
and so the identity is the only symmetry which takes each facet to itself. 
\end{proof}

\begin{example}\labell{ex:rat4}\rm
Consider the triangle  with vertices
$(0,0)$, $(3,0)$, and 
$(0,2)$.
Note that this polytope is rational but not smooth.
The linear map $\Hat a: (x_1,x_2) \mapsto 
(-x_1 - \frac{3}{2} x_2, \frac{2}{3} x_1)$
is a robust symmetry that gives rise to the symmetry
$\Hat a + (3,0)$ of $\De$. Since its transpose $\Hat a^*$
takes $(-3,0)$ to $(0,-2)$, $(0,-2)$ to $(3,2)$ and $(3,2)$ to $(-3,0)$,
 normalized outward conormals have the form 
 $\la(-3,0), \la (0,- 2), $ and $\la (3,2)$ for some $\la >0$.
 \end{example}

We will also need the following lemma.

\begin{lemma}\label{le:smoothxi}
Let $F_i$ and $F_j$ be facets of a smooth polytope $\Delta \subset \ft^*$.
If there exists $\xi \in \ft$ such that
$$ \langle\eta_i,\xi\rangle > 0 > 
\langle\eta'_j,\xi\rangle
\quad \mbox{and otherwise }  \quad
\langle \eta_k , \xi \rangle = 0 \ \forall \, k,
$$
then there exists $\xi' \in \ell^*$ such that

$$
\langle\eta_i,\xi' \rangle = 1 = - \langle\eta_j,\xi' \rangle
\quad \mbox{and otherwise }  \quad
\langle\eta_k,\xi' \rangle = 0\;\, \forall \, k.
$$
\end{lemma}

\begin{proof}
By the definition of smooth, the primitive outward conormals at
any vertex of $\Delta$ form a basis for $\ell$.  
By looking at vertices in  $F_i\less F_j$ and $F_j\less F_i$ 
one sees that there is a basis of $\ell$ that contains $\eta_i$ but not $\eta_j$ 
and also a basis that contains $\eta_j$ but not $\eta_i$.  
Therefore if we choose $\la>0$ so that  $\xi': = \la\xi$ satisfies 
$\langle\eta_i,\xi'\rangle=1$, then $\xi'$ lies in $\ell^*$ and is primitive;  hence we must
also have $\langle\eta_j,\xi'\rangle=-1$.
\end{proof}

We shall say that a linear map $\Hat a:\ft^*\to \ft^*$ is  a {\bf reflection} 
if its fixed point set is a hyperplane and  
$\Hat a\circ \Hat a = \id$.

\begin{lemma}\labell{le:xi}
Let $F_i$ and $F_j$ be distinct facets of 
a simple  polytope $\Delta \subset \ft^*$.
The following are equivalent:
\begin{enumerate}
\item [(i)]
$F_i$ is equivalent to $F_j$.
 
\item [(ii)]
There exists a vector $\xi \in \ft^*$ 
satisfying 
$$
\langle\eta_i,\xi \rangle > 0 >  \langle\eta_j,\xi \rangle
\quad \mbox{and otherwise }  \quad
\langle\eta_k,\xi \rangle = 0\;\, \forall \, k. $$

\item [(iii)]
There exists a vector $\xi \in \ft^*$ 
satisfying 
$$
\langle\eta_i,\xi \rangle = 1 = - \langle\eta_j,\xi \rangle
\quad \mbox{and otherwise }  \quad
\langle\eta_k,\xi \rangle = 0\;\, \forall \, k.
$$

\item [(iv)]
There exists a reflection $\Hat{a} \colon \ft^* \to \ft^*$ which is a robust symmetry of $\Delta$;
for all $\kappa \in \Cc_\Delta$ the
associated symmetry $a_\kappa$ 
exchanges $F_i$ and $F_j$ and otherwise takes each facet to itself.
\end{enumerate}
\end{lemma}

\begin{proof}
It is obvious that (iii) implies (ii) and that (iv) implies (i).

Our first step is to show that  (ii) implies  (iv).
Assume that there exists $\xi \in \ft^*$ 
such that $ \langle  \eta_1,\xi \rangle >0 >  \langle  \eta_2,\xi \rangle$ 
and $\langle \eta_k, \xi \rangle = 0$ for all $k > 2$.
Let $$
\eta'_1 = \frac{\eta_1}{ \langle \eta_1, \xi \rangle } 
\quad \mbox{and}  \quad
\eta'_2 = \frac{\eta_2}{- \langle \eta_2, \xi \rangle };
$$
clearly
$
\langle\eta'_1,\xi \rangle = 1 = - \langle\eta'_2,\xi \rangle.
$
Define a reflection $\Hat{a} \colon \ft^* \to \ft^*$ by 
$$
\Hat{a}(x) = x - \langle \eta'_1 - \eta'_2, x \rangle \,\xi.
$$
Given $\kappa \in \Cc_\kappa$, define an affine
map
$a_\kappa \colon \ft^* \to \ft^*$ by
$$
a_\kappa(x) = 
\Hat{a}(x)   +\left(\kappa'_1 - \kappa'_2 \right) \xi,
\ \  \mbox{where}  \ \
\kappa'_1 = \frac{\kappa_1}{ \langle \eta_1, \xi \rangle }
\ \  \mbox{and}  \ \ 
\kappa'_2 = \frac{\kappa_2}{- \langle \eta_2, \xi \rangle } .$$
It is easy to check that $a_\kappa$ carries $\Delta(\kappa)$ to
itself, exchanging $F_1$ and $F_2$ and otherwise taking each facet to
itself.  
Moreover, if $\Delta$ is smooth, then by Lemma~\ref{le:smoothxi} we may
assume that $\xi$ lies in $\ell^*$ and that $\langle \eta_1, \xi \rangle = 1 = - \langle \eta_2, \xi \rangle$.
Hence, $\Hat{a}$ induces an isomorphism of the integer lattice,
and so it is a robust symmetry of $\Delta$.

Finally, we will prove that 
(i) implies (iii).
Let
 $\Hat{a} \colon \ft^* \to \ft^*$ be a 
robust symmetry
 and let $\Hat{a}^* \colon \ft \to \ft$ denote
the dual map.  
For each $\kappa \in \Cc_\De$,  
let
$a_\kappa \colon \ft^* \to \ft^*$ be the associated
symmetry of $\Delta(\kappa)$.
We can naturally partition the set of facets of $\Delta$ into 
$a_\kappa$-orbits;
let $\Oo$ denote the set of such orbits.
Assume that each $a_\kappa$ takes $F_1$ to $F_2$.

Since $a_\kappa(\Delta) = \Delta$, $a_\kappa^k = \id$ for some $k > 0$.
Define $p \colon \ft \to \ft$ by
$$p = 1 + \Hat{a}^* + \cdots + (\Hat{a}^*)^{k-1}.$$
We shall find the required vector $\xi$ by studying the properties of  
this map $p$.

By our normalization convention (see Remark~\ref{rmk:normalize}),
$p(\eta_j) = \frac{k}{|O|} \sum_{i \in O} \eta_i$ for all $j$,
where $O$ is the orbit containing $F_j$.
Moreover, since $\Delta$ is bounded,
the $\eta_i$ span $\ft$.  Hence,
\begin{gather}
\labell{eq:kerp}
\ker p  = \left\{  \sum \beta_i \eta_i  \in \ft \left| \ \sum_{i \in O} \beta_i = 0
\ \ \forall \, O \in \Oo  \right. \right\}, \quad \mbox{and} \\
\labell{eq:imp}
p(\ft)= \left\{ \left. \sum_{O \in \Oo}
\beta_O
\bigg( \sum_{i \in O} \eta_i\bigg)
\ \right| 
\ \beta_O \in \R \ \forall \, O \in \Oo \right\}.
\end{gather}
Now fix   
a nonzero 
$\beta \in \R^N$ such
that $\sum_{i \in O}\beta_i = 0 $ for all $O \in \Oo$,
and let $H  = \sum \beta_i \eta_i$.
Then  Proposition~\ref{prop:inessential} implies that
\begin{equation}\labell{eq:Hsemi}
\langle H, c_{\Delta}(\kappa) \rangle =  \sum \beta_i \kappa_i.
\end{equation}
Since this equality holds for  all $\kappa$ in the open set 
 $\Cc_\Delta$, 
this implies that $H \neq 0$.
Therefore, $\ker p$ has dimension $N - |\Oo|$ by \eqref{eq:kerp},
and so   $p(\ft) \subset \ft$ has codimension $N - |\Oo|$. 
On the other hand, 
let $W \subset \ft$ be the subspace spanned by $\eta_1 + \eta_2$
and $\{\eta_k\}_{k > 2}$.
By \eqref{eq:imp}, 
the fact that  $F_1$ and $F_2$ are in the same orbit
implies that
$p(\ft)$ is a subspace of $W$ of codimension  at most $N - |\Oo| - 1$.
Therefore, $W$ is a proper subspace of $\ft$. 
Since the $\eta_i$ span $\ft$, this implies that there exists $\xi \in \ft^*$ so that
$$
\langle\eta_1,\xi\rangle = 1 = - \langle\eta_2,\xi\rangle
\quad \mbox{and otherwise }  \quad
\langle\eta_k,\xi\rangle = 0\;\, \forall \, k.
$$
This completes the proof.
\end{proof}

The above lemma implies that the equivalence relation $\sim$ depends only on the set of conormals  of $\De$ and not on the intersection pattern of 
its facets.  In other words, it depends only on the $1$-skeleton of the corresponding fan, not on its higher dimensional simplices.

This lemma (or more precisely, its proof) has the following corollary,
which will be useful in \S\ref{s:geometry}.

\begin{cor} \labell{cor:xi} 
Let $\Delta \subset \ft^*$ be a simple polytope.
Let  $\Hat a \colon \ft^* \to \ft^*$ be a linear map so
that for all $\kappa'$ in some neighborhood $\Uu$ of $\ka$ in 
$\Cc_\De$ there exists an affine map $a_\kappa$ that
satisfies conditions (1) and (2) of Definition~\ref{def:robust}.
Then $\Hat a$ is a robust symmetry of $\Delta$.
\end{cor}

\begin{proof} 
Suppose that the $a_\kappa$ take $F_i$ to $F_j$.
The above proof that (i) implies 
(iii)
goes through in this situation
since all we use is that equation (\ref{eq:Hsemi}) holds  for all $\ka$ in some open set
$\Uu$. Since 
(iii) implies (ii) and (ii) implies (iv),
there is a robust symmetry that 
exchanges $F_i$ and $F_j$ but takes every other facet to itself.  
By composing such maps, we can construct a robust symmetry
$\Hat b$  so that $\Hat a$ and $\Hat b$ induce the same permutation on
the facets.
By Lemma~\ref{le:symmetries},
this implies that $\Hat a = \Hat b$.
\end{proof}

\begin{rmk}\labell{rmk:p}\rm 
Let $\Delta$ be a simple polytope.
We claim that $H \in \ft$ is inessential on $\Delta$ exactly
if there exists a robust symmetry $\Hat{a}$ of $\Delta$ so that
$H$ is in the kernel of the associated projection
$$p = 1 + \Hat{a}^* + \cdots + (\Hat{a}^*)^{k-1},$$
that is,
so that $H$ has zero average with respect to the action
of $\Hat{a}^*$.
Here, $\Hat{a}^* \colon \ft^* \to \ft^*$ is the dual
map and $k> 0$ is chosen so that ${a}^k = \id$. 

To see this, note that the lemma above immediately implies
that there exists a robust symmetry $\Hat{a}$ of $\Delta$
whose set $\Oo$ of orbits of facets is precisely the set $\Ii$
of equivalence classes of facets. 
Moreover, given any robust symmetry $\Hat{a}$ of $\Delta$,
 all the facets
 in any orbit $O \in \Oo$ are equivalent.
The claim now follows from the formulas (\ref{eq:iness}) and
(\ref{eq:kerp}). 
 \end{rmk}

The next lemma gives a very useful 
criterion for a set of facets to be equivalent.

\begin{lemma}\labell{le:equiv}
Let $\Delta \subset \ft^*$ 
be a simple polytope.
Given a subset $I \subset \{1,\ldots,N\}$,
$F_i \sim F_j$ for all $i$ and $j$ in $I$
exactly if the subspace 
$W\subset \ft$ spanned
by the outward conormals $\eta_k$ for $k \not\in I$
has codimension $|I| - 1$. 
Moreover, in this case 
the linear combination
$\sum_{i \in I} c_i \eta_i$ lies in $W$
if and only if $c_i = c_j$ for all $i$ and $j$.
\end{lemma}

\begin{proof}
Assume that $F_i \sim F_j$ for all $i$ and $j$  in $I$.
By Lemma~\ref{le:xi},
for each such pair 
there exists a vector $\xi_{ij} \in \ft^*$  satisfying 
$\langle  \eta_i,\xi_{ij} \rangle = 1 = -
\langle  \eta_j,\xi_{ij} \rangle$ and
otherwise $\langle \eta_k, \xi_{ij} \rangle = 0$
for all $k$.
Fix $i\in I$.  Then the vectors $\xi_{ij}$, $j\in I\less\{i\}$,
span an
$(|I| -1)$-dimensional subspace of $\ft^*$;
let 
$W\subset \ft$ be its annihilator. 
By construction, $\eta_k$ lies in 
$W$ for  all $k \not\in I$.
On the other hand, since $\Delta$ is bounded
the positive span
of the $\eta_k$ contains all of $\ft$.
(Recall that the {\bf positive span} of $\eta_1,\ldots,\eta_N$ is the set
$\{ \sum a_i \eta_i \mid a_i > 0 \ \forall \, i \}.$)
In particular, any $\eta_k$
can be written as a linear
combination of the other outward conormals. Hence, the subspace spanned
by the $\eta_k$ for $k \not \in I$ has
codimension at most $|I| -1$; so it must be  
$W$.
Finally, since 
$\langle \sum_{k \in I} c_k \eta_k, \xi_{ij} \rangle =
c_i - c_j$,  
$\sum_{k \in I} c_k \eta_k$ 
lies in 
$W$ exactly if $c_i = c_j$ for all $i$ and $j$.

Conversely,
assume that the plane 
$W \subset \ft$ spanned by the
outward conormals $\eta_k$ for $k \not\in I$ has codimension
$|I| - 1$.  
Hence for any pair $i \neq j$ in $I$,
the subspace spanned by $W$ {\em and}  the remaining $|I|-2$ elements
in  $\{\eta_k, k\in I\}$   
is a proper subspace of $\ft$.
Hence there exists 
$\xi \in \ft^*$ so that 
$\langle  \eta_k,\xi \rangle = 0$ for 
all $k$ except $i$ and $j$. 
Since the positive span of the $\eta_k$
contains all of $\ft$, the positive span of
the 
numbers $\langle  \eta_k,\xi \rangle$, namely
$$
\Big\{\sum a_k \langle  \eta_k,\xi \rangle \Big| \ a_k  > 0 \Big\} = \bigl\{a \langle  \eta_i,\xi \rangle +
b \langle  \eta_j,\xi \rangle \mid a,b>0\bigr\},
$$
contains all of $\R$.
In particular,  
after possibly replacing $\xi$ by $-\xi$,
we have $\langle  \eta_i,\xi \rangle > 0 > \langle  \eta_j ,\xi\rangle $.
By Lemma~\ref{le:xi}, this implies that $F_i$ and $F_j$ are
equivalent.
\end{proof}

\begin{example}\rm \labell{ex:product}
Given simple polytopes $\Tilde \Delta \subset \Tilde \ft^*$
and $\Hat \Delta \subset \Hat \ft^*$, consider the
product $\Delta = \Tilde \Delta \times \Hat \Delta$.
It follows immediately from 
Lemma \ref{le:equiv} above
 that two facets ${\Tilde F}_i \times \Hat \Delta$
and ${\Tilde F}_j \times \Hat \Delta$ of $\Delta$
are equivalent exactly if ${\Tilde F}_i$ and
${\Tilde F}_j$ are equivalent as facets of $\Tilde \Delta$.
Similarly,
two facets
 $\Tilde \Delta \times {\Hat F}_\ell$
and ${\Tilde \Delta} \times {\Hat F}_k$ of $\Delta$
are equivalent exactly if 
${\Hat F}_\ell$ and
${\Hat F}_k$ are equivalent as facets of $\Hat \Delta$.
Therefore, $(\Tilde H, \Hat H) \in \Tilde \ft \times \Hat \ft$
is inessential on $\Delta$ exactly if $\Tilde H$ is inessential
on $\Tilde \Delta$ and $\Hat H$ is inessential on $\Hat \Delta$.
\end{example}

The second part of the next result shows that, in the smooth case, 
if an element of the integer lattice  defines an  inessential function then this function has rational coefficients.

\begin{lemma}\labell{le:rational2}
Let $\Delta \subset \ft^*$ be a simple polytope;
let $\Ii$ denote the set of equivalence classes of facets.
The
set of inessential functions is
an $(N - |\Ii|)$-dimensional  subspace  of $\ft$.
Moreover, if $\Delta$ is smooth then 
this is a rational subspace, and
every
inessential $H \in \ell$
can be written 
$$
H = \sum \beta_i \eta_i, \quad \mbox{where}\ \beta_i \in \Q
\ \forall \, i \ \mbox{and} \ \sum_{i \in I} \beta_i = 0 \ \forall \, I \in \Ii.
$$
\end{lemma}

\begin{proof} 
Let
$H \in \ft$ be inessential. By  
Definition~\ref{def:iness} we may write
 $H = \sum \beta_i \eta_i$, where the $\be_i$ satisfy all the above conditions except that they may be in $\R$ not $\Q$.
Fix $i \in I \in \Ii$.  By Lemma~\ref{le:xi}, for all 
$j \in I \smallsetminus \{i\}$, 
there exists  $\xi_{ij} \in \ft^*$ 
so that
$$ 
\langle\eta_i,\xi_{ij}\rangle = 1 = - \langle\eta_j,\xi_{ij}\rangle
\quad \mbox{and otherwise }  \quad
\langle\eta_k,\xi_{ij}\rangle = 0\;\, \forall \, k.
$$
Then 
$$
|I|\be_i = |I|\be_i -  \sum_{j\in I}\be_j = 
\sum_{j\in I} (\be_i - \be_j) 
= \sum_{j \in I \smallsetminus \{i\}} \langle H, \xi_{ij} \rangle.
$$
If $H = 0$, this implies immediately that $\beta_i = 0$;
this proves the first claim.\footnote
{
This also follows from Proposition~\ref{prop:inessential}, as we saw while proving Lemma~\ref{le:xi}.}

On the other hand, if $\Delta$ is smooth
then, because
the $\{\eta_k\}$ form a basis for $\ell$,
$\xi_{ij}$ lies in the lattice $\ell^*$ for all $j \in I \smallsetminus \{i\}$.
Hence, if $H \in \ell$, then
$|I| \beta_i$ is an integer
for every  $i \in I \in \Ii$; in particular, $\beta_i \in \Q$.
\end{proof}

%%%%%%%%%%%%%%%%%%%%%%%%%%%%%%%%%%%%%%%%%
\subsection{Bundles and expansions}
%%%%%%%%%%%%%%%%%%%%%%%%%%%%%%%%%%%%%%%%%

We will now  define 
two important  class of polytopes: bundles and $k$-fold expansions.
Bundles over the simplex $\De_k$ and $k$-fold expansions 
are very similar; in particular,  both admit
nonconstant
inessential functions.   
Our main result  is
that these two classes of polytopes are
the only ones which admit nonconstant inessential functions.

\begin{defn}\labell{def:bund} 
Let $\Tilde{\Delta} =
\bigcap_{j = 1}^{\Tilde{N}} \{ x \in \Tilde{\ft}^* \mid \langle \Tilde{\eta}_j,x
\rangle \leq \Tilde{\kappa}_j \}$ and
 $\Hat{\Delta} =
\bigcap_{i = 1}^{\Hat{N}} \{ y \in \Hat{\ft}^* \mid \langle \Hat{\eta}_i,y
\rangle \leq \Hat{\kappa}_i \}$ be simple polytopes.
We say that a simple 
polytope 
$\Delta \subset \ft^*$
is a {\bf bundle} with {\bf fiber} $\Tilde{\Delta}$ over the {\bf
base} $\Hat{\Delta}$ if 
there exists  a short exact sequence
$$ 
0 \to \Tilde{\ft} \stackrel{\iota}{\to} \ft \stackrel{\pi}{\to} \Hat{\ft} \to 0
$$
so that
\begin{itemize}
\item 
$\Delta$ is combinatorially equivalent to the product $\Tilde{\Delta}
\times \Hat{\Delta}$. (See Definition~\ref{def:combeq}.)
\item 
If $\Tilde{\eta}_j\,\!'$   
 denotes the outward conormal to
the  facet 
$\Tilde{F}_j\,\!'$
of $\Delta$ which corresponds to  $\Tilde{F}_j \times
\Hat{\Delta} \subset \Tilde{\Delta} \times \Hat{\Delta}$,
then $\Tilde{\eta}_j\,\!' = \iota(\Tilde{\eta}_j)$ 
for all $1 \leq j \leq \Tilde{N}.$
\item 
If $\Hat{\eta}_i\,\!'$  
denotes the outward conormal to
the  facet 
$\Hat F_i\,\!'$
of $\Delta$ which corresponds to  $\Tilde{\Delta} \times
\Hat{F_i} \subset \Tilde{\Delta} \times \Hat{\Delta}$,
then $\pi(\Hat{\eta}_i\,\!') = \Hat{\eta}_i$
for all $1 \leq i \leq \Hat{N}.$
\end{itemize}
 The facets 
$\Tilde{F}_1\,\!', \ldots, \Tilde{F}_{\Tilde{N}}\,\!'$ 
will be called {\bf fiber facets},  
and the facets 
$\Hat F_1\,\!' \ldots, \Hat F_{\Hat N} \,\!'$
will be called {\bf base facets}.
\end{defn}

If $\De$ is smooth then
$\Tilde\De$ and $\Hat\De$ are both smooth.
However, the converse does not hold; see Example~\ref{ex:1}.

\begin{rmk}\rm
This terminology is justified by the following fact. 
As we explain in greater detail 
in \S\ref{s:geometry},
there is a manifold
$M_\Delta$ associated to every smooth polytope $\Delta$.
By Remark~\ref{rmk:geobun},  if $\Delta$ is a 
bundle with  fiber $\Tilde{\Delta}$ over the 
base $\Hat{\Delta}$ in the sense defined above,
then $M_\Delta$ is a bundle with
fiber $M_{\Tilde{\Delta}}$ over the base $M_{\Hat{\Delta}}$
in the usual sense.
\end{rmk}

The example below illustrates an important, but slightly confusing,
point.  The fiber facets are almost never 
equivalent in any sense to 
the fiber polytope $\Tilde\De$. Indeed, if the fiber
$\Tilde\De$ is a one-simplex,
the fiber facets
$\Tilde F_j\,\!'$
are analogous to the base polytope $\Hat\De$.
In other words we may identify 
$P(\Tilde F_j\,\!')$
 with $\Hat\ft^*$ in such a way as to set up a combinatorial equivalence between  
 $\Tilde F_j\,\!'$
and $\Hat\De$ in which corresponding facets are parallel; cf. the discussion just before Example~\ref{ex:1}.
More generally, 
if the fiber $\Tilde\De$ has dimension $k$ and the base $\Hat\De$ has dimension $n$, then the
non-empty intersection of any  $k$ fiber facets
is analogous to the base polytope and may be considered as a section of the bundle.  Similarly,
the non-empty intersection of any $n$
base facets is affine equivalent to the fiber polytope.

\begin{figure}[htbp] 
   \centering
\includegraphics[width=4in]{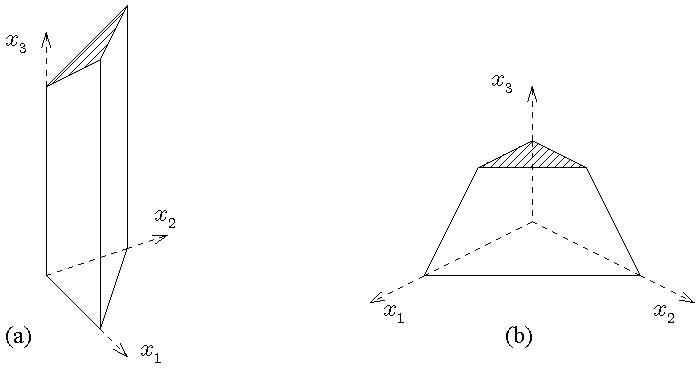} 
   \caption{(a) is a $\Delta_2$ bundle over $\Delta_1$ with a base facet shaded; (b) is a
   $\Delta_1$ bundle over $\Delta_2$ with a fiber facet shaded.
   Notice that the fiber facets in (a) contain a common vector that is  vertical, while the fiber facets in (b) are parallel. 
}
   \label{fig:1}
\end{figure}

\begin{example}\labell{ex:2} \rm
The  polytope $Y_a$ defined in Example~\ref{ex:1}
is a $\Delta_2$ bundle over $\Delta_1$ 
because $\Delta$ 
is combinatorially equivalent to $\Delta_1 \times \Delta_2$
and the conormals $(-1,0,0), (0,-1,0)$, and $(1,1,0)$ are
linearly dependent. 
However, it 
is not a $\Delta_1$ bundle over $\Delta_2$ unless $a_1 = a_2 = 0$.
 
In contrast, the
polytope 
$$  \{ x \in \R^3 \mid x_i \geq 0 \ \forall \,  i, \
 x_1 + x_2  \leq \lambda + b x_3, \mbox{ and } x_3 \leq h \}$$
is a $\Delta_1$ bundle over $\Delta_2$ for all $h > 0$
and $\lambda > \max(0,-bh)$, but is not
a $\Delta_2$ bundle over $\Delta_1$ unless $b = 0$. 
See Figure~\ref{fig:1}.
\end{example}

\begin{defn}\labell{def:expand}
Let $\Tilde{\Delta} =
\bigcap_{j = 1}^{\Tilde{N}} \{ x \in \Tilde{\ft}^* \mid \langle  
\Tilde{\eta}_j, x
\rangle \leq \Tilde{\kappa}_j \}$ be a simple polytope. 
Given a natural number $k$,
a polytope $\Delta \subset \ft^*$ is the {\bf $\mathbf k$-fold expansion}
of $\Tilde{\Delta}$ along the facet $\Tilde{F}_1$ if there
is an identification
$\ft = \Tilde{\ft} \oplus \R^k$ so that
$$
\Delta =
\bigcap_{j = 2}^{\Tilde{N}}\bigl \{ x \in {\ft}^* \mid 
\langle \Tilde{\eta}_j',x \rangle \leq \Tilde{\kappa}_j\bigr\} 
\cap\bigcap_{i = 1}^{k+1} \bigl\{ x \in {\ft}^* \mid \langle  \Hat\eta_i,x
\rangle \leq \Hat\kappa_i \bigr\},
$$ 
where 
$\Tilde{\eta}_j'= (\Tilde{\eta}_j,0)$
for $j\ge 2$,
  $\Hat\eta_i = (0, - e_i)$ and $\Hat\kappa_i = 0$ for $1 \leq i \leq k$,
$\Hat\eta_{k+1} = (\Tilde{\eta}_1, \sum e_i)$ and
$\Hat\kappa_{k+1} = \Tilde{\kappa}_1.$
We shall call the facet $\Tilde{F}_j'$ of $\Delta$ with  outward conormal 
$\Tilde{\eta}_j'$
the  {\bf fiber-type facet
associated to $\mathbf{ \Tilde{F}_j} $  }
and the facet
$\Hat{F}_i$ with outward conormal $\Hat\eta_i$ a  {\bf base-type facet}.
\end{defn}

\begin{figure}[htbp] 
   \centering
\includegraphics[width=4in]{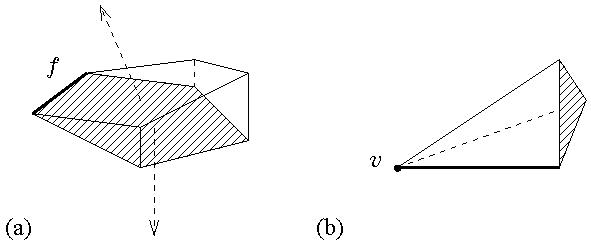} 
   \caption{(a) is the $1$-fold expansion of the shaded polygon along $f = \Tilde{F}_1$; the two dotted arrows are the conormals to the base-type facets. (b) is the $2$-fold expansion of the heavy line at the vertex $v= \Tilde{F}_1$. Its one fiber-type facet is dotted.}
   \label{fig:example}
\end{figure}

It is easy to check that $\De$ is simple, 
and is smooth 
exactly 
if $\Tilde\De$ is.  
The base-type facets intersect.  Indeed, the face
$ \bigcap_{i=1}^{k+1} \Hat{F}_i$
can be identified with $\Tilde F_1$. 
Similarly, 
the face
$\bigcap_{i \neq n} \Hat{F}_i$ can be identified with $\Tilde\Delta$
for all $1 \leq n \leq k+1$.   See Remark~\ref{rmk:geoexpand}
for a geometric interpretation of $1$-fold expansions.

\begin{example}\rm
The $k$-fold expansion of an $n$-simplex along
one of its facets is 
an $(n+k)$-simplex.
\end{example}

\begin{rmk}\labell{rmk:blowexp} \rm
For $\epsilon > 0$ sufficiently small, the
$\eps$-blow up $\De'$ of a polytope $\De$ along a face 
$f = F_I = \cap_{i \in I} F_i$
of codimension at least $2$ is
defined to be
 the intersection
$$
\Delta' =  \Delta \cap  
\{  x \in \ft^* \mid
\langle  \eta_0', x\rangle \leq  \kappa_0'  \}, 
$$
where
$\eta_0': = \sum_{i \in I} \eta_i$ and 
$\kappa_0' = \sum_{i \in I} \kappa_i - \epsilon$.
Thus $\Delta'$  
is obtained from $\De$ by cutting out a 
small neighborhood of $f$ by a suitably angled cut.\footnote
{
We give more details of this construction in Part II. 
}
Suppose now that $\De$ is the $k$-fold expansion of $\Tilde\De$  
 along $\Tilde F_1$, and consider its blow up along
 the face
$f = \bigcap_{i=1}^{k+1} \Hat{F}_i$. 
In this case, the conormal to the new facet $F_0'$ 
is $\eta_0'= 
(\Tilde\eta_1,0)\in \ft=\Tilde\ft \oplus \R^k$.  Therefore in the blow up there is a fiber type facet with conormal in $\Tilde\ft \oplus \{0\}$ corresponding to {\it each} of the facets of $\Tilde\De$.  
It is 
now
straightforward to check directly that $\Delta'$ is a
$\Tilde \Delta$ bundle over $\Delta_k$ and that the base facets
are $\Hat{F}_1 \cap \Delta',\ldots,\Hat{F}_{k+1} \cap \Delta'$; this justifies
our terminology.
\end{rmk}

The next result is an immediate consequence of Lemma~\ref{le:equiv}.

\begin{lemma} \labell{le:equiv0}
 The base facets of a bundle over
a simplex are equivalent, and the base-type facets of
a $k$-fold expansion are equivalent.
\end{lemma}

Conversely, 
the next 
result  shows that
if a
polytope has equivalent facets, it  is 
either
an expansion or a bundle over a simplex.

\begin{prop}\labell{prop:Iexpan} 
Let $\Delta \subset \ft^*$ be a simple  
(or smooth)
polytope. 
Let $I \in \Ii$ be an equivalence class of facets and 
define $I': = I \ssminus \{n\}$ for some $n \in I$.

\begin{itemize}\item[(i)] If $F_I  = \emptyset$, then $\Delta$ is a $F_{I'}$ bundle
over $\Delta_{|I|-1}$ 
with base facets  $\{F_i\}_{i \in I}$.

\item[(ii)] If $F_I \neq \emptyset$, then $\Delta$ is the
$(|I|-1)$-fold expansion of $F_{I'}$ along
 $F_I = F_n \cap F_{I'}$ 
with base-type facets $\{F_i\}_{i \in I}$.
\end{itemize}
\end{prop}

\begin{proof}{}
Renumber the facets so that $I = \{1,\ldots,|I|\}$ and $n = |I|$.
By Lemma~\ref{le:equiv},  the
plane $\Tilde\ft \subset \ft$ spanned by the $\eta_j$ for all $j > |I|$
has codimension $|I| - 1$.

Fix any $J \subset \{|I|+1, \ldots,N\}$
so that $F_J$ is not empty and is minimal, in
the sense that it has no lower dimensional
face $F_{J'} \varsubsetneq F_J$, 
where $J' \subset \{|I|+1, \ldots,N\}$.
We claim that $F_J$ is 
affine equivalent to the (standard) 
$(|I|-1)$-simplex with facets 
$F_i \cap F_J$ for $1 \leq i \leq |I|$; moreover, 
the $\eta_j$ for $j \in J$ span  $\Tilde \ft$.
To see this, note that  the only possible facets of $F_{J}$  
have the form $F_i \cap F_J$, where $1 \leq i \leq |I|$.  In particular,
$F_{J}$ has at most $|I|$ facets.
On the other hand, the plane
$P(F_{J})$ contains a translate of the annihilator of $\Tilde\ft$,
which implies that  $\dim (F_J) \geq |I| - 1$.
The claim follows.

By the paragraph above, 
$\sum_{1 \leq i \leq |I|} \eta_i$ lies in\footnote
{
This holds even in the nonsmooth case by our  conventions
 for normalizing the $\eta_i$; see Remark \ref{rmk:normalize} and Example \ref{ex:simplex}.}
$\Tilde \ft$
and $\{\eta_i\}_{1 \leq i < |I|}$
descends to a basis for 
the
 quotient of  $\ft$ by $\Tilde \ft$.
Hence, we can identify $\ft$ with $\Tilde \ft \oplus \R^{|I| -1}$ so that 
$\eta_i  = (0,-e_i)$  for all $1 \leq i < |I|$
and $\eta_{|I|} = (\alpha, \sum_{i =1}^{|I|-1} e_i)$ for some $\alpha \in \Tilde \ft$.
Moreover, let $J \subset \{ |I| + 1, \dots, N \}$  be any subset so
that the face $F_J$ is not empty.
By the paragraph above, the intersection
$F_J \cap F_{\Hat{I}}$ is also not empty for any proper subset  $\Hat{I} \subsetneq I$.
Hence, if $F_I = \emptyset$, then $\De$ is 
combinatorially equivalent to the product 
$F_{I'} \times \De_{|I|-1}$, and so $\De$ is a $F_{I'}$ bundle over $\De_{|I|-1}$.
On the other hand, if $F_I \ne \emptyset$
then
the outward conormal to the facet 
$F_I = F_{|I|} \cap F_{I'}$
of the polytope $F_{I'}$ is the image of $\eta_{|I|}$
in the quotient of $\ft$ by the span  of $\{\eta_i\}_{i \in I'}$.
This is $\alpha \in \Tilde\ft$, as required.
\end{proof}

\begin{cor}\labell{cor:Iexpan}
Let $\De$ be a simple polytope such that all its facets are equivalent.  Then $\De=\De_k$.
\end{cor}
\begin{proof}  In this case there is just one equivalence class of facets $I$.  Since $\De$ is bounded, $F_I=\emptyset$. Now apply part (i) of the above proposition.
\end{proof}

%%%%%%%%%%%%%%%%%%%%%%%%%%%%%%%%%%%%%%%%%%%%%%%%%%%%%%%%%%%
\subsection{Subtracting inessential functions}\labell{ss:bs}
%%%%%%%%%%%%%%%%%%%%%%%%%%%%%%%%%%%%%%%%%%%%%%%%%%%%%%%%%%

In this subsection, we show that it is often possible to use
Proposition~\ref{prop:inessential} to   
simplify a mass linear function by subtracting an inessential
function; the  resulting  function will have fewer
asymmetric facets.

\begin{lemma}\labell{le:ea}
Let $H \in \ft$ be a mass linear function on a simple polytope
$\Delta \subset \ft^*$. 
If $F_1,\ldots,F_m$ are equivalent facets,
there exists an inessential function $H' \in \ft$ 
so that the mass linear function
$\Tilde{H} = H - H'$  
has the following properties:
\begin{itemize}
\item 
For all $i < m$,
the facet $F_i$ is  $\Tilde{H}$-symmetric. 
\item   For all $i > m$,  
the facet $F_i$ is $\Tilde{H}$-symmetric iff it is $H$-symmetric. 
\end{itemize}
\end{lemma}

\begin{proof}{}
Since $H$ is mass linear,
$\langle H,c_\Delta\rangle = 
\sum_{i=1}^N \be_i \kappa_i$ for some $\be_i \in \R$.
Define $\al_i = \be_i$ for all $1 \leq i < m$ and 
define $\al_m = - \sum_{i = 1}^{m-1} \be_i$.
Then $\sum_{i=1}^m \al_i = 0$, and so by definition
$H' = \sum_{i=1}^m \al_i \eta_i$
is inessential. 
By Proposition~\ref{prop:inessential}, $\langle H',c_\Delta\rangle  = 
\sum_{i = 1}^m \al_i \kappa_i$.  
Therefore, if $\Tilde{H} = H - H'$ then
$$
\langle \Tilde H, c_\Delta \rangle  = \Bigl(\sum_{i=1}^m \beta_i\Bigr) \kappa_m
+ \sum_{i=m+1}^N \beta_i \kappa_i,
$$
as required.
\end{proof}

\begin{cor}\labell{cor:ea}
Let $H \in \ft$ be a mass linear function on a 
simple polytope $\Delta \subset \ft^*$.
If the asymmetric facets 
are all equivalent, then $H$  is inessential.
\end{cor}

\begin{proof}
By Lemma~\ref{le:ea}, there exists an inessential function
$H' \in \ft$ so that 
the mass linear function
$\Tilde{H} = H - H'$ 
has at most one asymmetric facet.
By Lemma~\ref{le:twoasym}, this implies that $\Tilde{H} = 0$.
\end{proof}

This method 
of simplification 
works particularly well when the equivalent facets
do not
intersect, that is, when our polytope is a bundle over
a simplex.  To prove this, we will need the following 
lemma.

\begin{lemma}\labell{le:prebund}  
Let $H \in \ft$ be a mass linear
function on a simple polytope $\De \subset \ft^*$ .
Assume that $\De$ is a
$\Tilde \De$ bundle over the simplex $\Delta_k$,
and let $\Tilde \ft \subset \ft$
be the subspace spanned by the conormals to the fiber facets.
Then $H$ lies in 
$\Tilde \ft$
exactly if the base facets are symmetric. 
\end{lemma}

\begin{proof}\,
Let $\Tilde\eta_1, \ldots, \Tilde\eta_{\Tilde{N}}$ be the outward
conormals to the fiber facets,  and let
$\Hat\eta_1,\ldots,\Hat\eta_{k+1}$ be the outward conormals to
the base facets.  Let $\Tilde\kappa_1,\ldots,\Tilde\kappa_{\Tilde{N}}$
and $\Hat\kappa_1,\ldots,\Hat\kappa_{k+1}$ be the associated
support numbers.  

If the base facets are  symmetric then 
$H$ lies in $\Tilde \ft$
by Lemma~\ref{le:Hsum}.
To prove the converse, 
assume that $H$ lies in  $\Tilde \ft$. 
Because $\Tilde \De$ is bounded,  
the positive span
of the $\Tilde \eta_i$ contains all of $\Tilde \ft$.
In particular, $H = \sum \be_i \Tilde \eta_i$, where $\be_i > 0$
for all $i$.
Because $\langle \Tilde\eta_j, x\rangle \leq \Tilde\kappa_j$
for all $x \in \Delta$, we have
$\langle 
\Tilde\eta_j, c_\Delta \rangle \leq \Tilde\kappa_j$ for all $j$.  Hence 
$\langle  H,c_\Delta \rangle \leq \sum \be_j \Tilde \kappa_j$.
The number $\sum \be_j \Tilde\kappa_j$  does not depend 
on the $\Hat\kappa_i$.
On the other hand,
because the base is a simplex,
we can 
make any $\Hat\kappa_i$  arbitrarily
large without changing  the combinatorics of the polytope $\De$.
Since $\langle H, c_\Delta \rangle$ is a linear function on $\Cc_\Delta$, 
this implies  that it 
cannot depend positively 
 on any $\Hat\kappa_i$.
The same  argument for $-H$ shows that 
$\langle H,c_\Delta \rangle$ cannot depend negatively 
on any $\Hat\kappa_i$. 
\end{proof}

\begin{prop}\labell{prop:bund}  
Let $H \in \ft$ be a mass linear function on a simple polytope
$\Delta \subset \ft^*$ which is
a bundle over the 
simplex $\De_k$. 
There exists  
an inessential function $H' \in \ft$ 
 so that 
the mass linear function $\Tilde{H} = H - H'$ has the following properties:
\begin{itemize}
\item The base facets are $\Tilde{H}$-symmetric
\item A fiber facet is $\Tilde{H}$-symmetric iff it is $H$-symmetric.
\end{itemize}
\end{prop}

\begin{proof}\,  
Let $\Hat\eta_1, \ldots,\Hat{\eta}_{k+1}$ denote 
the outward conormals to the base facets.
By Lemma~\ref{le:equiv} 
and the definition of a bundle, 
the base facets are equivalent.
Moreover,
we may decompose the vector $H\in \ft$ uniquely as 
$$
H= \Tilde{H} +  \sum_{i=1}^{k+1} \al_i\Hat\eta_i,  
\quad\mbox{where}\quad
 \sum \al_i = 0
\ \ \mbox{and}\ \
\Tilde H
\in \Tilde{\ft},
$$
where  $\Tilde{\ft} \subset \ft$ is the subspace  spanned
by the conormals of the fiber facets.
The function $H' = \sum_i \al_i\Hat\eta_i$ is 
inessential by definition.
Hence, by Proposition~\ref{prop:inessential}, $H'$ is mass linear and 
$\langle H',c_\De \rangle$ depends only on the support numbers of the 
base facets.  Thus
$\Tilde H$ is mass linear and
the fiber facets 
are $\Tilde{H}$-symmetric exactly
if they are $H$-symmetric.
The base facets are $\Tilde{H}$-symmetric by  Lemma~\ref{le:prebund}.
\end{proof}

%%%%%%%%%%%%%%%%%%%%%%%%%%%%%%%%%%%%%%%%%%%%%%%%%%%%%%%%%
\subsection{Flat facets}\labell{ss:flat}
%%%%%%%%%%%%%%%%%%%%%%%%%%%%%%%%%%%%%%%%%%%%%%%%%%%%%%%%%

In this subsection, we consider  
flat facets.
The main result (Proposition~\ref{prop:flat}) 
allows us to reduce 
our discussion of the structure of mass linear pairs $(\De,H)$
to the 
case that 
every asymmetric facet is pervasive.  

If $\Delta$ is a 
bundle over $\Delta_1$
then the base facets are flat and do not intersect.
The lemma  below proves  
the converse.

\begin{lemma}\labell{le:flata}
Let $\Delta \subset \ft^*$ be a simple polytope.
Let $F_1$ and $F_2$ be flat facets of $\Delta$ which do not intersect.
Then $\Delta$ is a bundle over the simplex $\Delta_1$ with base
facets  $F_1$ and $F_2$.
\end{lemma}

\begin{proof}
By definition,
since  $F_1$ is flat there exists
a nonzero
$\xi_1 \in \ft^*$
so that $\langle  \eta_j,\xi_1 \rangle = 0$  for every  facet
$F_j$ (other than $F_1$ itself) which meets $F_1$.
Since $F_1$ is bounded, 
these $\eta_j$ span $\ft/\eta_1$;
hence 
$\langle \eta_1, \xi_1 \rangle \ne 0$ and 
so
we can renormalize $\xi_1$ so that  $\langle \eta_1, \xi_1 \rangle = 1$. 
Similarly, there exists $\xi_2 \in \ft^*$ so that
$\langle \eta_2, \xi_2 \rangle = 1$ 
but $\langle  \eta_j,\xi_2 \rangle = 0$  for every (other) facet
$F_j$ which meets $F_2$.

We claim that
\begin{equation}\labell{flataeq}
\langle  \eta_k,\xi_1 \rangle \leq 0 \quad \forall \, k \neq 1 
\qquad \mbox{and} \qquad
\langle  \eta_k,\xi_2 \rangle \leq 0 \quad \forall \, k \neq 2. 
\end{equation}
To see this,
pick a point $x $ which lies in 
some
$F_k$  but not in $F_1$.
By definition $\langle  \eta_k, x  \rangle = \kappa_k$,
$\langle  \eta_1, x  \rangle < \kappa_1$,
and $\langle  \eta_i, x \rangle \leq  \kappa_i$  for all $i$.
Define 
$y: = x + \bigl(\kappa_1 - \langle  \eta_1,x \rangle \bigr) \xi_1$ .
An easy calculation shows that $\langle  \eta_1,y \rangle = \kappa_1$
and
$\langle  \eta_k,y \rangle = \kappa_k +
\bigl(\kappa_1 - \langle  \eta_1,x \rangle \bigr) \langle \eta_k, \xi_1
\rangle. $ 
If $F_i\ne F_1$ is any facet which meets $F_1$ 
then 
$\langle \eta_i, \xi_1 \rangle = 0$,
and hence
$\langle  \eta_i,y \rangle = 
\langle  \eta_i,x \rangle \leq \kappa_i$.
Therefore, $y \in F_1 \subset \Delta$.
In particular, $\langle  \eta_k,y \rangle \leq \kappa_k$
for all $k$.
Since $\langle \eta_1,  x \rangle <  \kappa_1$,
this implies that $\langle \eta_k, \xi_1 \rangle \leq 0$,
as required.

Next, we claim that $\xi_2$ is a  multiple of $\xi_1$. 
To see this, let 
$J$ denote the set of $j \in \{2,3, \ldots,N\}$ 
so that $F_j$ meets $F_1$. 
On the one hand, since $\langle \eta_j, \xi_1  \rangle =  0$ for all  $j \in J$
and since $F_1$ is bounded,  
the positive
span of the $\eta_j$ for $j \in J$ contains 
$\xi_1^\circ$, 
the annihilator of $\xi_1$.
On the other hand, since  $F_1 \cap F_2 = \emptyset$, $2$ is not in $J$.
Hence \eqref{flataeq} implies that
 $\langle \eta_j, \xi_2 \rangle \leq 0$ for all  $j \in J$.
Therefore,  $\langle \alpha, \xi_2 \rangle \leq 0$ for all
$\alpha \in \xi_1^\circ$.
Since $\xi_1^\circ$ is a linear subspace, this is possible only if 
 $\langle \alpha, \xi_2 \rangle = 0$ for all   
$\al \in \xi_1^\circ$.  
Thus $\xi_2$ is a multiple of $\xi_1$ as claimed.

Since $\langle \eta_2, \xi_2 \rangle = 1$
and $\langle \eta_1,\xi_2 \rangle \leq 0$ by \eqref{flataeq},
$\xi_2$ cannot be a positive multiple of $\xi_1$; so it
must be a negative multiple.
Hence, for all $k$ except $1$ and $2$,
we have $\langle \eta_k ,\xi_1 \rangle \leq 0$ and $\langle \eta_k, \xi_1 \rangle \geq 0$,
that is,  $\langle \eta_k, \xi_1  \rangle =  0$.

To see that $F_1$ and $F_2$  
are the base facets of a bundle over $\De_1$ it now suffices to check that 
$\De$ is combinatorially equivalent to
the product $F_1\times \De_1$.  Equivalently, we must check that every 
non-empty face $F_J$, where $J \subset \{3,\ldots,N\}$,
meets both $F_1$ and $F_2$. 
To see this, 
pick $x \in F_J$ and 
define $y:= x + (\kappa_1 - \langle \eta_1, x \rangle) \xi_1$
as before.
Since $\langle \eta_1, \xi_1 \rangle = 1$, $\langle \eta_2, \xi_1 \rangle < 0$,
and otherwise $\langle \eta_j, \xi_1 \rangle = 0$,
$y \in F_1 \cap F_J$.
A similar argument shows that $F_2 \cap F_J \neq \emptyset$.
\end{proof}

\begin{cor}\labell{cor:newflat}
Let $H \in \ft$ be a mass linear function
on a simple polytope  $\Delta \subset \ft^*$. 
If $F$ is an asymmetric facet which is not pervasive,
then $\Delta$ is an $F$ bundle over 
$\Delta_1$. 
\end{cor}

\begin{proof}
Let $G$ be any facet which is disjoint from $F$.
By Lemma~\ref{le:symasym}, 
$G$ is also  asymmetric.
Hence, Proposition~\ref{prop:asym} implies that 
$F$ and $G$ are flat.
By Lemma~\ref{le:flata}, 
this implies that
$\Delta$ is an $F$ bundle over $\Delta_1$.
\end{proof}

Here is the main result of this section.

\begin{prop}\labell{prop:flat}
Let $H \in \ft$ be a mass linear function
on a simple  polytope $\Delta \subset \ft^*$.
There exists an inessential function $H' \in \ft$ so that 
the mass linear function $\Tilde H = H - H'$ has the following property:
every $\Tilde H$-asymmetric
facet is pervasive.
\end{prop}

\begin{proof}\,
We will argue this by 
induction on the number 
of asymmetric facets which are not pervasive.
Let $F$ be an $H$-asymmetric facet which is not pervasive.
By Corollary~\ref{cor:newflat}, $\Delta$ is an $F$ bundle over $\De_1$.
By Proposition~\ref{prop:bund} 
there exists an inessential function $H' \in \ft$ so that
the mass linear function    
$\Tilde{H} =  H - H'$ 
has the following properties:
the base facets (including  $F$) are $\Tilde{H}$-symmetric, and 
the fiber facets are $\Tilde H$-symmetric exactly if they are $H$-symmetric. 
Hence, there are fewer $\Tilde{H}$-asymmetric 
facets which are not pervasive.
Since the sum of inessential functions is inessential,
this completes the proof.
\end{proof}

%%%%%%%%%%%%%%%%%%%%%%%%%%%%%%%%%%%%%%%%%%%%%%%%%%%%%%%%%%%%%
\subsection{Polytopes such that all linear functions are mass linear}\labell{ss:allmass}
%%%%%%%%%%%%%%%%%%%%%%%%%%%%%%%%%%%%%%%%%%%%%%%%%%%%%%%%%%%%%

This subsection is devoted to a proof of
Theorem~\ref{thm:allmass}.
First, note that every inessential function is mass linear by Proposition~\ref{prop:inessential}.
Next, 
note that by 
Example~\ref{ex:square},
if $\Delta \subset \ft^*$ is a product of simplices, then
every $H \in \ft$ is an inessential mass linear function on $\Delta$,
and  $\bigcap_{i \in I} F_i = \emptyset$ for every equivalence class
of facets $I$.
Thus  (iii) $\Longrightarrow $  (ii) $\Longrightarrow $ (i) and (iii) $\Longrightarrow$ (iv).

To complete the proof, we
will repeatedly need the following elementary lemmas.

\begin{lemma} \labell{le:etaasym}
Let $\Delta \subset \ft^*$ be a simple polytope.
Each
facet $F_i$ is $\eta_i$-asymmetric.
\end{lemma}

\begin{proof}
Assume that $F_i$ is $\eta_i$-symmetric.  By Lemma~\ref{le:symcent},
this implies that $\langle \eta_i, c_{F_i} \rangle = \langle \eta_i, c_\Delta \rangle$.
Since $\langle \eta_i, c_{F_i} \rangle = \kappa_i$, this
contradicts the assumption  that $F_i$ is $\eta_i$-symmetric.
\end{proof}

\begin{lemma}\labell{le:oldIII}
Let $F_i$ be a pervasive facet of a simple polytope $\Delta \subset \ft^*$,
and let $j,k\ne i$.
Then $F_{ji} = F_j \cap F_i$ and $F_{ki} = F_k \cap F_i$ 
are equivalent as facets of $F_i$ exactly if
 $F_j$ and $F_k$ are equivalent as facets of $\Delta$.
\end{lemma}

\begin{proof}\,
By Lemma~\ref{le:equiv},
$F_j\sim F_k$ exactly if  there is a vector $\xi \in \ft^*$ 
that lies in the annihilator ${\eta_\ell}^\circ$ of the conormal $\eta_\ell$ for all $\ell\ne j,k$. 

Since $F_i$ is pervasive, $F_{\ell i}$ is a facet of $F_i$ for every 
$\ell  \neq i$.  Moreover, $P(F_i)$ is naturally
affine isomorphic to ${\eta_i}^\circ $.  
Therefore, $F_{ji}$ and $F_{ki}$ are facets of $F_i$,
and they are equivalent facets exactly if there exists $\xi \in \eta_i^\circ \subset \ft^*$
that lies in the annihilator of ${\eta_\ell}^\circ$ for all $\ell \ne i,j,k$.

But these two conditions are clearly equivalent.
\end{proof}

\NI
{\bf Proof that (i) $\Longrightarrow$ (iii) in Theorem~\ref{thm:allmass}.}

Let $\Delta \subset \ft^*$ be a simple polytope;
let $\Ii$ denote the set of equivalence classes of facets.
Assume that every $H \in \ft$ is mass linear on $\Delta$.
We will prove that $\Delta$ is 
a product of simplices.
The proof is divided into a sequence of steps.
\MS

\NI {\bf Claim I:} 
{\it If $\dim \Delta \leq 2$, then $\Delta$ is 
a product of simplices.}
\SSS

The claim is  trivial in dimension $1$.
So  assume that  $\dim \Delta = 2 $  and that
$\Delta \neq \Delta_2$. 
Since $\Delta$ has more than three edges,  no edge is pervasive.
Moreover,   each outward conormal $\eta_i$ is mass linear by assumption,
and  each edge $F_i$ is $\eta_i$-asymmetric 
by Lemma~\ref{le:etaasym}.
Therefore, Proposition~\ref{prop:asym} implies that 
each edge $F_i$ of $\De$ is flat.  But this is possible only if
$\Delta$ is the product $\Delta_1 \times \Delta_1$.
\MS

We will now assume that $\dim \Delta \geq 3$ and that 
(i) $\Longrightarrow$ (iii)
for all lower dimensional polytopes.

\MS

\NI {\bf Claim  II:} 
{\it Each proper face of $\Delta$ is 
a product of simplices.}
\SSS

Since we have assumed that 
(i) $\Longrightarrow$ (iii) holds 
for lower dimensional polytopes,
it is enough to prove every facet $F_i$ of $\Delta$
satisfies
condition (i);
that is, we
must show that every $H'$ in 
$\ft/ \left< \eta_i \right>$ -- which is naturally the dual space to $P(F_i)$ --  is mass linear on $F_i$. 
To see this,
let $\Hh_i$ denote the set of $H \in \ft$ so that the facet $F_i$ is $H$-symmetric.
By Lemma \ref{le:etaasym},
$\eta_i \not\in \Hh_i$.
On the other hand, since every $H \in \ft$ is mass linear on $\Delta$,
there is $\la\in \R$ such that the coefficient of $\ka_i$ in
$H-\la\eta_i$ is zero, that is, $H-\la\eta_i \in \Hh_i$
Hence $\Hh_i \subset \ft$ is a linear subspace of codimension 
one and so  the  natural projection from  $\Hh_i$ to  $\ft/ \left< \eta_i \right>$ is surjective. 
Finally, by 
Proposition~\ref{prop:symface},
 each $H \in \Hh_i$ restricts to a mass linear function
on $F_i$.  
\MS

\NI{\bf Claim III:} {\it If $F_1$ is not pervasive, $\Delta$ is a bundle
over the $1$-simplex $\Delta_1$.}

\MS
Since $\eta_1 \in \ft$ is a mass linear function on $\Delta$
and $F_1$ is $\eta_1$-asymmetric, $\Delta$ is a bundle
over $\Delta_1$ by Corollary~\ref{cor:newflat}.

\MS

\NI {\bf Claim IV:}
{\it If every facet of $\Delta$ is pervasive, 
$\bigcap_{i \in I} F_i = \emptyset$ for all $I \in \Ii$.}
\SSS

Let $\{F_1,\ldots,F_{k+1}\}$ be an equivalence class of facets.
If $\De$ has no other facets, then the conclusion holds because
$\De$ is bounded. (Alternatively,  
$\De=\De_k$ by Corollary~\ref{cor:Iexpan}.) Otherwise, fix
$j > k+1$.  Since $F_j$ is pervasive,
$\{F_{ij}\}_{i=1}^{k+1}$ is a nonempty  equivalence class of facets
of $F_j$ by Lemma~\ref{le:oldIII}.  
Moreover, by  Claim  II, $F_j$ is the product of simplices.
Hence, 
by Example~\ref{ex:square},
$\bigcap_{i=1}^{k+1} F_{ij} = \emptyset.$
In particular, $k > 0$.
Therefore, 
a nearly identical argument shows that
$\bigcap_{i=2}^{k+1} F_{i1}  = \emptyset$.
But  $\bigcap_{i=2}^{k+1} F_{i1} =\bigcap_{i=1}^{k+1} F_{i}$.

\MS

\NI {\bf Claim V:} 
{\it $\Delta$ is a bundle over the $k$-simplex
$\Delta_k$, where $k \neq \dim \Delta - 1$.}
\SSS

If some facet of $\Delta$ is not pervasive, then
$\Delta$ is a bundle over $\Delta_1$ by Claim III. 
Since $\dim \Delta \geq 3$, $\dim \Delta - 1 \neq 1$.
So assume instead that every facet is pervasive.
By Claim IV,  $\cap_{i \in I} F_{i} = \emptyset$ for all $I \in \Ii$.
Therefore, applying
Proposition~\ref{prop:Iexpan}, 
we see
that $\Delta$ is a bundle over the $k$-simplex $\Delta_k$.
Finally, if $k = \dim \Delta - 1$ then the fiber is $1$-dimensional;
hence it is $\Delta_1$.  This contradicts the claim that every
facet is pervasive.
\MS

\NI {\bf Claim VI:}  
{\it $\Delta$ is a product of simplices.}
\SSS

By Claim V, $\Delta$ is a bundle over the $k$-simplex $\Delta_k$,
where $k \neq \dim \Delta - 1$.
Label the base facets  $\{F_1,\ldots,F_{k+1}\}$.
If $k = \dim \Delta$ we are done, so assume that 
$k \leq \dim \Delta  - 2$.

Since $k+1 < \dim \Delta$ and $\Delta$ is bounded, 
there exists 
$i > k+1$
so that $\eta_i$  does not lie in 
the span of 
$\eta_1,\ldots,\eta_{k+1}$. 
Since $\Delta$ is a bundle over $\Delta_k$ with base facets $F_1,\dots,F_{k+1}$,
the facet
$F_i$ is also a bundle over $\Delta_k$
with base facets $F_1 \cap F_i,\ldots,F_{k+1} \cap F_i$.
Moreover, identifying the dual space to $P(F_i)$  with $\ft / \left< \eta_i \right>$ as in Claim II,
we see that
the projections 
of $\eta_1,\ldots,\eta_{k+1}$ to $\ft / \left< \eta_i \right>$
are the conormals to the base facets 
$F_1 \cap F_i,\ldots,F_{k+1} \cap F_i$ of $F_i$.  But the bundle is trivial by Claim II.  Hence these conormals
are linearly dependent 
in the quotient space.  
Since $\eta_i$ does not lie in the span of
$\eta_1,\ldots,\eta_{k+1}$, this implies 
that $\eta_1,\ldots,\eta_{k+1}$  are themselves linearly dependent,
that is, $\Delta$ is a trivial bundle over $\Delta_k$.
Finally, the fiber of this bundle
is naturally isomorphic to the face $\cap_{i =2}^{k+1} F_i$.
By Claim II, this face is 
a product of simplices.
\QED\MS

\NI {\bf Proof that (iv) $\Longrightarrow$ (iii) in Theorem~\ref{thm:allmass}.}

Let $\Delta \subset \ft^*$ be a simple polytope; let $\Ii$
denote the set of equivalence classes of facets.
Assume that $\bigcap_{i \in I} F_i = \emptyset$ for
all $I \in \Ii$; in particular, $|I| > 1$.
We will  prove that $\Delta$ is a product of
simplices.  As before, 
we
divide the proof into a sequence of steps.

\MS
\NI {\bf Claim I:} 
{\it If $\dim \Delta \leq 2$, then $\Delta$ is 
a product of simplices.}
\MS

If $\dim \Delta = 2$ 
and there is an equivalence class with 
more than two
elements, 
or if there is an equivalence class with exactly two edges which meet, 
then Proposition~\ref{prop:Iexpan} implies that
$\De$ is $\De_2$.
Otherwise $\De$ is a $\De_1$-bundle over $\De_1$ such that
its two fiber facets are equivalent. 
By Lemma~\ref{le:equiv},  this implies that $\De$ is a product.
\MS

We will now assume that
(iv) $\Longrightarrow$ (iii) 
for all lower dimensional polytopes.

\MS
\NI{ \bf Claim II:}
{\it Every proper face of $\Delta$ is a product of simplices.}
\MS

Since we have assumed that (iv) $\Longrightarrow$ (iii) holds for
lower dimensional polytopes, it is enough to prove that every
facet $F_i$ satisfies
condition (iv); in other words,
 we must show that $\bigcap_{i \in I} F_{i1} = \emptyset$
for every equivalence class $\{F_{i1}\}_{i \in I}$ of facets of  $F_1$.
Even if $F_1$ is not pervasive, it 
is clear from the definitions 
that if $F_i$ and $F_j$
are equivalent facets of $\Delta$
that intersect $F_1$
then $F_{i1}$ and $F_{j1}$
are equivalent facets of $\Delta$.
Therefore, 
for each equivalence class $I$ of facets of $F_1$, the set $\{F_1\} \cup \{F_i\}_{i \in I} $ contains
an equivalence class of facets of $\De$.
By assumption, this
implies that  $F_1 \cap \bigcap_{i \in I} F_i  = \emptyset$.
But $F_1 \cap \bigcap_{i \in I} F_i = \bigcap_{i \in I} F_{i1}.$
\MS

\NI {\bf Claim III:} {\it If $F_1$ is not pervasive, $\Delta$
is a bundle over the $1$-simplex $\Delta_1$.}
\MS

Let $I \in \Ii$ be the equivalence class of facets which
contains $F_1$.  By assumption, $\bigcap_{i \in I} F_i =
\emptyset$.   By Proposition~\ref{prop:Iexpan}, this
is only possible if $\Delta$ is a bundle over $\Delta_1$.

\MS
The rest of the proof now follows exactly as in the previous
case, although  in fact Claim IV is true by assumption.
\QED

%%%%%%%%%%%%%%%%%%%%%%%%%%%%%%%%%%%%%%%%%%%%%%%%%%%%%%%%%%%%%%%%
\section{Polytopes in $2$ and $3$ dimensions}\labell{s:23dim}
%%%%%%%%%%%%%%%%%%%%%%%%%%%%%%%%%%%%%%%%%%%%%%%%%%%%%%%%%%%%%%%%

In this section, we give a complete classification of polygons
and smooth $3$-dimensional polytopes which admit mass linear functions.
Here, a {\bf (convex) polygon} is a simple $2$-dimensional polytope.

%%%%%%%%%%%%%%%%%%%%%%%%%%%%%%%%%%%%%%%%%%%%%%%%%%%%%%%%%%%%%%%%
\subsection{Polygons}\labell{ss:2dim}
%%%%%%%%%%%%%%%%%%%%%%%%%%%%%%%%%%%%%%%%%%%%%%%%%%%%%%%%%%%%%%%%%%%

In this subsection, we describe  mass linear functions on polygons,
showing that they are all inessential.

We shall need the following lemma; since
every edge has only two facets, it is an immediate consequence
of  Lemmas~\ref{le:symasym} and \ref{le:twoasym}. 

\begin{lemma}\labell{le:symedge}
Fix  
a nonzero
$H \in \ft$ and  a simple
polytope $\Delta \subset \ft^*$.
If $\Delta$ has a symmetric edge then it has 
exactly 
 two asymmetric facets.
\end{lemma}

Here is our main result;
it is an elaboration 
of Theorem~\ref{thm:2}.

\begin{prop}\labell{prop:2dim}  
Let $H \in \ft$ be a nonzero mass linear function on a 
polygon $\De \subset \ft^*$.
Then one of the following statements holds:
\begin{itemize}
\item $\Delta$ is a triangle;
at most one edge is symmetric.
\item $\Delta$ is a $\Delta_1$ bundle over $\Delta_1$; 
the base facets are the asymmetric edges.
\item $\Delta$ is the product $\Delta_1 \times \Delta_1$; 
each edge is asymmetric.
\end{itemize} 
In any case,
$H$ is inessential.
\end{prop}

\begin{proof}  
Assume that $\Delta$ has a symmetric edge.
By Lemma~\ref{le:symedge},
$\Delta$ has two asymmetric edges.
Moreover,  by  Lemma~\ref{le:symasym},
every symmetric edge must meet both  asymmetric edges;
in particular, there are at most two symmetric edges.
If there is one symmetric edge, $\Delta$ is a triangle.
If there are two symmetric edges, then the asymmetric edges are not pervasive
and so by Proposition~\ref{prop:asym} they are flat.
This implies that $\Delta$ is a $\Delta_1$ bundle over $\Delta_1$,
and that the fiber facets are symmetric.
In either case, the asymmetric facets are equivalent.
Hence, $H$ is inessential by Corollary~\ref{cor:ea}.

So assume that every edge is asymmetric.
If $\Delta$ has three edges, it is a triangle.
Otherwise, since none of the edges are pervasive
Proposition~\ref{prop:asym} implies that they are all flat. 
By definition, this is impossible unless $\Delta$ is the product $\Delta_1 \times \Delta_1$.
Moreover, all the facets of
a triangle
are equivalent,
and the two pairs of  opposite  facets of $\Delta_1 \times \Delta_1$ are equivalent.
It is easy to check
that in either 
case 
every $H' \in \ft$ is inessential.
\end{proof}

The next corollary will be useful in Part II.

\begin{cor}\labell{cor:2dim}
Let $H \in \ft$ be a 
mass linear function on a 
polygon $\De \subset \ft^*$.
If two edges $F_i$ and $F_j$ do not intersect
then $\gamma_i + \gamma_j = 0$, where $\gamma_k$ is the coefficient
of the support number of $F_k$ in the linear function $\langle H, c_\Delta \rangle$.
\end{cor}

\begin{proof} 
If $H \neq 0$, then 
the lemma above implies that $H$ is inessential
and $F_i$ and $F_j$ 
are opposite edges of  a quadrilateral. 
Hence while they may (or may not) be equivalent to each other, 
neither is equivalent to any other edge. 
Therefore Proposition~\ref{prop:inessential} implies that either
$\ga_i=\ga_j=0$ (if $F_i\not\sim F_j$) or $\ga_i + \ga_j = 0$.
\end{proof}

%%%%%%%%%%%%%%%%%%%%%%%%%%%%%%%%%%%%%
\subsection{Polytopes with two asymmetric facets}
%%%%%%%%%%%%%%%%%%%%%%%%%%%%%%%%%%%%%

We now show that any mass linear function on a simple polytope
with at most two asymmetric
facets is inessential.  In this paper,
we need this result for the $3$-dimensional
case, but it is valid in all dimensions.

\begin{prop}\labell{prop:2asym}
Let $H \in \ft$ be a mass linear function on a simple 
polytope $\Delta \subset \ft^*$ with exactly two asymmetric facets.
Then the asymmetric facets are equivalent. 
\end{prop}

\begin{proof}{}
Let $F_1$ and $F_2$ be the asymmetric facets.

First consider a $2$-dimensional symmetric face   
$Y$ with two symmetric edges $e$ and $e'$. By  
Proposition~\ref{prop:symface}, $Y$  has exactly
two asymmetric edges $F_1 \cap Y$ and $F_2 \cap Y$,
and the restriction of $H$ to $Y$ is mass linear. 
Therefore by Proposition~\ref{prop:2dim}, 
$Y$ is a $\De_1$ bundle over $\De_1$ 
and $e$ and $e'$ are fiber facets.    Hence they are parallel.

More generally, define a graph $\Ga$ as follows:
let $V: = V(\Ga)$ be the set of vertices in $F_1 \ssminus (F_1 \cap F_2)$
and let $E: = E(\Ga)$ be the set of edges which join such vertices.
By  Lemma~\ref{le:symasym} each symmetric edge intersects $F_1$
and $F_2$ at two different points; hence,
it must intersect $F_1$ in a point of $V$.
On the other hand,  given any $v \in V$,
since $F_1$ is the only asymmetric facet through $v$
and since $\Delta$ is simple, $v$ must lie on a unique symmetric edge.
Hence,
intersection with $F_1$ induces 
a one-to-one correspondence between the set of
symmetric edges and $V$.
By a similar argument, intersection with $F_1$  also a one-to-one correspondence
between  the set of $2$-dimensional symmetric faces
$Y$
 with
two symmetric edges and $E$. 
(Here, we use the fact 
proved above that the two symmetric edges 
in $Y$ are parallel and so cannot
meet.)
Hence, if 
$e$ and $e'$ are symmetric edges so that
$e  \cap F_1$ and 
$e' \cap F_1$ lie in the same component
of $\Gamma$, then $e$ and $e'$ are parallel.
Since $\Ga$ is clearly connected, this implies that
all symmetric edges are parallel.

Every symmetric facet $G$ contains
a symmetric face $g$ which is minimal in 
the sense that it has  no symmetric faces. By 
Proposition~\ref{prop:symface}, 
$g$ is a polytope
with exactly two facets  $F_1 \cap g$ and $F_2 \cap g$.
This is only possible if $g$ is a symmetric edge.
Therefore, the  conormals to the symmetric facets lie
in a codimension $1$ subspace.
By Lemma~\ref{le:equiv}, this implies that the asymmetric 
facets  are equivalent.
\end{proof}

This proposition has the following important consequence.

\begin{corollary}\labell{cor:2asym}
Let $H \in \ft$ be a mass linear function on a simple 
polytope $\Delta \subset \ft^*$ with
exactly two asymmetric facets $F_1$ and $F_2$. 
Then $H$ is inessential, and
exactly one of the following occurs.
\begin{itemize}
\item
$\De$ is a $F_1$ bundle over $\De_1$; 
the base facets are the asymmetric facets, or
\item $\De$ is the expansion of $F_1$ along $F_{12}$;
the base-type facets are the asymmetric facets.
\end{itemize}
\end{corollary}

\begin{proof} 
By Proposition~\ref{prop:2asym}, $F_1$ and $F_2$  are equivalent.
Thus the claim follows immediately from
Corollary~\ref{cor:ea} and
Proposition~\ref{prop:Iexpan}.
\end{proof}

%%%%%%%%%%%%%%%%%%%%%%%%%%%%%%%%%%%%%%%%%%%%%%%%%%%%%%%%%%
\subsection{Mass linear functions on the polytopes defined in Example~\ref{ex:1}}\label{ss:calcul}
%%%%%%%%%%%%%%%%%%%%%%%%%%%%%%%%%%%%%%%%%%%%%%%%%%%%%%%%%%

In this subsection, we calculate all mass linear functions
on the polytopes defined in
Example~\ref{ex:1}. 
In particular, we show that these functions are not all inessential.
Let's recall this example.
Given $a = (a_1,a_2) \in \R^2$,
define
\begin{gather*}
\eta_1 = -e_1, \  \eta_2 = -e_2, \  \eta_3 = e_1 + e_2, \ 
\eta_4 = - e_3, \ \ \eta_5 = e_3 + a_1 e_1 + a_2 e_2, \ \  \mbox{and} \\
\Cc_a = \left\{ \kappa \in \R^5 \left| \
 \sum_{i=1}^3 \kappa_i > 0 \mbox{ and }
\kappa_4 + \kappa_5 >  - a_1 \kappa_1  - a_2 \kappa_2   
+  \max ( 0,  a_1 ,a_2) \sum_{i=1}^3 \kappa_i \right. \right\}.
\end{gather*}
Given $\kappa \in \Cc_a$, let 
\vspace{-.1in}
$$
Y = Y(\kappa) = \bigcap_{i = 1}^5 \{ x \in (\R^3)^* \mid \langle \eta_i, x\rangle 
\leq \kappa_i  \}.
$$
Here is 
the main result of this section.

\begin{proposition}\labell{prop:Mab}
Let $Y$ be the polytope defined in Example~\ref{ex:1}.
Then $H \in \ft$ is mass linear function on $Y$ exactly if
$$
H = \sum_{i=1}^5 \ga_i \eta_i,
\ \  \mbox{where} \ \ 
\ga_1+\ga_2+\ga_3=0,\ \   \ga_4 +  \ga_5 = 0,  
\  \mbox{and} \ \ a_1 \ga_1 + a_2 \ga_2 = 0. 
$$
In this case
$\langle H, c_Y \rangle = \sum_{i=1}^5 \ga_i \kappa_i;$ moreover, 
$H$ is inessential exactly if  
at least one of the following conditions holds:
$a_1  a_2 (a_1 - a_2) = 0$,
or $\gamma_1 = \gamma_2 = \gamma_3 = 0.$
\end{proposition}

\begin{rmk}\labell{rmk:Mab}\rm
In particular, all vectors $H \in \ft$  
are mass linear on $Y$ (and are inessential) exactly if 
$a_1=a_2=0$ and  $Y$ is the product $\De_1\times \De_2$. 
This is consistent with  Theorem~\ref{thm:allmass}.
Otherwise,
$\Delta$ admits a $2$-dimensional family of
mass linear functions.
If $a_1 a_2 (a_1 - a_2) = 0$, these functions are all inessential.
Otherwise, there is a $1$-dimensional family
of inessential functions; the rest are essential.
\end{rmk}

The proof of this proposition  
rests mainly on the
following direct calculation.

\begin{lemma}\labell{le:Mab}    
Let $Y$ be the polytope defined in Example~\ref{ex:1}, and
let $H = \ga_1\eta_1+\ga_2\eta_2+\ga_3\eta_3$, where 
$\ga_1+\ga_2+\ga_3=0$. 
Then $H$ is mass linear on $Y$ 
if and only if 
$$
\ga_1a_1+\ga_2a_2 = 0; \quad \mbox{in this case} \quad
\langle H, c_Y \rangle = \sum_{i=1}^3 \ga_i \kappa_i.$$
\end{lemma}

\begin{proof}  
As a first step, fix $\kappa_1 = \kappa_2 = \kappa_4 = 0$,
and let $\kappa_3 = \lambda$ and $\kappa_5 = h$.
Let 
$\Delta_2^\la \subset \R^2$ 
denote the  $2$-simplex given by the inequalities
$$
 x_1 \geq 0, \quad x_2 \geq 0, \quad 
\mbox{and} \quad x_1+x_2 \leq  \lambda.
$$
An elementary calculation shows that for any non-negative 
integers $i_1$ and $i_2$ 
$$ 
\int_{\Delta_2^\la} x_1^{i_1} x_2^{i_2}  = 
\frac{ i_1!\, i_2! \, \lambda^{I + 2}}{  (I + 2)!},$$
where $I  = i_1+i_2$ and 
where by convention $0! = 1$. Further, both here and below we integrate with respect to the standard measure $dx_1 dx_2$ on $\R^2$.

Since $Y$ is a $\Delta^\la_2$ bundle over $\Delta_1$,
 $Y$ has volume 
$$
V = \int_{\Delta_2^\la} \bigl( h - a_1 x_1 - a_2 x_2  \bigr)=
 \frac{3 h \lambda^2 -  (a_1+a_2) \lambda^{3} }{3!}.
$$
For $j \neq 3$,
the moment 
$\mu_j$ of $Y$ 
along the $x_j$ axis is
$$
\mu_j= \int_{\Delta_2^\la} \bigl((h - a_1 x_1 - a_2 x_2) x_j\bigr)
= \frac{4 h \lambda^{3} - 
 (a_j+ a_1+a_2) \lambda^{4} }{4!}.
 $$
Let $c_j :=
\mu_j/V$  denote the $j$'th component of the center of mass.
For $j \neq 3$, 
$$
c_j = 
\frac{\lambda}{4}\  \frac{4h  - 
\lambda (a_j + a_1+a_2)}
{3h - \lambda (a_1+a_2)}.
$$
Since $ \gamma_1+\ga_2+\ga_3 = 0$, 
a straightforward calculation shows that 
$$ 
\langle H, c_{Y} \rangle = \sum_{i=1}^2 
(\gamma_{3} - \gamma_i) c_i
= \lambda \left(\gamma_{3}  +
\frac{\lambda(\gamma_1 a_1+\gamma_2 a_2)}
{12 h  -4 \lambda (a_1+a_2)}\right).
$$
This is a linear function of $h$ and $\lambda$ exactly
if $ \gamma_1 a_1+\gamma_2 a_2 = 0$.    
Hence if $H$ is mass linear, this condition must be satisfied.

To prove the converse, assume that $\gamma_1 a_1 + \gamma_2 a_2 = 0$.
Given 
$\kappa \in \Cc_a$,
note that by Equation~\eqref{shift} 
in Lemma~\ref{le:Hsum},
\begin{gather*}
Y(\kappa) = Y(0,0,\lambda, 0,h) -
(\kappa_1, \kappa_2, \kappa_4), \\  \mbox{ where } 
\lambda = \kappa_1 + \kappa_2 + \kappa_3 \mbox{ and } 
h = \kappa_4 + \kappa_5 + a_1 \ka_1 + a_2 \ka_2.
\end{gather*}
Hence,
$$
\langle H, c_Y(\kappa) \rangle = \langle H, c_Y(0,0,\lambda,0,h) \rangle - 
\langle H, (\kappa_1,\kappa_2,\kappa_4) \rangle = \sum_{i=1}^3 \kappa_i \gamma_i,$$
as required.
\end{proof}

\begin{rmk}\rm\labell{rmk:easier}
Lemma~\ref{le:Mab} implies the following fact.
If  $H = \gamma_1 \eta_1 + \gamma_2 \eta_2 + \gamma_3 \eta_3$
is mass linear on $Y$, where $\gamma_1 + \gamma_2 + \gamma_3 = 0$,
then $\langle H, c_Y \rangle = \sum_{i=1}^3 \gamma_i \kappa_i$.
In fact, there is an  
easier proof  of this fact.
As we saw in Example~\ref{ex:2}, $\Delta$ is a 
$\Delta_2$ bundle over $\Delta_1$.
By Lemma~\ref{le:prebund}, the base facets $F_4$ and $F_5$ are symmetric.
By Proposition~\ref{prop:symface}
the restriction of $H$ to  $F_4$ is mass linear and the 
coefficient of the support number of $F_i$ in $\langle H, c_Y \rangle$
is the coefficient of the support number of $F_i \cap F_4$
in $\langle H, c_{F_4} \rangle $ for all $1 \leq i \leq 3$.
But since $F_4$ is the $2$-simplex,  every $H$ is inessential.
The claim now follows from Proposition~\ref{prop:inessential}.
\end{rmk}

We are now ready to prove our main proposition.
\MS

\NI {\bf Proof of Proposition~\ref{prop:Mab}.}
First, note
that  every $H \in \ft$ can be written uniquely as
$H = \sum_{i=1}^5 \gamma_i \eta_i $, where 
$\gamma_1 + \gamma_2 + \gamma_3 = 0$ and $\gamma_4 + \gamma_5 = 0$.
(This holds because $\eta_1+\eta_2+\eta_3=0$ and $\eta_4+\eta_5$ is in the span of 
$\eta_1$ and $\eta_2$.)
By Lemma~\ref{le:equiv} (or by Example~\ref{ex:1b}), 
$F_4$ and $F_5$ are equivalent, and so
$\gamma_4 \eta_4 + \gamma_5 \eta_5$ is inessential and
$$
\langle \gamma_4 \eta_4 + \gamma_5 \eta_5 , c_Y \rangle
= \gamma_4 \kappa_4 + \gamma_5 \kappa_5
$$
by Proposition~\ref{prop:inessential}. 
On the other hand, by Lemma~\ref{le:Mab},
$\Tilde H := \ga_1 \eta_1 + \ga_2 \eta_2 + \ga_3 \eta_3$
is mass linear exactly if $a_1 \ga_1 + a_2 \ga_2 = 0$, in
which case
$$
\langle  \Tilde H, c_Y \rangle
= \ga_1 \kappa_1 + \ga_2 \kappa_2 + \ga_3 \kappa_3.
$$ 
The first claim follows immediately.

To prove the last statement, 
it is convenient to introduce a dummy variable 
$a_3$ that we set equal to $0$. 
Then, by Lemma~\ref{le:equiv},
for any $\{i,j\} \subset  \{1,2,3\}$,
$F_i$ is equivalent to $ F_j$ exactly if $a_i = a_j$.
Consider a mass linear 
function
$H = \sum_{i=1}^5 \gamma_i \eta_i $, where 
$\gamma_1 + \gamma_2 + \gamma_3 = 0$ and $\gamma_4 + \gamma_5 = 0$
as above. 
Since
$\ga_4 \eta_4 + \ga_5 \eta_5$
is inessential, 
$H$ is inessential exactly  if
$\Tilde H: = 
 \ga_1 \eta_1 + \ga_2 \eta_2 + \ga_3 \eta_3$ is inessential.   But 
  if $a_1 = a_2$, for example, then
we must have
 $\ga_1 + \ga_2 = 0$ and hence also $\ga_3=0$ so that
 $\Tilde H$ is inessential.
The cases where $a_1 = a_3$ or $a_2 = a_3$ follow similarly. 
In contrast, if $a_1 a_2 (a_1 - a_2) \neq 0$
then every mass linear 
$H$ is essential 
unless $\gamma_1 = \gamma_2 = \gamma_3$.
\QED

%%%%%%%%%%%%%%%%%%%%%%%%%%%%%%%%%%%%%%%%%%%%%%%%%%%%%%%%%%%%%%%%%%%%%%%
\subsection{Smooth $3$-dimensional polytopes}\labell{ss:3dim}
%%%%%%%%%%%%%%%%%%%%%%%%%%%%%%%%%%%%%%%%%%%%%%%%%%%%%%%%%%%%%

In this subsection 
we describe mass linear functions on $3$-dimensional smooth polytopes.
In particular,
we show that the only $3$-dimensional smooth polytopes which 
admit essential mass linear functions are $\Delta_2$ bundles over $\Delta_1$.

We will need the following lemma.
Note that the claim  
is false without the smoothness assumption;
this lemma 
is the one extra piece of information that we use in the smooth case.

\begin{lemma}\labell{prodsimp}
Let $\Delta$ be a smooth polytope which is combinatorially
equivalent to $\Delta_k \times \Delta_n$.
Then $\Delta$ is either a $\Delta_k$ bundle over $\Delta_n$,
or  a $\Delta_n$ bundle over $\Delta_k$.
\end{lemma}

\begin{proof}
By assumption, we can label the facets $F_1,\dots,F_{k+1}$ and $F'_1,\dots,F'_{n+1}$ so
that $F_I \cap F'_J \neq \emptyset$ for any proper subsets $I \subset \{1,\dots,k+1\}$
and $J \subset \{1,\dots,n+1\}$.
Denote the corresponding conormals by $\eta_1,\dots,\eta_{k+1}$ and
$\eta'_1,\dots,\eta'_{n+1}$.
Given $i \in \{1, \dots, k + 1\}$ and $j \in \{1,\dots, n+1\}$,  
let $I_i = \{1, \dots, k+1\} \smallsetminus \{i\}$ and
$J_j = \{1,\dots,n+1 \} \smallsetminus \{j\}$.
Since $F_{I_{k+1}} \cap  F'_{J_{n+1}} \neq \emptyset$,
the vectors $\eta_1,\dots,\eta_k$ and $\eta'_1,\dots,\eta'_n$ form a basis for $\ft$
(and thus identify $\ft$ with $\R^{k+n}$).
Write
$$
\eta_{k+1} = \sum_{i=1}^k a_i\eta_i + \sum_{j=1}^n a_j'\eta_j'
\quad   \mbox{and} \quad
\eta_{n+1}' = \sum_{i=1}^k b_i\eta_i + \sum_{j=1}^n 
b_j'\eta_j'.
$$

Given $i \in \{1, \dots, k + 1\}$ and $j \in \{1,\dots, n+1\}$,  
let  $A_{i,j}$ denote the matrix with
columns $\eta_{1},\dots,\eta_{i-1}, \eta_{i+1},\dots,\eta_{k+1}, \eta'_1,\dots,
\eta'_{j-1}, \eta'_{j+1},\dots,\eta'_{n+1}$. 
Since $F_{I_i} \cap F'_{J_{j}} \neq \emptyset$ and $\Delta$ is smooth,
$\det A_{i,j} =  \pm 1$.
In fact, $\det A_{i,j} = (-1)^{(k+1- i) + (n+1 - j)} $.
In particular, for all $i \in \{1,\dots,k\}$,
$\det A_{i,n+1} = (-1)^{k+1-i}$;
therefore, $a_i = -1$ for all such $i$.
A similar argument shows that $b_j' = -1$ for all $j$.

It remains to show that either $a'_j=0$ for all $j \in \{1,\dots,n\}$ or 
$b_i=0$ for all $i \in \{1,\dots,k\}$.  Suppose not.  Then
there exists $i\in \{1,\dots,k\}$
and $j\in \{1,\dots,n\}$  such that  $b_i a_j'\ne 0$.
On the other hand, since the matrix 
$A_{i,j}$ has determinant 
$(-1)^{(k-i) + (n-j)}$,
$$
\det \left(\begin{array}{cc} -1&a_j' \\
b_i& -1\end{array}\right) = 1.
$$
But then $b_i a_j'= 0$, a contradiction.
\end{proof}

\begin{lemma}\labell{le:3da}
Let $H \in \ft$ be a mass linear function on a smooth
$3$-dimensional polytope $\Delta \subset \ft^*$
with more than two asymmetric facets.
If every asymmetric facet is pervasive, 
then one of the following statements holds:
\begin{itemize}
\item $\Delta$ is the simplex $\Delta_3$.
\item $\Delta$ is a $\Delta_1$ bundle over $\Delta_2$; the base facets
are the  asymmetric facets.
\item $\Delta$ is a $\Delta_2$ bundle over $\Delta_1$; the fiber facets
are the asymmetric facets.
\end{itemize}
\end{lemma}

\begin{proof}
Let $A$ be the set of asymmetric facets 
and let $S$ be the set of symmetric facets.
By Lemma~\ref{le:symedge},
$\Delta$ has no symmetric edges. 
Since 
the intersection of two symmetric facets is by definition
a symmetric face, 
this implies that
the symmetric facets are disjoint. 
Since  
every asymmetric facet is pervasive,
there are $\frac{1}{2}|A|(|A| - 1) + |A| |S|$ edges.
Since the polytope is simple and $3$-dimensional, 
there are two-thirds as many vertices.
Since the Euler number of $\Delta$  is  $2$,
$$
|A| + |S| - \frac{1}{6}|A|(|A| -1) - \frac{1}{3}|A| |S| = 2,
\quad \mbox{that is,}
$$
$$ \left( |A| - 3 \right) (\left (|A| + 2|S| - 4 \right) = 0.
$$
Assume first that $|A| + 2|S| = 4$.
Since every $3$-dimensional polytope has at least four facets,
this implies that $|A| = 4$ and $|S| = 0$;
hence $\Delta$ is the simplex $\Delta_3$.

So assume that  $|A| = 3$;  
label the asymmetric facets $F_1, F_2$, and $F_3$.
Since no $3$-dimensional polytope has only
three facets, $|S| \neq 0$.
Since every symmetric facet $G$ 
intersects  each asymmetric facet,
$G$ 
is a triangle, and hence also intersects $F_1 \cap F_2$, 
$F_1 \cap F_3$, and $F_2 \cap F_3$. Since the edge $F_1 \cap F_2$
has two vertices, $|S|$ is $1$ or $2$. 
In the first case, 
$\Delta$ is again a tetrahedron.
In the second case,
$\De$ is combinatorially equivalent to the product $\De_2\times \De_1$. 
The claim now follows from  Lemma~\ref{prodsimp}. 
\end{proof}

\begin{rmk}\rm\labell{rmk:3dimperv}
More generally, let  $\Delta$ be a $3$-dimensional
simple polytope; assume 
that every facet is pervasive. 
As in the argument above,
the fact that there are two-thirds as many vertices
as edges and that the Euler number of $\Delta$ is $2$ implies that $\Delta$
is combinatorially equivalent to a $3$-simplex.

In contrast, there are many $4$-dimensional simple polytopes 
with the property 
that every facet is pervasive.
For example,  $C_4(N)$,  the $4$-dimensional cyclic polytope
with $N$ vertices is simplicial (each facet contains exactly $4$ vertices) and
$2$-neighborly (each  pair of vertices is joined by an edge); cf. \cite{Gr}.
Therefore, its dual polytope is a 
$4$-dimensional simple polytope with $N$
facets, each of which is pervasive.
\end{rmk}

We are now ready to classify $3$-dimensional smooth polytopes
with essential mass linear functions.

\MS

\NI
{\bf Proof of Theorem~\ref{thm:3dim}.}\,\,
By Lemma~\ref{le:twoasym}, if $H \neq 0$,
then $\Delta$ has at least two asymmetric facets.
If it has exactly two, then the  claim
follows by Corollary~\ref{cor:2asym}.  
 So assume that
$\Delta$ has more than two asymmetric facets.
By  Proposition~\ref{prop:flat},
after possibly subtracting an inessential function,
we may assume that every asymmetric facet is  pervasive.
Hence, by Lemma~\ref{le:3da}, $\Delta$ is either
the simplex $\Delta_3$, a $\Delta_1$ bundle over $\Delta_2$ with
symmetric fiber facets, or a $\Delta_2$ bundle over $\Delta_1$ with
symmetric base facets.
In either of the first two 
cases,
the asymmetric facets are equivalent.
Hence  $H$ is inessential by Corollary~\ref{cor:ea}. 
Therefore,  $\Delta$ is
a $\Delta_2$ bundle over $\Delta_1$,
that is, $\Delta$ is the polytope described in Example~\ref{ex:1}.
\QED
\MS

With a little more effort, we can give a complete
list of all smooth $3$-dimensional polytopes which
admit nonconstant mass linear functions.

\begin{thm}\labell{thm:3d1}
Let $H \in \ft$ be a 
nonzero 
mass linear function on a smooth
$3$-dimensional polytope $\Delta \subset \ft^*$.
One 
of the following statements holds:
\begin{itemize}
\item $\Delta$ is a bundle over a simplex.
\item $\Delta$ is a $1$-fold expansion.
\end{itemize}
\end{thm}

\begin{proof}
The $3$-simplex $\Delta_3$ is the expansion of the $2$-simplex $\Delta_2$ along an edge.
Moreover,
a $\Delta_2$ bundle over $\Delta_1$, a $\Delta_1$ bundle over $\Delta_1 \times \Delta_1$,
and the product $\Delta_1 \times \Delta_1 \times \Delta_1$ are all bundles over the $1$-simplex
$\Delta_1$.
Therefore, this result follows immediately from 
Proposition~\ref{prop:3d} below.
\end{proof}

\begin{prop}\labell{prop:3d}
Let $H \in \ft$ be a 
nonzero 
mass linear function on a smooth
$3$-dimensional polytope $\Delta \subset \ft^*$.
One 
(or more)
of the following statements holds:
\begin{itemize}
\item $\Delta$ is a bundle over $\Delta_1$; the base facets are
the  asymmetric facets.
\item $\Delta$ is a $1$-fold expansion; the base-type facets are 
the asymmetric facets.
\item $\Delta$ is the simplex $\Delta_3$.
\item $\Delta$ is a $\Delta_1$ bundle over $\Delta_2$; 
the base facets are the asymmetric facets.
\item $\Delta$ is a $\Delta_2$ bundle over $\Delta_1$; if either base facet
is asymmetric then both are.
\item $\Delta$ is a $\Delta_1$ bundle over $\Delta_1 \times \Delta_1$;
the base facets are the asymmetric facets.
\item $\Delta$ is the product $\Delta_1 \times \Delta_1 \times \Delta_1$;
every facet is asymmetric.
\end{itemize}
\end{prop}

\begin{proof}
By Lemma~\ref{le:twoasym}, $\Delta$ has at least two asymmetric facets.
If it has exactly two, then the  claim
follows by Corollary~\ref{cor:2asym}. 
So assume that
$\Delta$ has more than two asymmetric facets.

If all the asymmetric facets are pervasive, then the
claim holds by Lemma~\ref{le:3da}.
So assume that there is an asymmetric  facet $F$ which is not pervasive.

Corollary~\ref{cor:newflat}
implies that $\Delta$ is bundle over $\Delta_1$ with base facets
$F$ and $G$.  By Lemma~\ref{le:symasym}, $G$ is also asymmetric.
By Proposition~\ref{prop:bund},  there exists a mass linear function
$\Tilde{H} \in \ft$ so that
$\Tilde{H} - H$ is inessential,
the base facets  are $\Tilde{H}$-symmetric,
and at least one  fiber  facet is still 
$\Tilde{H}$-asymmetric. By  
Proposition~\ref{prop:symface},
the restriction of $\Tilde{H}$ to the $2$-dimensional
polygon $F$ is mass linear 
and at least one edge  of $F$ is $\Tilde{H}$-asymmetric.
By 
Proposition~\ref{prop:2dim}, this implies that $F$
has at most four edges.

If $F$ is a triangle, then  $\Delta$ is a $\Delta_2$ bundle over $\Delta_1$
and both base facets are $H$-asymmetric 
(although they are $\Tilde H$-symmetric).

If $F$ is a quadrilateral, 
$\Delta$ is combinatorially equivalent
to the product $\Delta_1 \times \Delta_1 \times \Delta_1$.
Suppose that $\Delta$ has an $H$-symmetric facet.
Since every symmetric facet intersects every asymmetric facet
and there are at least three asymmetric facets,
this implies that $\Delta$ has exactly two symmetric facets 
and they do not intersect.
Since none of the facets are pervasive,
Proposition~\ref{prop:asym} implies that every asymmetric facet 
is flat, that is, there exists a hyperplane which contains the 
conormal of every facet that meets it.
If $F$ and $F'$ are asymmetric facets which intersect, then the associated
hyperplanes must be 
distinct, and so  
they intersect in a line $L$.
Since both symmetric facets
intersect $F$ and $F'$, the two symmetric faces 
have conormals in $L$ and hence are parallel.  Moreover 
each is affine equivalent to the product $\Delta_1 \times \Delta_1$.
Hence, $\Delta$ is a $\Delta_1$ bundle over $\Delta_1 \times \Delta_1$
and the fiber facets are symmetric.
On the other hand, if none of the facets are $H$-symmetric then
a very similar argument proves that every pair of opposite faces
is parallel.  
Hence, $\Delta$ is the product $\Delta_1 \times \Delta_1 \times \Delta_1$.
\end{proof}

We will need the following result in the next paper.

\begin{lemma}\labell{sum0}
Let $H \in \ft$ be a mass linear function on a smooth $3$-dimensional
polytope $\Delta \subset \ft^*$.
Then $$\sum_{i=1}^N \gamma_i = 0,$$
where $\gamma_i$ is the coefficient of the support number of the
facet $F_i$ in the linear function $\langle H, c_\Delta \rangle.$
\end{lemma}

\begin{proof}
If  $H$ is inessential, then $\sum_i \gamma_i = 0$ by 
Proposition~\ref{prop:inessential}.
Otherwise, by Theorem~\ref{thm:3dim}, $\Delta$ is a $\Delta_2$ bundle over $\Delta_1$.
The result now follows immediately from Proposition~\ref{prop:Mab}.
(See also Remark~\ref{rmk:easier}.)	
\end{proof}

%%%%%%%%%%%%%%%%%%%%%%%%%%%%%%%%%%%%%%%%%%%%%%%%%%%%%%%%%%%%%%%%%%%%%%%% 
\section{Relationship to Geometry} \labell{s:geometry}
%%%%%%%%%%%%%%%%%%%%%%%%%%%%%%%%%%%%%%%%%%%%%%%%%%%%%%%%%%%%%%%%%%%%%%%% 

We begin by constructing the 
symplectic
toric manifold 
associated to
each smooth polytope.
Let $ \Delta  = 
\bigcap_{i = 1}^N \{ x \in \ft^* \mid \langle \eta_i, x\rangle \leq \kappa_i  \} $
be a smooth polytope
with facets $F_1,\dots,F_N$.
Let  $e_1,\ldots,e_N$  
denote the  standard basis for $\R^N$,
and  $e_1^*, \ldots, e_N^*$ 
denote the dual basis for 
$(\R^N)^*$.
Let $\pi \colon  \R^N \to \ft$ denote the linear function  given by
$\pi(e_i) = -\eta_i$
for all $i$; let $\fk$ denote its kernel.
Since $\Delta$ is bounded, 
$\pi$ is surjective.
Since  $\eta_i \in \ell$ for all $i$, $\pi$ induces a surjective map
from $(S^1)^N$ to 
$T: = \ft/\ell$;
since $\eta_1,\ldots,\eta_N$ generate $\ell$,
the kernel $K$ of this map is connected.
Thus we have dual short exact sequences
of vector spaces
and a short exact sequence of groups
\begin{align}
\labell{eq:pik}
& 0 \to \fk \stackrel{\iota}{\to}  \R^N \stackrel{\pi}{\to} \ft \to 0,
\\
\labell{eq:dual}
& 0 \to \ft^* \stackrel{\pi^*}{\to} (\R^N)^* \stackrel{\iota^*}{\to} \fk^* \to 0,
\qquad \mbox{and} 
\\
\labell{eq:piK}
&{1} \to K \stackrel{\iota}{\to}  (S^1)^N \to T \to {1}.
\end{align}

The  moment map 
$\mu \colon  \C^N\to \fk^*$ 
for  the 
natural $K$ action on $\C^N$
is given by
$\mu(z) = \io^*\left( \frac12 |z_1|^2,\dots,\frac12 |z_N|^2\right)$. 
Thus the image $\mu(\C^N)$ 
is the projection of the positive quadrant in $(\R^N)^*$ 
 and so has outward conormals
$\pi(-e_1) = \eta_1, \dots, \pi(-e_N) = \eta_N$.
Let 
$(M_\Delta,\omega_\Delta)$ be 
the symplectic quotient 
$\C^N/\!\!/K: = \mu^{-1}(\io^*(\ka))/K$. 
The $(S^1)^N$ action on
$\C^N$ descends to an effective
$T$ action on $M_\Delta$; the moment map $\Phi_\Delta \colon
M_\Delta \to \ft^*$ for this action is given  by
$\Phi_\Delta([z]) = (\pi^*)^{-1}\left( \frac12 |z_1|^2 - \kappa_1,
\ldots, \frac12 |z_N|^2 - \kappa_N \right)$ for all $[z] \in M_\Delta$.
Hence, $\Phi_\Delta(M_\Delta) = \Delta$, as required.

Note also that $M_\De$ has a canonical $\om_\De$-compatible complex structure $J_\De$, and hence has a natural K\"ahler structure.  
To define $J_\De$, denote by $(\om_0,J_0)$ the standard symplectic form and complex structure on $\C^N$ and define $V: = \mu^{-1}(\io^*(\ka))$. Then the bundle $TV\cap J_0(TV)$ is $\om_0$-orthogonal to the $K$-orbits and its complex structure descends to $J_\De$ on $T(M_\De)$.
 Moreover, $J_\De$ is unique up to a $T$-equivariant K\"ahler isometry.\footnote{Throughout 
  we have assumed that our presentation of $\De$ is minimal, i.e. that each facet $F_i$ is nonempty.  If one inserted \lq\lq ghost" facets, then the above 
  procedure might give a different K\"ahler structure.
For further details, see  for example \cite{CDG}.}

\MS

\begin{rmk}\labell{alternate}\rm
There is an alternate way to construct 
the manifold $M_\Delta$ described above.
Given $z \in \C^N$, let  $I_z$ be the
set of 
$i \in \{1,\ldots, N\}$ 
so that $z_i = 0$. 
Define a set
$$
\Uu := 
\Big\{ z \in \C^N \Big\vert \bigcap_{i \in I_z}\; F_i \neq \emptyset  \Big\}.
$$
As a complex manifold,  $M_\Delta = \Uu/K_\C$,
where  $K_\C$ denotes the complexification of $K$.
Moreover, under this quotient map,
the $i$'th coordinate hyperplane in $\C^N$ corresponds
to $\Phi^{-1}(F_i)$, the preimage of the $i$'th facet
under the moment map, 
For more details, see Audin \cite[Ch~VII]{Aud} or Cox--Katz~\cite{CK}.
From this description, it is clear that
$J_\De$ does not depend on $\ka$.
\end{rmk}

The cohomology $H^2(M;\R)$  
is naturally isomorphic to $\fk^*$.
Under this isomorphism, the
cohomology class of the symplectic form
$\omega$ corresponds to  $\iota^*(\kappa) \in \fk^*$.  
The integral $\frac 1{n!}\int_{M} \omega^n$ is (up to a normalizing constant) the volume of $\Delta$
with respect to standard Lebesgue measure.

Many of the terms that we introduce for polytopes have 
geometric interpretations. For example,
a facet $F$ of a smooth polytope $\Delta$
is flat exactly if the corresponding 
divisor $\Phi^{-1}(F)$ in $M$ has self-intersection zero.

\begin{rmk}\labell{rmk:geobun} \rm
Polytopes which are bundles in the sense of Definition~\ref{def:bund}
correspond to bundles of symplectic toric manifolds.
To see this, let 
$\Delta \subset \ft^*$ be a bundle with fiber
$\Tilde{\Delta} \subset \Tilde{\ft}^*$ over the base $\Hat\Delta
\subset \Hat{\ft}^*$;
we will use the notation of Definition~\ref{def:bund}.
Construct  the associated toric manifolds  $M_\Delta = \Uu/K_\C$,  
$M_{\Tilde\Delta} = \Tilde{\Uu}/\Tilde{K}_\C$, and
$M_{\Hat\Delta} = \Hat{\Uu}/\Hat{K}_\C$
as in Remark~\ref{alternate};
identify $\C^{\Tilde N + \Hat N}$ with
$\C^{\Tilde N} \times \C^{\Hat N}$.
Since $\Delta$ is combinatorially equivalent to 
$\Tilde \Delta \times \Hat \Delta$, $\Uu = \Tilde \Uu \times \Hat \Uu$.
Moreover, since 
$\Tilde{\eta}_j\,\!' = \iota(\Tilde{\eta}_j)$ for all $j$ and 
$\pi({\Hat{\eta}_i}\,\!') = \Hat{\eta}_i$ for all $i$,
$\Tilde K$ is the
intersection of $K$ with 
$(S^1)^{\Tilde N}  \subset
(S^1)^{\Tilde N} \times (S^1)^{\Hat N}$,
and $\Hat K$ is the image of $K$ under the
natural projection $(S^1)^{\Tilde N} \times (S^1)^{\Hat N} \to
(S^1)^{\Hat N}$;
in particular,  $\Hat K = K/ \Tilde K$.
Therefore, the natural projection
from $K \subset (S^1)^{\Tilde N} \times (S^1)^{\Hat N}$ to
$(S^1)^{\Tilde N}$
induces a homomorphism 
$\rho$ from
$\Hat K \simeq K/ \Tilde K$ to $\Tilde T = (S^1)^{\Tilde N}/ \Tilde{K}.$
The toric manifold
$M_\Delta$ is the associated bundle
$$ 
M_\Delta = M_{\Tilde{\Delta}} \times_{\Hat{K}_\C} \Hat \Uu ,
$$
where $\rho$ describes the action of $\Hat{K}_\C$  on  $M_{\Tilde{\Delta}}$.
\end{rmk}

\begin{example}\labell{ex:ZK} \rm 
Given 
$a =(a_1,a_2) \in \R^2$, 
let $Y_a$ 
be the polytope
defined in Example~\ref{ex:1}. 
As we saw in Example~\ref{ex:2},
$Y_a$ is a $\Delta_2$ bundle over $\Delta_1$.
Recall that the outward conormals
$\eta_1,\ldots,\eta_5$ 
are given by
$$
(-1,0,0),\;\;(0,-1,0),\;\;(1,1,0),\;\;(0,0,-1),\;\;\mbox{and }\;\; 
(a_1,a_2,1);
$$
these satisfy
$$ 
\eta_1+\eta_2+\eta_3=0
\quad \mbox{and} \quad
a_1 \eta_1 +a_2 \eta_2 + \eta_4 + \eta_5=0.
$$
Thus,
by  definition of the map $\pi$ in Equation (\ref{eq:pik}),
 $K$ is the subtorus of $(S^1)^5$ generated by the 
elements
$(1,1,1,0,0)$ and 
$(a_1,a_2,0,1,1)$ in $\Z^5$.
Moreover, by Definition~\ref{alternate}, 
$$
\Uu = \left( \C^{3} \smallsetminus \{0\} \right) \times \left(\C^2 \smallsetminus \{0\} \right) .
$$
Therefore $M_{Y_a}$ 
is the $\C P^2$ bundle 
$$ 
\CP^2 
\times_{\C^*} \left( \C^2 \smallsetminus \{0\} \right)
$$
over $\C P^1$, where 
$\C^*: = \C\less \{0\}$
acts on $\C P^2$  by 
$$
\lambda \cdot [z_1:z_2:z_3] = [\lambda^{a_1} z_1: \lambda^{a_2} z_2: z_3], 
$$
that is,
$M_{Y_a} = 
\Proj\bigl(\Oo(a_1) \oplus \Oo(a_2) \oplus \Oo(0)\bigr)$.
\end{example}

\begin{rmk}\labell{rmk:geoexpand} \rm
Polytopes which are $1$-fold expansions in the sense of 
Definition~\ref{def:expand}
correspond
to symplectic toric manifolds which are 
smooth symplectic pencils;
that is, they are swept out by a family 
$\Tilde M_\la, \la \in \C P^1$, of (real) codimension $2$ symplectic  submanifolds (often called fibers)   that intersect 
transversely 
along a codimension $4$ symplectic  submanifold called the axis or base locus. 
To see this,
let $\Delta \subset \ft^*$ be the $1$-fold expansion of a polytope 
$\Tilde \Delta
\subset \Tilde \ft^*$. 
Construct the associated toric manifolds
$M_\Delta = \Uu/K_\C$ and $M_{\Tilde{\Delta}} = \Tilde{\Uu}/\Tilde{K}_\C$,
as in Remark~\ref{alternate}.
We will use the notation of Definition~\ref{def:expand}
and identify $\C^N$ with $\C^2 \times \C^{N-2}$.
There is a well-defined map from  $\{ z \in \Uu  \mid z_1 z_2 \neq 0 \}$
to $\C^2 \smallsetminus \{0\}$ given by
$z \mapsto (z_1,z_2)$.
Since the $\Tilde{\eta_j}$ all lie in the hyperplane $\Tilde \ft$,
this  descends to an equivariant 
holomorphic map from 
$M_\Delta \smallsetminus 
\Phi^{-1}(\Hat{F}_1 \cap \Hat{F}_2)$ to $\CP^1$.
Therefore, $M_\Delta$ is a pencil with axis 
$\Phi^{-1}(\Hat{F}_1 \cap \Hat{F}_2)$.  
Notice that the toric 
manifold obtained by blowing up  $M_\De$ along this axis 
as in Remark~\ref{rmk:blowexp} is a fibration over $\C P^1$ with fiber 
$M_{\Tilde \De}$.
\end{rmk}

Our proof of Proposition~\ref{prop:symp} is based on
Weinstein's action homomorphism
$$
\Aa_{\om} \colon  \pi_1(\Ham(M,\om))\to \R/\Po.
$$
Here,  
$\Ham(M,\om)$ denotes the group of Hamiltonian 
symplectomorphisms of $(M,\om)$
and $\Po\subset \R$, the period group of $\om$, is the image of $\pi_2(M)$ 
under the homomorphism $\al\mapsto \int_{\al}\om$.
The homomorphism  $\Aa_\om$ is defined as follows:    If the loop
$\{\phi_t\}_{t\in S^1}$ 
in $\Ham(M,\om)$ is generated by the
mean normalized Hamiltonian\footnote{
We define the sign of $K_t$ by requiring that $\dot\phi_t = 
X_{K_t}$
where $\om(X_{K_t},\cdot) = -dK_t$.
We say that $K$ is {\bf mean normalized} iff
$\int_M K \om^n = 0$.} 
$K_t$ for $ t\in [0,1]$, then
$$
\Aa_{\om}(\{\phi_t\}) : = \int_0^1 K_t(\phi_t(p)) dt - \int_D u^*(\om)
$$
where $p$ is any point in $M$ and $u \colon D^2\to M$ is any smooth map 
such that $u(e^{2\pi it}) = \phi_t(p)$.
If the loop $\{\phi_t\}_{t\in S^1}$ is  a circle subgroup
$\Lambda \subset \Ham(M,\om)$
generated by a mean 
normalized Hamiltonian $K \colon M\to \R$ then 
we may take $p$ to be a fixed point  
of the action
and $D$ to be the constant disc.  Hence 
$\Aa_{\om}(\La)$ is the image in $\R/\Po$ of any critical value of $K$.
This is well defined because 
for any fixed points $p$ and $p'$
the difference $K(p)-K(p')$  
is the integral of $\om$ over the 
$2$-sphere formed by rotating an arc from $p$ to $p'$ 
by the $S^1$-action.
\MS

\NI
{\bf Proof of Proposition \ref{prop:symp}.}\,\,
First, note that since toric manifolds are simply connected, 
$\Symp_0(M,\omega)$ is $\Ham(M,\omega)$.
Given $H \in \ell$, 
the mean
normalized Hamiltonian
for the action of $\Lambda_H$ on $M$ is 
$$
K =  \langle H,  \Phi  \rangle  - \langle H, c_\Delta \rangle,
$$ 
where $c_\Delta$ denotes the center of mass of $\Delta.$
Since vertices of $\Delta$ correspond to 
fixed points in $M$, 
this implies that
$$
\Aa_{\om}(\La_H) =  \langle H ,   v  \rangle - \langle H, c_{\Delta} \rangle 
\in 
\R/\Po 
$$
for any vertex  $v$  of $\Delta$.

Clearly if $\Lambda_H$ vanishes in $\pi_1(\Symp_0(M,\omega))$,
then $\Aa_{\om}(\Lambda_H) = 0$.
In fact, it
is immediate from Moser's homotopy argument that  
if $\La_H$ contracts in $\Symp_0(M,\om)$,  
then it also contracts for all sufficiently small perturbations of 
$\om$.    
Since changing $\kappa$ corresponds to changing the
symplectic form  on $M$,  
for any vertex $v$ of $\De$
the image of $\La_H$ under the action homomorphism
$$
\Aa_{\om(\ka)}(\La_H) 
= \langle H, v \rangle  - \langle H,c_{\Delta}(\ka) \rangle
$$
must lie in ${\mathcal P}_{\om(\ka)}$ for all $\ka$ in some open set 
of $\R^N$.
Since ${\mathcal P}_{\om(\ka)}$ is generated by the lengths of the  edges
of the polytope $\Delta$ (considered as a multiple of a primitive vector),
it is a finitely generated subgroup of $\R$ 
whose generators are linear functions of the $\kappa_i$ with
integer coefficients.
Similarly, $ \langle H, v \rangle$ is a linear function
of the $\kappa_i$ with integer coefficients.
Since the function $\ka\mapsto \Aa_{\om(\ka)}(\La_H)$ is 
continuous, it follows that the 
function  $\ka \mapsto \langle H,c_{\Delta}(\ka) \rangle$
is also a linear function of the $\kappa_i$ with integer coefficients
 as $\ka$ varies in some open set.  
By Lemma~\ref{le:local}, this implies that
it is a linear function of $\kappa_i$ with integer 
coefficients 
for all $\ka \in \Cc_K$.
\QED\MS

The next proposition   
gives a geometric interpretation of the equivalence 
relation on the facets of 
$\De$ in terms of  $\Isom_0(M)$, 
 the identity component of the group of isometries of $(M,g_J)$.
Recall that these act by biholomorphisms of $(M,J)$ and hence also preserve $\om$.  
To see this, note that 
 the K\"ahler condition implies that $J$ is invariant under parallel translation with respect to the Levi--Civita connection of $g_J$, while, for each path $\phi_t, t\in [0,1],$ in $\Isom_0(M)$ with $\phi_0=id$ and each $x\in M$, the derivative $d\phi_t(x): T_xM\to T_{\phi(x)}M$ is the linear map given by parallel translation along the path $\phi_s(x), s\in [0,t]$.

\begin{prop}\labell{prop:preaut}
Let $(M,\omega,T,\Phi)$ be an $n$-dimensional 
 symplectic toric manifold
with moment polytope
$\Delta  \subset \ft^*$,  and 
let $\Isom_0(M)$ be the identity component of the
associated K\"ahler isometry group.
Let $\Ii$ denote the set of
equivalence classes of facets of $\Delta$,
and let $K$  be defined as in \eqref{eq:piK}.
Then 
$$
\Isom_0(M)  
\cong \Bigl( \prod_{I \in \Ii} \U(|I|) \Bigr)/K,
$$
where for each equivalence class $I \in \Ii$,  
the unitary group $\U(|I|) \subset \U(N)$ acts on the subspace 
of $\C^N$ spanned by 
$\{e_i\}_{i \in I}$.
Under this identification, 
$T = (S^1)^N/K$
and $\prod_{I \in \Ii} \U(|I|)$ is the centralizer of $K$ in $U(N)$.
Moreover, $\Isom_0(M)$ is a maximal connected  compact subgroup of 
$\Symp_0(M,\omega)$.
\end{prop}

\begin{proof}
Let $(M_\Delta, \omega_\Delta, J_\Delta,\Phi_\Delta)
= \C^N /\!\!/K$ be the K\"ahler toric manifold
associated to the moment polytope  $\Delta \subset \ft^*$,
as above.  Since $M$ and $M_\Delta$ are equivariantly
symplectomorphic by \cite{Del}, we may identify  $M$ and $M_\Delta$.

Let $Z(K)$ denote the 
centralizer
of $K \subset (S^1)^N$ in $U(N)$. 
Our first claim is that  $Z(K)$ 
is the product $\prod_{I\in \Ii} \U(|I|)$.
To see this,  note first that since $\eta_k = \pi(e_k)$ for all $k$,
$\xi \in \ft^*$ satisfies
$$
\langle  \eta_i,\xi \rangle = 1 = - \langle  \eta_j,\xi \rangle
\quad \mbox{and otherwise} \quad 
\langle  \eta_k,\xi \rangle = 0 \ \forall \, k
$$
exactly if 
$$
\langle  e_i,\pi^*(\xi) \rangle = 1 = - \langle  e_j,\pi^*(\xi) \rangle
\quad \mbox{and otherwise} \quad 
\langle  e_k, \pi^*(\xi) \rangle = 0  \ \forall \, k,
$$
that is, 
exactly if $\pi(\xi) = e^*_i - e^*_j$. 
Therefore,  Lemma~\ref{le:equiv} implies that 
two facets $F_i$ and $F_j$ of $\De$ are equivalent exactly if
$e^*_i - e^*_j \in \pi^*(\ft^*)$.
But this holds exactly if $\iota^*( e^*_i - e^*_j) = 0$,
that is, iff $\iota^*(e^*_i) = \iota^*(e^*_j)$.
Since the weights for the $K$ action on $\C^N$
are exactly
$\iota^*(e^*_1),\ldots, \iota^*(e^*_N)$, the claim follows immediately.

Now observe that 
$Z(K)\subset U(N)$
acts on $\C^N$, 
preserving the standard symplectic and complex structures on $\C^N$. 
This implies that 
the quotient group $Z(K)/K$ acts on $M_\Delta$, preserving $\omega_\Delta$ and $J_\Delta$
and hence also the K\"ahler metric $g_J$.
Thus, 
because $Z(K)$ is connected, $Z(K)/K \subset \Isom_0(M_\Delta)$.
 But
$\Isom_0(M_\Delta)$ 
is a compact connected Lie group.
Therefore, to complete the proof it
is enough to prove that $Z(K)/K$ is a maximal compact connected 
subgroup of $\Symp_0(M_\Delta,\omega_\Delta)$.

To this end, 
let $G$ be any connected compact  
subgroup of $\Symp_0(M_\Delta,\omega_\Delta)$ that contains
$Z(K)/K$. 
Then $G$ is a Lie group.
By Lemma~\ref{le:weyl} below, $T$ is the maximal torus of $G$
and the Weyl group of $G$ acts on $\ft^*$ as a subgroup of
the group of robust symmetries of $\Delta$.
By the same argument, $T$ is the maximal torus of $Z(K)/K$
and the Weyl group of $Z(K)/K$ acts as a subgroup
of the group of robust symmetries of $\Delta$.
Moreover, given any two equivalent facets $F_i$ and $F_j$,
it is easy to check that there is an element of the Weyl group
of $Z(K)/K$ which exchanges $F_i$ and $F_j$ and
otherwise takes each facet to itself.
Therefore,
the Weyl group of $Z(K)/K$ acts on $\ft^*$ as the group of robust symmetries of $\Delta$.
Since, $Z(K)/K \subset G$, this implies that
the two groups have the same Weyl group.  
But, because  $Z(K)/K = \prod_{I\in \Ii} \U(|I|)/K$, 
 each of its simple factors 
has a root system of type $A_{|I| - 1}$. 
Therefore, there is no way to extend its root system  without either increasing the rank or enlarging the Weyl group.
It follows that $G$ has the same roots as $Z(K)/K$, so that the two groups must coincide.
\end{proof}

\begin{lemma}\labell{le:weyl}
Let $(M,\omega,T,\Phi)$ be a symplectic toric manifold.
Let $G \subset \Symp(M,\omega)$ be a compact connected Lie subgroup
that contains $T$.
Then $T$ is a maximal torus of $G$, and the
Weyl group $W$ of $G$ acts on $\ft^*$ as  a subgroup 
of the group of robust
symmetries of the moment polytope $\Delta \subset \ft^*$.
\end{lemma}

\begin{proof}
Since $M$ is simply connected, 
there is a moment map $\Phi_G \colon M \to \fg^*$
for the $G$ action.
Since $\Phi_G$ is equivariant, 
the moment image $\Phi_G(M) \subset \fg^*$ 
is invariant under  
the coadjoint action of $G$ on $\fg^*$. 

Let $\iota \colon \ft \to \fg$ denote the inclusion map.
The dual map $\iota^* \colon \fg^* \to \ft^*$  is
equivariant with respect to the the coadjoint action
of $N(T)$, 
where $N(T)$ denotes the normalizer of $T$ in $G$.
Moreover, 
by
 adding a constant if necessary, we may assume that
$\Phi_T = \iota^* \circ \Phi_G$. 
Therefore, the moment polytope 
$\Delta = \Phi_T(M) \subset \ft^*$
is also invariant under the coadjoint action of 
$N(T)$; hence,
$\Phi_T(M)$ is invariant under the action
of the Weyl group $W = N(T)/T$.

Finally, write the moment polytope as
$\Delta = 
\bigcap_{i = 1}^N \{ x \in \ft^* \mid \langle \eta_i, x\rangle \leq \kappa_i\}$,
where $N$ is the number of its (nonempty) facets. 
By the construction explained in
the beginning  
of this section, for all $\kappa' \in \Cc_\Delta$, 
there exists a $T$-invariant symplectic form $\omega_{\kappa'}$ 
with corresponding  moment polytope $\Delta(\kappa')$.
Let $\alpha \in \Omega^2(M;\R)$ be  any $G$-invariant
$2$-form in the cohomology class 
$[\omega' - \omega]$.
For $\epsilon > 0$ sufficiently small,
the form $\omega + \epsilon \alpha$ is a
$G$-invariant symplectic form with moment polytope 
$\Delta((1-\eps)\kappa+ \epsilon\kappa')$.
By the previous paragraph, this implies that, up to translation,
this polytope is also invariant  
under the action of the Weyl group. 
The conclusion now follows from Corollary~\ref{cor:xi}.
\end{proof}

\NI
{\bf Proof of Lemma~\ref{le:aut}.}\,\,  
Let $Z(K) =  \prod_{I \in \Ii} U(|I|)$. Then
 the natural map from 
$\pi_1((S^1)^N)$ to $\pi_1 \left(Z(K) \right)$
is surjective.
Moreover,
$H'=\sum_i\be_i e_i\in \Z^N$ generates a trivial element of 
$\pi_1(Z(K)) = \pi_1 \left( \prod_{I \in \Ii} \U(|I|) \right)$
exactly if it lies in the Lie algebra of 
$\prod_{I \in \Ii} \SU(|I|)$,
that is, exactly if $\sum_{i \in I} \be_i  = 0$ for all $i \in I$. 

Since $K$ is connected, the natural
maps from $\pi_1( (S^1)^N) $ to $\pi_1(T)$ and
from
$\pi_1 \left(Z(K) \right)$ to 
$\pi_1 \left(Z(K)/K\right)$ 
are  surjective.
Therefore,  the natural map
from $\pi_1(T)$  to  $\pi_1 \left(Z(K)/K\right)$
is surjective.
Moreover, $H  \in \ell$ generates a trivial element
of $\pi_1 \left(Z(K)/K\right)$ 
exactly if there exists 
$H' = \sum_i\beta_i e_i \in \pi^{-1}(H) \subset \Z^N$  that 
generates a trivial element of $\pi_1 \left(Z(K)\right)$,
where $\pi: \R^N\to \ft$ is as defined in equation (\ref{eq:pik}).
But this happens
 exactly if  $H 
= \sum_i \beta_i \eta_i \in \ell$, where $\beta_i \in \Z$ for
all $i$ and $\sum_{i \in I} \beta_i = 0$ for all $I \in \Ii$.
\QED\MS

\begin{example}\rm
Given $a = (a_1,a_2) \in \R^2$, let $Y_a$ be the polytope
defined in Example~\ref{ex:1}
Construct $M_{Y_a} = \Uu/K$ as in Example~\ref{ex:ZK}; 
let $Z(K)$
be the centralizer of $K$ in $U(N)$.
Recall that $K$ is the subtorus of $(S^1)^5$ generated
by the elements $(1,1,1,0,0)$ and $(a_1,a_2,0,1,1)$.
Therefore, 
\begin{itemize}\item 
If $a_1 a_2 (a_1 - a_2) \neq 0$, then
$Z(K) = S^1\times S^1\times S^1\times \U(2)$,  and so 
$\pi_1(Z(K)/K)$ has rank $2$. 
\item
If $a_1 = a_2 \neq 0$ (or if $a_1 = 0 \neq a_2$ or
$a_1 \neq 0 = a_2$) then
$Z(K) = \U(2) \times S^1\times \U(2)$ and  
$\pi_1(Z(K)/K)$ has rank $1$. 
\item
If $a_1 = a_2 = 0$
then $Z(K) = U(3) \times U(2)$ and $\pi_1(Z(K)/K)$ has rank $0$.
\end{itemize}
Thus the structure of  $Z(K)/K$ 
depends on the coefficients $(a_1,a_2)$.  
This is consistent with the calculation of the equivalence relation on the facets via Lemma~\ref{le:equiv} given in the proof of 
Proposition~\ref{prop:Mab}. Also,
by Proposition~\ref{prop:preaut}, $\Isom_0(M) = Z(K)/K$.  
It is hard to check this independently except in the case 
$a_1 = a_2 = 0$ when $\De$ is a product.  However, one can verify  the consistency of
Lemma~\ref{le:aut}, which implies that the set of inessential functions should have a corresponding dependence on $(a_1,a_2)$.  
This is the case, as the reader may check by 
looking at Remark~\ref{rmk:Mab}.
\end{example}

\begin{rmk}\labell{rmk:Cox}\rm (i)  A rational polytope $\De$ defines an
algebraic variety $X: = X_\De$.  In this case, it is customary to consider
a refinement $\sim'$ of the equivalence relation $\sim$ that takes multiplicities into account.  One way of doing this is to consider the
Chow group $A_{n-1}(X)$ of divisors in $X$, and to set 
$F_i\sim'F_j$ exactly if 
$F_i$ and $F_j$ 
have the same image in
$A_{n-1}(X)$; cf.
Cox~\cite[\S4]{Cox}. If $\De$ is smooth then $X=M$, $A_{n-1}(X)$ is 
naturally identified with $H_{2n-2}(M;\Z)$ and $F_i\sim' F_j$ exactly if the 
submanifolds $\Phi^{-1}(F_i)$ and $\Phi^{-1}(F_j)$ of $M$ are homologous.
It then follows from standard facts about toric manifolds that 
  the two equivalence relations are the same.  
But they differ in general, for example for the polygon with conormals 
$(-1,0), (0,-\frac 23), (1,\frac23)$ corresponding to the weighted projective space $\CP^2(1,2,3)$; cf. Example \ref{ex:rat4}.

\MS

\NI (ii)  Consider the group 
$\Aut(M)$ of biholomorphisms of the complex manifold $(M,J_\De)$.
Its structure
was worked out by Demazure \cite{De}.
As explained in Cox \cite{Cox}, for any rational polytope $\De$, its identity component
$\Aut_0(X)$ is isomorphic to a semidirect product 
  $R_u\rtimes G_s$ where the unipotent radical $R_u$ 
  is isomorphic to $\C^k$ for some $k\ge 0$ and 
  the reductive part $G_s$ is isomorphic to the product 
  $\left(\prod_{I\in \Ii'} GL(|I|,\C)\right)/K_\C$, where  
  $K_C$ is the complexification of $K$ and $\Ii'$ is the set of  equivalence classes with respect to the relation $\sim'$ discussed in (i).   
  Moreover, the elements of $\pi_0(\Aut(X))$ act nontrivially on homology.

Let us call a facet $F_i$ of $\De$ {\bf convex} 
 if all the conormals except for $\eta_i$ lie in
a closed half-space of $ \ft$.  
It follows easily from the description of Demazure roots in
\cite[\S4]{Cox} that
 in the smooth case the unipotent radical $R_u$ is trivial (i.e. 
 $\Aut(M)$ is reductive)  exactly if 
 each convex  facet of $\De$  is equivalent to at least one other facet.
 This fact about the structure of $\Aut(M)$ does not seem relevant in our situation, although it is in other geometric contexts; cf. 
 Remark~\ref{rmk:Nill}.
\end{rmk}

\begin{rmk}\rm  The first claim in Proposition~\ref{prop:preaut} is that  $\Isom_0(M) \cong Z(K)/K$.  One can also prove this
 by thinking of
$\Isom_0(M)$ as the identity component of $\Symp_0(M)\cap \Aut_0(M)$, where $\Aut_0(M)$ is as in the previous remark. 
Since $R_u$ contains no 
nontrivial compact subgroup,  $\Isom_0(M)$ consists of the
unitary elements in 
$\left(\prod_{I\in \Ii} GL(|I|,\C)\right)/K_\C$ and so is 
$Z(K)/K$.
\end{rmk}

%%%%%%%%%%%%%%%%%%%%%%%%%%%%%%%%%%%%%%%%%%%%%%%%%%%%%%%%%%%%%%%
\begin{appendix}
\section{Powerful facets\\
Appendix written with Vladlen Timorin}
%%%%%%%%%%%%%%%%%%%%%%%%%%%%%%%%%%%%%%%%%%%%%%%%%%%%%%%%%%%%%%%

The purpose of this appendix is to 
analyze mass linear function on $4$-dimensional simple polytopes
where every facet is asymmetric. 
In part II, we will need  these results in order to
classify mass linear functions on $4$-dimensional smooth polytopes.
We need one last important definition.

\begin{definition}
We say that a facet 
$F$ of a simple polytope $\Delta$ is {\bf powerful}
if it is connected to each vertex  of $\Delta$ by an edge.
\end{definition}

We will adopt the notations of \S \ref{ss:volume}.
In particular, given $H \in \ft$ and
a simple polytope $\Delta \subset \ft^*$, fix an
identification of $\ft^*$ with Euclidean space and define
$$ V = \int_\Delta dx,  
\quad \mu = \int_\Delta H(x) dx, \quad
\mbox{and} \quad \Hat H(\kappa) = \langle H, c_\Delta(\kappa) \rangle.$$
Let $\p_i$ denote the operator of differentiation with respect
to $\kappa_i$.

\begin{prop}\labell{prop:powerful}
Let $H \in \ft$ be a mass linear function on a simple polytope
$\Delta \subset \ft^*$.
Every 
asymmetric facet is powerful.
\end{prop}

\begin{proof}
Let $F_1$ be an 
asymmetric facet 
and let $v = F_J$ be a vertex
which is not connected to $F_1$ by any edge;
number the facets so that $J = \{2,\dots, n+1
\}$.
Since  $\Hat H$ is linear, its second derivatives all vanish.
By assumption, 
the intersection $F_1 \cap e_j$ is empty
for each $j \in J$, 
where $e_j$ is the edge $\bigcap_{i \in J \ssminus \{j\}} F_i$.
Hence, by Theorem~\ref{der}, if
we apply the differential operator $ \p_1\cdots\p_{n+1}$ to the formula 
 $\mu=\Hat H V$
 we obtain
$$
0 = K_v \, \p_1 \Hat H,
$$
where $K_v$ is a positive real number.
Since $F_1$ is an asymmetric facet $\p_1 \Hat H \neq 0$; this gives a contradiction.
\end{proof}

With this motivation, we analyze polytopes with powerful facets.

\begin{lemma}\labell{le:power}
Let $\Delta$ be a simple polytope.
If two powerful facets $F$ and $G$ do not intersect,
then  $\Delta$ is combinatorially equivalent to the product
$\Delta_1 \times F$.
If every facet is powerful, then
$\Delta$ is combinatorially equivalent to a product 
$\Delta_1 \times \cdots \times \Delta_1 \times \Delta'$,
where every facet of $\Delta'$ is both powerful and pervasive.
\end{lemma}

\begin{proof}
Assume that two powerful facets $F$ and $G$ do not intersect.
Then every vertex in $F$ has a unique edge which does not
lie in $F$, and that edge must meet $G$. Conversely
every edge which meets (but does not lie in) $G$ must meet $F$.
Since $\Delta$ is connected, this implies that $\Delta$ is combinatorially 
equivalent to $\Delta_1 \times F$.

The second claim follows by induction.
\end{proof}

\begin{rmk}\rm
We can now give
an alternate proof of Corollary~\ref{cor:newflat}.
Let $H \in \ft$ be a mass linear function on a simple
polytope $\Delta \subset \ft^*$.
Let $F$ be an asymmetric facet,
and let $G$ be any facet which is disjoint from $F$.
By Lemma~\ref{le:symasym},
$G$ is also asymmetric.
Hence  
Proposition~\ref{prop:powerful} implies that $F$ and $G$ are powerful.
Next, 
Lemma~\ref{le:power} implies that $\Delta$ is
combinatorially equivalent to $\Delta_1 \times F$.
\end{rmk}

In the proof of our main proposition we need the following
well known technical lemma. We include a proof for the 
convenience of the reader.

\begin{lemma}\label{le:posspan}
Fix vectors $\eta_1,\ldots,\eta_N$ in a $n$-dimensional vector space
$\ft$. Assume that 
their positive span contains $\ft$
and that no subset of  $n$ vectors
is linearly dependent.
After possibly renumbering,
the positive span of $\eta_1,\ldots,\eta_{n+1}$ contains $\ft$.
\end{lemma}

\begin{proof}
After possibly renumbering, there exists $k$ so that
all of $\ft$ lies in 
the positive span of $\eta_1,\ldots,\eta_k$,
but 
not in the positive span of any proper subset.
Assume that  $k > n+1$.
Write 
$0 = \sum_{i=1}^k a_i \eta_i$, where $a_i >  0$ for all $i$.
There exists a linear relation
$\sum_{i=1}^{n+1} d_i \eta_i = 0$;
let $d_i = 0$ for all $n+1 < i \leq k$.
Since $k > n+1$, by adding an appropriate multiple of 
$(d_1,\ldots,d_k)$ to $(a_1,\ldots,a_k)$,
we can find non-negative numbers $b_1,\ldots,b_k$  -- at least one 
of which is zero and at least one of which  is not --
so that $\sum_{i=1}^k b_i \eta_i = 0$.
After renumbering again, this implies
that $\sum_{i=1}^\ell b_i \eta_i = 0$, where 
$\ell < k$ and $b_i > 0$ for all $i$.
Since no set of $n$ vectors is linearly dependent,
this implies that 
every $\alpha \in \ft$ can be written as a
(not necessarily positive) linear combination 
$\alpha = \sum_{i=1}^\ell c_i \eta_i$.
Hence, for large enough $t$, $\alpha$ can be written as a positive
sum
$\alpha = \sum_{i =1}^\ell  (c_i + t b_i)$.
This contradicts the claim that $k$ is minimal.
\end{proof}

\begin{prop}  \labell{prop:powerful2}
Let $\Delta$ be a $n$-dimensional simple polytope; 
assume that every facet is powerful.
If $n \leq 4$, then $\Delta$ is
combinatorially equivalent to the product of simplices.
\end{prop}

\begin{proof}{} 
The polytope $\De$ has at most $2n$ facets. 
To see this, pick a vertex $v$ of $\De$.  Since
$\De$ is simple, $v$ lies on $n$ facets and $n$ edges.
Since every facet of $\De$ is powerful, 
there is an edge from $v$ to every facet that does not contain $v$. 

By Lemma~\ref{le:power}, 
we may assume that every facet of $\Delta$ is pervasive.
By Remark~\ref{rmk:3dimperv}, if 
$n = 3$ 
this implies that 
$\Delta$ is combinatorially equivalent to a  simplex;
the same claim is obvious if
$n < 3$  or if  $n = 4$
 and  $\Delta$ has 
five facets.
So assume that $\Delta$ is $4$-dimensional and  has
 $N$ facets,  where $6\le N \le 8$.

Write $\Delta = \bigcap_{i=1}^N 
\{ x \in \ft^* \mid  \langle \eta_i, x \rangle \leq \kappa_i \}$. 
Since our statement is purely combinatorial
we may slightly perturb the  facets so that the polytope is 
generic.  In particular, no four outward conormals are linearly  
dependent;
hence, after possibly renumbering,
the lemma above implies that
the positive span of $\eta_1,\ldots,\eta_5$
contains all of $\ft$. 
Now define
$$
P: = \bigcap_{i=1}^6
\{ x \in \ft^* \mid  \langle \eta_i, x \rangle \leq \kappa_i \}
\quad \mbox{and} \quad 
H_\ell^+ :=  \{ x \in \ft^* \mid  \langle \eta_\ell, x \rangle \leq \kappa_\ell \}
\quad \mbox{for} \quad 7 \leq \ell \leq N.$$
Clearly, $$
\Delta = P \cap \bigcap_{\ell=7}^N H_\ell^+.
$$
By construction, $P$ is bounded; 
by the genericity assumption, $P$ is simple, as is
$P \cap H_\ell^+$ for any $\ell$.
Since every  facet of $\Delta$ is pervasive,
and every facet of $P$ contains a facet of $\Delta$, 
the facets of $P$ are also pervasive.
Moreover, it is well known that any $n$-dimensional 
polytope with $n+2$ facets, such as $P$, 
is 
combinatorially equivalent to
a product of simplices;
for a proof in the current setting see ~\cite[Prop~1.1.1]{Tim}. 
Thus $P$ is combinatorially equivalent
either to  $\Delta_2\times \Delta_2$ or to $\De_1\times \De_3$.
Finally, some of the facets  of the latter polytope are nonpervasive, and so
$P$ is combinatorially equivalent
to  $\Delta_2\times \Delta_2$.

This solves the case $N = 6$, so assume that $N > 6$.
Label the facets of $P$ as $F_k$ and $F'_k$ for $0 \leq k \leq 2$,
so that $F_0 \cap F_1 \cap F_2 = F'_0 \cap F'_1 \cap F'_2 = \emptyset$.
Let $V$ be the set of all vertices of $P$;
label them 
$$
v_{ij} = \bigcap_{m \neq i}  F_m \cap\bigcap_{n \neq j} F'_n.
$$
We will think of the vertices as lying on a $3 \times 3$ grid,
where the first subscript determines the row and the second
determines the column.

Now, fix 
$\ell$ such that 
$7 \leq \ell \leq N$.  Let  $V_\ell^+ = V\cap H_\ell^+$, 
and let $V^-_\ell \subset V$ be the complement of  
$V_\ell^+$.  
Since the facets of $\Delta$ are pervasive, the
facets $F_0 \cap \Delta$ and $F_1 \cap \Delta$ must intersect.
A fortiori, $F_0 \cap F_1 \cap H_\ell^+ \neq \emptyset$.  Since $H_\ell^+$
is a half space, this implies that (at least) one of the
three vertices $v_{20}, v_{21}$ or $v_{22}$
of $F_0\cap F_1$
lies in $V_\ell^+$.  
A similar argument implies that $V_\ell^+$ contains at least one
vertex in each row and at least one in each column.

Moreover,
suppose that $v_{00}$ and $v_{11}$ are in $V^+_\ell.$
Consider the
$2$-dimensional face $F_2 \cap F_2'$; it  is a quadrilateral with vertices
$v_{00}, v_{01}, v_{11},$ and $v_{10}$.   The intersection of
this face with the half space $H_\ell^+$ cannot contain just the
two opposite vertices $v_{00}$ and $v_{11}$; it must also contain
(at least) one of the vertices $v_{10}$ or $v_{01}$.
More generally, if $V^+_\ell$ contains
the vertices $v_{ij}$ and $v_{mn}$ 
where $i \neq m$ and $j \neq n$, then 
it also contains $v_{in}$ or $v_{mj}$.
A brief analysis shows that this,
together with the paragraph above, implies that
$V_\ell^+$ contains at least five vertices.

On the other hand, the facets  $F_i \cap \Delta$ 
and $F_i' \cap \Delta$  must also
intersect the new facet $F_\ell$ for each $i$.
Therefore, 
the facets $F_i$ and $F_i'$
cannot lie entirely in $V_\ell^+$.
Hence at
least one vertex of  $F_i$ and at least one vertex
of $F_i'$ lies in $V^-_\ell$ for each $i$.  
It is straightforward to check
that this is only possible if
$V^-_\ell$ contains
two vertices $v_{ij}$ and $v_{mn}$,
where $i \neq m$ and $j \neq n$. 
As in the previous paragraph, this implies that 
$V^-_\ell$ also contains either $v_{in}$ or $v_{mj}$.
In either case,
$V_\ell^-$  contains two vertices which lie in the same column,
and two vertices which lie in the same row.

Now suppose that $v_{00}$ is a vertex of $\Delta$.
Every edge of $\Delta$ which meets $v_{00}$ is of the
form $e \cap \Delta$, where $e$ is an edge of $P$ which
meets $v_{00}$, that is,
where
$e$ is  $ \overline{v_{00} v_{10}}$, $\overline{v_{00} v_{20}}$,
$\overline{v_{00} v_{01}}$, or $\overline{v_{00} v_{02}}$.
Moreover, consider
-- for example -- 
 the edge  $e = \overline{v_{00} v_{10}}$.
If $e \cap \Delta$ connects 
$v_{00}$  to  
$F_0 \cap \Delta$,
then
$v_{10}$ must lie in $\Delta$;
otherwise  $e \cap \Delta$ will not intersect 
$F_0 \cap \Delta$.
On the other hand, if $e \cap \Delta$ connects $v_{00}$
to the new facet $F_\ell$
of $\Delta$,
then $v_{10}$ must lie in $V_\ell^-$.
Since $F_0 \cap \Delta$ is powerful, there must be an edge from
$v_{00}$ to $F_0 \cap \Delta$. 
This edge must be of the form $e \cap \Delta$, where $e$ is
an edge in $P$ from $v_{00}$ to $F_0$, that is,
$\overline{v_{00} v_{10}}$ or $\overline{v_{00} v_{20}}$.
Hence, (at least) one of the vertices $v_{10}$ and
$v_{20}$ lies in $\Delta$.  
On the other hand, the new facet $F_\ell$ is also powerful.
Hence, $V^-_\ell$ contains (at least) one of the vertices $v_{10}, v_{20},
v_{02}$, or $v_{01}$.
A similar argument
shows that  $\Delta$ cannot contain exactly
one vertex in any row or column, and that if
$\Delta$  contains $v_{ij}$ 
then $V^-_\ell$ contains some vertex which either lies
in the same row or in the same column.

We now specialize, and assume that $N = 7$.
In this case, $v_{ij}$  lies in $\Delta$ exactly if
$v_{ij} \in V_7^+$.
Therefore, since 
$V^+_7$ contains at least one vertex in
each column, the paragraph above implies
that it contains at least two vertices in
each column.
But this contradicts
the fact that $V_7^-$ contains two vertices in the same column.

So we may assume that $N = 8$.
In this case, $v_{ij}$  lies in $\Delta$ exactly if
$v_{ij} \in V_7^+ \cap V_8^+$.
Hence,  $V_7^+ \cap V_8^+$
cannot contain exactly one vertex in any row or column.
On the other hand, $V_7^-$ contains two vertices in the same
column and also contains two vertices in the same row, 
so $V^+_7 \cap V^+_8$ cannot contain two
or more
 vertices
in every column or in every row.
Together, these imply that $V^+_7 \cap V^+_8$ must not contain
any vertices in at least one row and at least one column.
Assume that these are the last row and column.
Since $V^+_7$ and $V^+_8$ each contain at least five
vertices, $V^+_7 \cap V^+_8$ is not empty.
This implies that,
$V^+_7 \cap V^+_8 = \{v_{00},v_{01},v_{10},v_{11}\}$.
But then $V^-_7$ must contain at least one vertex in
the same row or column as 
each of
$v_{00},v_{01},v_{10}$ and $v_{11}$.
After possibly swapping the columns with the rows,
this implies that $V_7^-$ contains $v_{02}$ and $v_{12}$.
Since $V_7^+$ contains at least one vertex in every column,
this implies that $V_7^+$ contains $v_{22}$.
A similar argument implies that $V_8^+$ contains $v_{22}$.
This contradicts the fact that $v_{22}$ does not
lie in $V^+_7 \cap V^+_8$.
\end{proof}

\begin{rmk}\labell{rmk:app}\rm It is natural to wonder 
if there is an analog of this result in higher dimensions. 
Chen~\cite{Chen} recently found a 
$7$-dimensional
polytope with all facets pervasive and powerful
that is not 
combinatorially equivalent to
a product.  However,  
since his example is not smooth,
the following question is still open:
\begin{quote}{\it 
Let $\Delta$ be a smooth 
polytope;  assume that every facet
is powerful.  Is $\Delta$ combinatorially equivalent to the
product of simplices? If not,
is there a structure theorem for such polytopes?}
\end{quote}
Note that, by Lemma~\ref{le:power}, it is enough to consider the
case that  every facet of $\Delta$ is pervasive.
Moreover,  the proof above shows that if 
$\Delta$ is $n$-dimensional, then it has at most $2n$ facets. 
\end{rmk}

Combined,  Propositions~\ref{prop:powerful} and \ref{prop:powerful2}
have the following important corollary.

\begin{corollary}\labell{cor:allasym}
Let $H \in \ft$ be a mass linear function on a simple 
$4$-dimensional polytope $\Delta \subset \ft^*$
with no symmetric facets.
Then  $\Delta$ is combinatorially equivalent to 
a product of simplices.
\end{corollary}

This allows us to complete the first step in
classification of mass linear function on smooth
$4$-dimensional polytopes.

\begin{theorem}
Let $H \in \ft$ be a mass linear function on a smooth
$4$-dimensional polytope $\Delta \subset \ft^*$. 
Then there exists an inessential function $H' \in \ft$
so that the mass linear function $\Tilde H = H - H'$
has the following property: at least one facet of
$\Delta$ is $\Tilde H$-symmetric.
\end{theorem}

\begin{proof}
Assume that $\Delta$ has no symmetric facets.
By Corollary~\ref{cor:allasym}, 
$\Delta$ is combinatorially equivalent to a product of simplices.
Moreover, by Proposition~\ref{prop:flat} we may assume that
every facet if pervasive;
hence, $\Delta$ is combinatorially equivalent  to $\Delta_2 \times \Delta_2$.
By Lemma~\ref{prodsimp}, this implies that $\Delta$ is a $\Delta_2$ bundle
over $\Delta_2$.  The result now follows immediately
from Proposition~\ref{prop:bund}.
\end{proof}
\end{appendix}


\begin{thebibliography}{cccc}

\bibitem{Aud} M. Audin, {\it Torus actions on symplectic manifolds}, 
Progress in Mathematics vol 93, 2nd revised edition, Birkh\"auser (2005).

\bibitem{CDG} D. Calderbank, L. David and P Gauduchon,
The Guillemin formula and K\"ahler metrics on 
toric symplectic manifolds,
{\it J. Symp. Geom.} {\bf 1} (2003), 767--84.

\bibitem{Chen} J.-W. Chen,  Neighborly properties of simple convex polytopes, Ph. D. thesis, Stony Brook (2007).

\bibitem{Cox} D. Cox, The 
homogeneous
coordinate ring of a toric variety, 
{\it J. Algebraic Geom.} {\bf 4} (1995), 17--50.

\bibitem{CK} D. Cox and S. Katz, {\it Mirror Symmetry and Algebraic Geometry}, Mathematical Surveys and Monographs vol 68, Amer. Math Soc. (1999).


\bibitem{Del} T. Delzant, 
      {\it Hamiltoniens p\'eriodique et images convexes de l'application
       moment}, Bull. Soc. Math. France 116 (1998) 315-339.

\bibitem{De}  M. Demazure, Sous-groupes alg\'ebriques de rang maximum du groupe de Cremona, {\it Ann. Ec. Norm. Sup.} {\bf 3} (1970), 507--88.


\bibitem{EP} M. Entov and L. Polterovich,
    Rigid subsets of symplectic manifolds,
      arXiv:math/0704.0105

%\bibitem{Fu}  W. Fulton, {\it Introduction to Toric Varieties}, 
%Annals of Math. Studies {\bf 131}, PUP (1993).


\bibitem{Gr} B. Gr\"unbaum, {\it Convex Polytopes}, Graduate Texts in Mathematics 221, Second edition, Springer (2003). 

\bibitem{JK} T. Januszkiewicz and J. K\c edra, Characteristic classes of 
       smooth fibrations,  SG/0209288.

\bibitem{KM}
J. K\c edra and D. McDuff, Homotopy properties of Hamiltonian group actions, {\it Geometry and Topology}, {\bf 9} (2005) 121--162.

\bibitem{MT1} D. McDuff and S. Tolman, Topological properties of Hamiltonian circle actions, SG/0404338,      {\it International Mathematics Research Papers}.
      {\bf vol. 2006},  Article ID 72826, 1--77.
      
\bibitem{NILL} B. Nill, Complete toric varieties with reductive         
       automorphism group, AG/0407491, {\it Math. Z.} {\bf 252} (2006), 767--86.

\bibitem{Rez} A. G. Reznikov, Characteristic classes in symplectic 
topology, {\it Selecta Math,} {\bf 3} (1997), 601--642.

\bibitem{Sei} P. Seidel, $\pi_1$  of symplectic automorphism groups
          and invertibles in quantum cohomology rings, {\it Geom. and Funct. 
         Anal.} {\bf 7} (1997), 1046 -1095.
         
\bibitem{Tim} V. Timorin, An analogue of the 
Hodge--Riemann relations for simple convex polytopes, 
{\it Russ. Math Surveys} {\bf 54:2} (1999), 381--426.

\bibitem{Vi} A. Vi\~na, Hamiltonian diffeomorphisms of toric manifolds, SG/0506176, {\it J. Geom. Phys.} {\bf 57} (2007), 943--65.

\bibitem{WaZ}  X. J. Wang and X. Zhu, K\"ahler Ricci solitons on toric manifolds
with positive first Chern class,  {\it Advances in Mathematics} {\bf 188}
 (2004) 87--103.


\end{thebibliography}
\end{document}